\documentclass[10pt]{article}

\usepackage{amsmath}
\usepackage{amssymb}
\usepackage{caption}
\usepackage{graphicx}
\usepackage{subfig}
\usepackage{authblk}
\usepackage{verbatim}
\usepackage{epstopdf}
\usepackage{comment}
\usepackage{color}

\usepackage[hypertexnames=false,colorlinks=true,linkcolor=blue,citecolor=blue]{
hyperref}
\usepackage[numbers,comma,square,sort&compress]{natbib}
\usepackage[a4paper,text={6.5in,10in},centering]{geometry}
\usepackage{verbatim}
\usepackage[mathlines]{lineno}

\graphicspath{{Figures/}}

\makeatletter\@addtoreset{equation}{section}\makeatother
\renewcommand{\theequation}{\arabic{section}.\arabic{equation}}

\newtheorem{theorem}{Theorem}[section]

 {\begin{trivlist}\item[]\textbf{Proof#1 }}%
 {\hspace*{\fill}$\rule{0.3\baselineskip}{0.35\baselineskip}$\end{trivlist}}

\def\phi{\varphi}
\def\R{\mathbb{R}}

\def\sgn{\operatorname{sgn}}
\def\({\left(}
\def\){\right)}
\def\eb{\epsilon}

\renewcommand{\Re}{{\rm Re}}
\renewcommand{\Im}{{\rm Im}}
\def\ri{{\rm i}}
\def\d{{\rm d}}
\def\Phihat{\widehat{\Phi}}
\def\Psihat{\widehat{\Psi}}
\def\hhat{\widehat{h}}
\def\htilde{\widetilde{h}}
\def\that{\widehat{t}}
\def\phat{\widehat{p}}
\def\qhat{\widehat{q}}
\def\shat{\widehat{s}}
\def\Ghat{\widehat{G}}
\def\what{\widehat{w}}
\def\gbar{\overline{g}}
\def\rbar{\overline{r}}
\def\qbar{\overline{q}}
\def\eq{{\rm eq}}

\def\ctext{\color{black}}

\begin{document}

\title{The dynamics of interacting multi-pulses in the one-dimensional quintic complex Ginzburg-Landau equation}
\author{T. Rossides}
\author{D. J. B. Lloyd}
\author{S. Zelik}
\author{M. R. Turner\footnote{Corresponding author: m.turner@surrey.ac.uk}}
\affil{\small School of Mathematics and Physics, University of Surrey, Guildford, GU2
7XH, UK}
\date{\today}
\maketitle

\begin{abstract}

We formulate an effective numerical scheme that can readily, and accurately, calculate the dynamics of weakly interacting multi-pulse solutions of the quintic complex Ginzburg-Landau equation (QCGLE) in one space dimension. The scheme is based on a global centre-manifold reduction where one considers the solution of the QCGLE as the composition of individual pulses plus a remainder function, which is orthogonal to the adjoint eigenfunctions of the linearised operator about a single pulse. This centre-manifold projection overcomes the difficulties of other, more orthodox, numerical schemes, by yielding a fast-slow system describing `slow' ordinary differential equations for the locations and phases of the individual pulses, and a `fast' partial differential equation for the remainder function. With small parameter ${\ctext \epsilon=e^{-\lambda_r d_0}}$ where $\lambda_r$ is a constant and $d_0>0$ is the minimal pulse separation distance, we write the fast-slow system in terms of first-order and
second-order correction terms only, a formulation which is solved more efficiently than the full
system. This fast-slow system is integrated numerically using adaptive time-stepping. Results are presented here for two- and three-pulse interactions. For the two-pulse problem, cells of periodic behaviour, separated by an infinite set of heteroclinic orbits, are shown to `split' under perturbation creating complex spiral behaviour. For the case of three pulse interaction a range of dynamics, including chaotic pulse interaction, are found. While results are presented for pulse interaction in the QCGLE, the numerical scheme can also be applied to a wider class of parabolic PDEs.

\end{abstract}

\section{Introduction}\label{sec:Introduction}
The study of single and multi-pulse states has provided many insights over the years for a range of physical systems \cite{dysthe1979,remoissenet1999,rs1997,dauxois2006}. When multiple pulses/solitons are placed close to one another they generally interact and may generate complicated, or even chaotic, space-time dynamics; see \cite{ZelMiel09,Turaev07,Turaev09,Turaev12} and references therein. However, this interaction is usually very weak since its appearance is due to nonlinear effects and the interaction of the pulse tails, which in most cases are exponentially decaying. Since the rate of interaction is exponentially small with respect to the distance between pulses, direct numerical integration methods are usually inefficient for computing the resulting interaction and their delicate dynamical properties as the computations can take very long times. Thus, more efficient numerical techniques and analytical methods are required in order to identify these dynamics.
\par
One of the most promising analytical methods for studying multi-pulses in {\it dissipative} systems is the so-called {\it centre-manifold reduction} developed in \cite{Bjorn02,ei2002,ZelMiel09}. Indeed, let us consider the dissipative system of the form
\begin{equation}\label{0.init}
\partial_t u=\mathcal{A}u+f(u),\ ~~~~~ \ f(0)=f'(0)=0,
\end{equation}
where $u=(u^1,\cdots,u^n)(t,x)$ is an unknown vector function, $\mathcal{A}$ is a  uniformly elliptic differential operator in $x\in\R^m$, $t\in[0,\infty)$ is time and $f$ is a given nonlinear function. If \eqref{0.init} possess a steady pulse equilibrium solution $V(x)$, then, due to the symmetries of the system, it simultaneously possess a whole manifold, $V_\xi(x)$, of such solutions parametrized by $\xi$, which normally lives in the Lie group of symmetries of the considered system. In particular, if~\eqref{0.init} is spatially homogeneous, $\xi=r\in\R^n$, where the parameter $r$ corresponds to spatial shifts of the pulse. It is then natural to seek a desired multi-pulse solution $u(t,x)$ in the form
\begin{equation}\label{0.mult}
u(t,x)=\sum_{i=1}^n V_{\xi_i(t)}(x)+w(t,x),
\end{equation}
where the reminder function $w(t,x)$ is assumed to be small in comparison to the pulses themselves. This construction suggests that the dynamics of the multi-pulse system are determined, up to a small corrector $w$, by the slow evolution of the internal parameters $\vec\xi(t):=(\xi_1(t),\cdots,\xi_n(t))$ belonging to the Lie group of symmetries of the considered system. Such a solution structure has also been used successfully in the study of excitation kinks in parabolic differential equations \cite{pauthier2021}.
\par
As shown in \cite{ZelMiel09}, the solution $u(t,x)$ can indeed be found in the form of \eqref{0.mult} under some natural assumptions on the dissipative system considered, provided the {\ctext minimal} separation distance
\begin{equation}\label{0.dist}
{\ctext d_0:=\min_t\left[\min_{i\ne j}|r_i(t)-r_j(t)|\right]},
\end{equation}
between pulses is large. Moreover, in this case the initial dissipative system can be re-written in the neighbourhood of a multi-pulse structure in terms of a fast-slow system with respect to the new dependent variables $\vec\xi$ and $w$ which has the following general form
\begin{equation}\label{0.slow-fast}
\begin{cases}
\partial_t w-\mathcal L(\vec\xi(t))w=-\epsilon h(\epsilon, \vec\xi(t),w)+\epsilon\Phi(\eb,\vec\xi(t))+G(\epsilon,\vec\xi(t),w),\ \ G(\epsilon,\xi,0)=G'_w(\epsilon,\xi,0)=0,\\
\frac \d{\d t}\vec\xi(t)=\epsilon H(\epsilon,\vec \xi(t),w),
\end{cases}
\end{equation}
where ${\ctext \epsilon=O(e^{-\lambda_r d_0})}$, with $\lambda_r>0$, is a small parameter, $\mathcal L(\vec\xi)$ is an invertible linear operator for all admissible $\vec\xi$ and $\Phi$, $G$, $H$ and $h$ are nonlinear functions which can be explicitly found.
\par
The structure of the fast-slow system \eqref{0.slow-fast} allows us to apply the standard centre-manifold reduction technique and verify the existence of an invariant manifold of the form $w(t)=W(\epsilon,\vec\xi(t))$ in the properly chosen phase space of~\eqref{0.slow-fast}. The reduction remains valid while the separation distance $d(t)>d_0$ between pulses remains large. Thus, the dynamics of the multi-pulse interaction will be governed by a system of ODEs
\begin{equation}\label{0.fullODE}
\frac \d{\d t}\vec\xi(t)=\epsilon H(\epsilon,\vec\xi(t),W(\epsilon,\vec\xi(t)),
\end{equation}
on the Lie group of symmetries of the initial dissipative system. It is also clear that $W=O(\epsilon)$, so the leading order approximation to \eqref{0.fullODE} has the form
\begin{equation}\label{0.redODE}
\frac \d{\d t}\vec\xi(t)=\epsilon H(0,\vec\xi(t),0).
\end{equation}
Rather surprisingly, the nonlinearity in the leading order approximation can be found analytically (up to a few constants which usually have to be calculated numerically) in many physically relevant cases. Thus, this approach allows us to approximate the multi-pulse dynamics by the explicitly written system of ODEs which can be further studied both analytically or numerically, see \cite{ZelMiel09} for more details.
\par
It was suggested in \cite{rossides2014} to use the fast-slow system \eqref{0.slow-fast} for numerical simulations, rather than directly integrating the original PDE, in cases where accuracy beyond the leading order approximation \eqref{0.redODE} is required. The advantage of this approach was demonstrated for the one-dimensional, real Ginzburg-Landau equation \cite{rossides2014} where the associated multi-kink interactions were computed. However, the method suggested in \cite{rossides2014} is based upon the straightforward discretization of the first equation in the fast-slow system~\eqref{0.slow-fast} (which is an evolutionary PDE) and did not utilize the fast-slow structure of this system and the associated centre-manifold reduction. Thus, the numerical scheme may still be slow when the separation distance of the pulses is large.
\par
In the present paper, we take the development of the scheme suggested in \cite{rossides2014} to the next level by taking advantage of the fast-slow structure of \eqref{0.slow-fast} and approximating the slow-manifold, {\ctext and we are expanding the manifold in terms of $\epsilon$ which does not affect the temporal dynamics on the manifold.} Namely, instead of solving the first PDE in \eqref{0.slow-fast} fully, we asymptotically determine the leading order terms in the expansion of the centre manifold $W$ with respect to the small parameter $\epsilon$, and solve this equation numerically. Indeed, the leading order term has the form
\begin{equation}\label{0.lead}
W_1(\vec\xi)=\epsilon[\mathcal{L}(\vec\xi)]^{-1}(\Phi(0,\vec\xi)-h(0,\vec\xi,0)),
\end{equation}
and thus its computation is reduced from a fully time-dependent PDE, to inverting the linear elliptic operator $\mathcal{L}(\vec\xi)$. Inserting this approximation into \eqref{0.fullODE}, we end up with the ODE system
\begin{equation}\label{0.3-ODE}
\frac \d{\d t}\vec\xi(t)=\epsilon H(\epsilon,\vec\xi(t),W_1(\vec\xi(t)),
\end{equation}
which captures the pulse interaction dynamics up to correction terms of $O(\epsilon^3)$. Note that while we terminate our approximation by including the $O(\epsilon^2)$ terms, we can compute higher-order approximate systems as desired, each with smaller and smaller error terms, leading to a hierarchy of numerical schemes.


\par
In this paper we apply the centre-manifold reduction to the model example of the 1D quintic complex Ginzburg-Landau equation (QCGLE), given by
\begin{equation}\label{e:CGLE}
 u_t = \alpha u_{xx} + \beta u + \gamma u|u|^2 + \delta u|u|^4 := \alpha u_{xx}+\beta u+u\theta(|u|^2),
\end{equation}
where $u(t,x)\in\mathbb{C}$, $\alpha,~\beta,~\gamma,~\delta$ are arbitrary complex parameters with $\Re(\alpha)>0$ and $\theta(z):=\gamma z+\delta z^2$ is the nonlinear component. Here $x\in\R$ is the single space dimension and $t\in[0,\infty)$ is time. This equation has been studied in nonlinear optics, and is used to describe phenomena related to pulse formation, such as mode locking in lasers \cite{Akhmediev972,Ding2011,Moores93}, light propagation in nonlinear fibers \cite{Akhmediev972} and transverse pattern formation in nonlinear optic systems \cite{Mandel04}.
\par
It is well-known that the QCGLE possesses stable, reflectional-symmetric pulses in a large region of parameter space of the form
\begin{equation}\label{e:pulseprof}
u(t,x)= V(x),\ ~~{\rm with~the~property}~~\ V(-x)=V(x).
\end{equation}
Moreover, equation \eqref{e:CGLE} possesses not only the shift symmetry $x\to x-r,r\in\mathbb{R}$ but also a phase symmetry $u\to e^{\ri g}u,~g\in\mathbb S_1$. Thus we have $\xi=(r,g)$, which together with the initial pulse \eqref{e:pulseprof}, gives the whole two-dimensional manifold of pulse solutions
\begin{equation}\label{e:pulseman}
\mathcal P:=\{V_\xi(x):=e^{\ri g}V(x-r),\ \ \xi\in\mathbb \R\times \mathbb S_1\},
\end{equation}
see for example~\cite{Deissler94,Turaev07}. Since these stable pulses can be arbitrarily shifted in their location and phase, one expects that multi-pulse solutions will affect both each others location and phase. Interaction of
well-separated pulses can produce spatial oscillations which are responsible for the formation
of bound states of the dissipative pulses; as seen in experimental studies such as \cite{Tang2001,Grelu2003}.
In order to demonstrate the type of dynamics we expect from multi-pulse systems we examine the case of two weakly-interacting pulses in more detail than was highlighted in \cite{ZelMiel09}. Following the general strategy described above, we seek a two-pulse solution of the form
\begin{equation}\label{e:2pans}
u(t,x) = V_{\xi_1(t)}(x) + V_{\xi_2(t)}(x) + w(t,x),
\end{equation}
where $\xi_k=(r_k,g_k)\in\mathbb \R\times \mathbb S_1$ and $w(t,x)$ is a small reminder function. Then, under the generic assumption that the initial one-pulse manifold $\mathcal P$ is normally hyperbolic, it is possible to justify the centre-manifold reduction, and verify that the two-pulse dynamics is determined by a system of four ODEs for the variables $\xi_1(t) = (r_1(t),g_1(t))$ and $\xi_2(t)=(r_2(t),g_2(t))$. Moreover, due to the symmetries of the this system, the problem can be reduced to a planar system of two ODEs with respect to the difference variables $\bar g=g_2-g_1$ and $\bar r=r_2-r_1$ (we assume here that $r_2-r_1>0$), see \cite{ZelMiel09} for the full details.
\par
The leading order planar system approximation to \eqref{0.redODE} can be expressed as the pair of ODEs
\begin{equation} \label{e:ode2}
\begin{cases}
\frac \d{\d t}\bar r =  Je^{-\lambda_r\bar r}\textrm{sin}(-\lambda_i \bar r + \kappa _1) \textrm{cos}(\bar g), \\
\frac \d{\d t}\bar g = -Ke^{-\lambda_r\bar r}\textrm{cos}(-\lambda_i \bar r + \kappa _2) \textrm{sin}(\bar g),
\end{cases}
\end{equation}
where the components of $\lambda=\lambda_r+\ri\lambda_i$ are determined by the dispersion relation
\begin{equation}\label{0.disp}
\alpha\lambda^2+\beta=0.
\end{equation}
 The other real parameters $\kappa_1$, $\kappa_2$, $J$ and $K$ can be determined analytically as long as the asymptotic form of the initial pulse tails, and the associated adjoint eigenfunctions tails, are analytically known; see \cite{ZelMiel09,Akhmediev97}. Typically these parameters are found numerically.

\begin{figure}[ht!]
     \begin{center}
(a)\includegraphics[width=2.0in, height=2.0in]{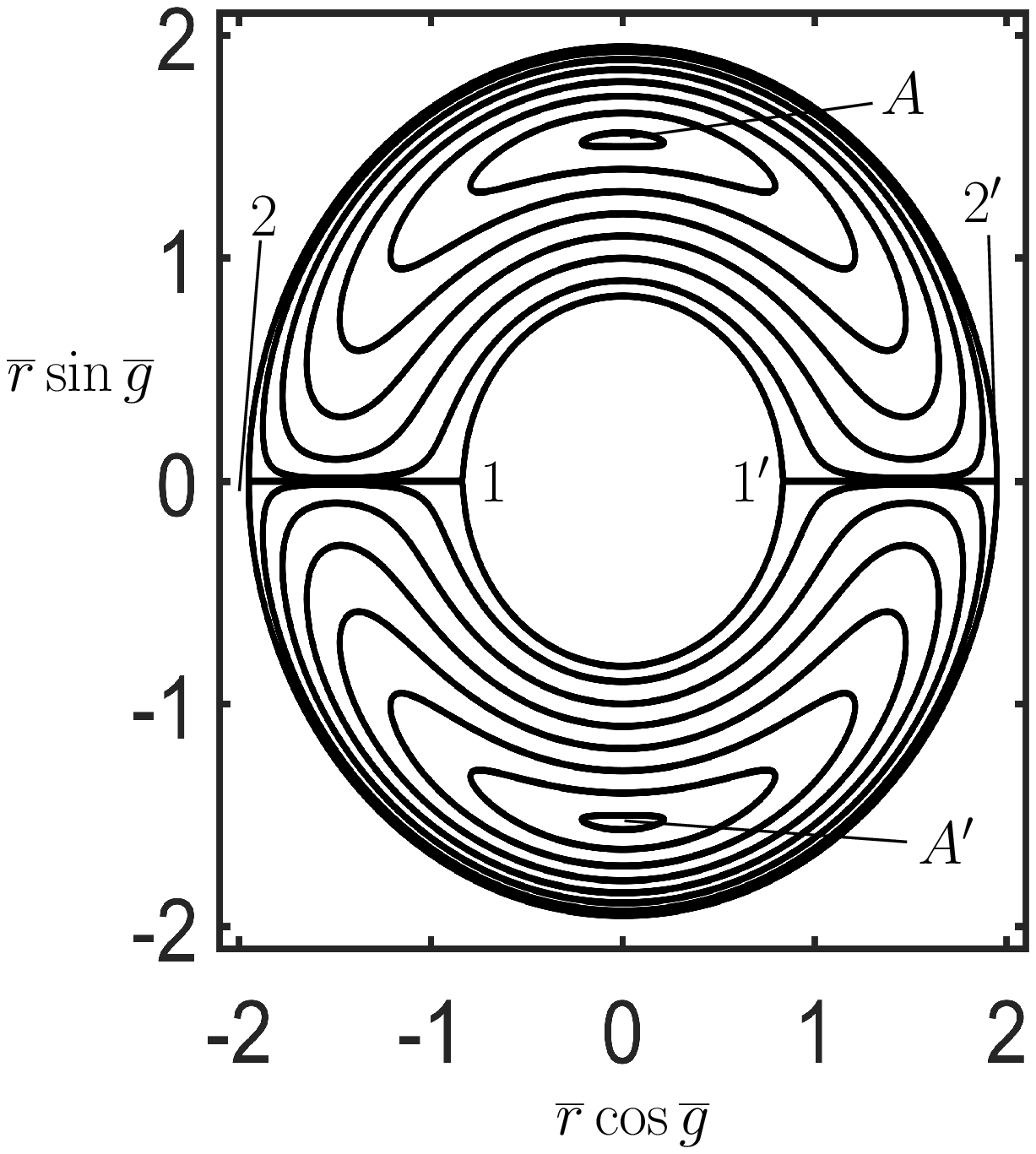}
(b)\includegraphics[width=2.5in]{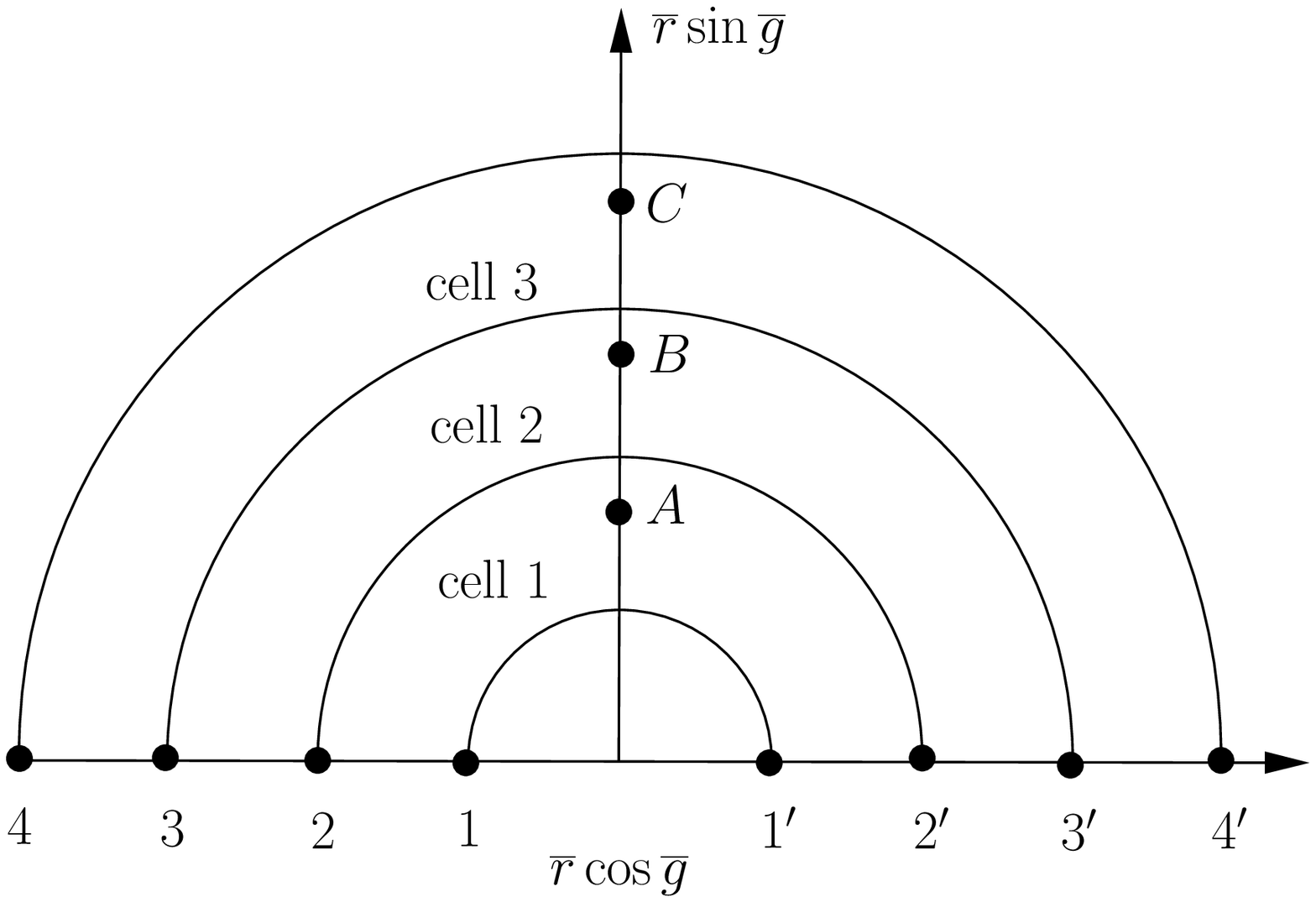}
    \end{center}
    \caption{(a) The Hamiltonian phase portrait for cell 1 dynamics in the reduced projected system \eqref{e:ode2}, describing the evolution of separation distance $\bar r =r_2-r_1,$ and phase difference $\bar g=g_2-g_1$. Note that all trajectories are closed and the equilibrium points $A$ $\&$  $A'$ with
$\rbar\cos(\gbar)=0$ are identified as centres, whereas points $1,2,1',2'$ with $\rbar\sin(\gbar)=0$ are saddles. (b) A schematic phase portrait of the upper half-plane for the reduced projected system, split into cells based on the separating heteroclinic orbits which connect the equilibrium saddle points $1,2,1',2'$ on the axis $\rbar\sin(\gbar)=0$.}
   \label{fig:Sergeyppred}
\end{figure}

Depending on the difference of the two phase parameters $\kappa_1-\kappa_2$, two principally different types of phase portrait for the planar system of ODEs \eqref{e:ode2} can arise. The first one has a gradient-like behaviour and the second one is Hamiltonian (which is reversible). In the first case, the system is structurally stable, so the leading order approximation determines completely the type of the dynamics observed, so we concentrate below on the more interesting Hamiltonian case only. An example phase portrait for this case is shown in figure~\ref{fig:Sergeyppred}(a).

For system parameters which give rise to the Hamiltonian case, the planar system \eqref{e:ode2} for $\bar r$ and $\bar g$ has two types of equilibria. The first type corresponds to in-phase ($\bar g=\bar g_{\eq}=0$) and anti-phase ($|\bar g_{\eq}|=\pi$) pulse pairs and they are all saddles found at $\rbar=\rbar_{\eq}=\pm(n\pi+\kappa_1)/\lambda_i$ for $n=1,2,...$. The second type corresponds to pulse pairs with a phase shift $|\bar g_{\eq}|=\pi/2$ and are centres found at $\rbar_{\eq}=\pm((2n-1)\pi/2+\kappa_2)/\lambda_i$ for $n=1,2,...$. Due to the Hamiltonian structure of the
system (see, \cite{Turaev07}), heteroclinic orbits connect the saddles in semi-circular arcs and form closed cells filled by periodic orbits oscillating around the centre equilibria. We show schematically in figure \ref{fig:Sergeyppred}(b) how the heteroclinic orbits connect two saddles forming closed cells which we call cell 1, cell 2 etc. The two heteroclinic orbits HO1 and HO2 form the edges of cell 1 and defines what we consider as weakly-interacting pulses, while heteroclinic orbits HO2 and HO3 form the edge of cell 2 defining the region that we consider as very-weakly-interacting pulses. We will refer to these heteroclinic orbits as `cell boundaries' in the current article.
\par
There is, however, no reason to expect that the Hamiltonian/reversible structure of the leading order reduced system~\eqref{e:ode2} to persist in the full reduced system on the centre manifold, since the QCGLE is fundamentally a dissipative PDE which is neither Hamiltonian nor time-reversible. Hence, one expects that the periodic orbits, centres and cell boundaries shown in figure~\ref{fig:Sergeyppred} are not robust and will `split' when higher order terms are included in the reduction. Despite the cells no longer remaining closed in the full reduced system~\eqref{0.slow-fast}, the splitting will be small (at least of $O(\epsilon^2)$) and we will still be able to recognize different regions of the phase plane as cell 1 and cell 2 etc, based on the location of the cell boundary in the reduced system~\eqref{e:ode2}. Hence we will still use this terminology when discussing the full system.
\par
Numerical simulation results of two-pulse interactions in the full QCGLE presented by authors such as \cite{Turaev07,Akhmediev97}, and references therein, show that
the phase portrait in figure~\ref{fig:Sergeyppred} is indeed not robust and that the Hamiltonian structure disappears if higher order corrections are taken into the account. In particular, the centres  $A$, $B$ and $C$ as seen in figure \ref{fig:Sergeyppred}(b) become   {\it foci} for the phase trajectories, while the robust equilibria $1,~2,~3,$ etc. remain as saddle points. Also, if the interaction of the pulses is particularly strong, then an additional `collision area' appears in which the pulses come together and one is annihilated.

These papers have highlighted, that due to the very small level of interaction between the pulses, conventional numerical schemes used to date are computationally very expensive and the results obtained in these works were not sufficient to determine the dynamics of the phase portrait. In this paper, the centre-manifold reduction technique will produce a more efficient numerical scheme which will allow us to identify the phase-plane dynamics (at least for the first few cells) for a range of system parameters. In particular, we aim to answer the following open questions for the QCGLE (to the best of our knowledge):
\begin{enumerate}
\par
\item Do stable or unstable {\it limit cycles} exist in the cells, and if so, how are they generated?;
\par
\item How exactly do the heteroclinic orbits split, and what do the dynamics close to the cell boundaries look like?;
\par
\item How are these results altered by varying the parameters of the QCGLE?
\par
\end{enumerate}

The current paper is laid out as follows. In \S\ref{sec:Method} we present an overview of the novel centre-manifold reduction technique, along with the approximations we will use in order to investigate the system dynamics for multi-pulses in the QCGLE. We will also present an asymptotic analysis in terms of the small parameter $\epsilon$ to rigorously identify the governing problem consisting of leading order and second order terms, to be solved. In \S\ref{sec:Impl_Just} we present the numerical justification of the centre-manifold projection system, as well as a further approximation where we reduce the problem to an ODE system only, for the case of very-weakly-interacting pulses. Here we compare results of these two systems to the standard time-stepping scheme used in works such as \cite{Turaev07,Akhmediev97}. In \S\ref{sec:Results} we present numerical results for the two- and three-pulse interaction solutions of the QCGLE, where in particular we highlight the phase-plane dynamics for two-pulses for a selection of system parameters. Concluding remarks and discussion is given in \S \ref{sec:Conclusion}.

\section{The analytic structure of the centre-manifold}\label{sec:Method}

In this section, we briefly recall the derivation and basic properties of the fast-slow system \eqref{0.slow-fast} for the example of the QCGLE (see \cite{ZelMiel09} for more details), and also present a formal asymptotic calculation to determine the leading order form of $w(t,x)$ given in \eqref{0.lead}.

\subsection{Single-pulse manifold and normal hyperbolicity}
We shall assume that the parameters in the QCGLE are fixed in such way that equation \eqref{e:CGLE} possesses a symmetric pulse equilibrium $V(x)\in\mathbb{C}$, such that $V(-x)=V(x)$, which solves the ODE
\begin{equation}\label{e:CGLEeq}
 \alpha V_{xx} + \beta V + \gamma V|V|^2 + \delta V|V|^4 := \alpha V_{xx} +\beta V+
f(V)=0,
\end{equation}
where $f(V)=V\theta(|V|^2)$ and $\theta(z)=\gamma z+\delta z^2$. Then, due to the phase and shift symmetry, it automatically possesses a 2D manifold \eqref{e:pulseman} of pulses parameterized by $\xi=(g,r)\in\mathbb S_1\times\R$. The linearised operator near the pulse $V$ is given by
\begin{equation}\label{1.lin}
\mathcal L:=\alpha\partial_x^2+\beta+f'(V),~~~{\rm with}~~~ \ \mathcal L z:=\alpha\partial_x^2z+\beta z+\theta(|V|^2)z+|V|^2\theta'(|V|^2)z+V^2\theta'(|V|^2)\overline{z},
\end{equation}
where the overbar denotes the complex conjugate. Note, that while we write the linearised operator in the complex form above, due to its concise notation and ease of manipulation in \S\ref{sec:multimanifold}, in order to solve for both $V$ and $z$ we write the system in terms of real and imaginary parts. This effectively doubles the order of the system. We also introduce the linear operator 
\begin{equation}\label{1.linshift}
\mathcal L_{\xi}:=\alpha\partial_x^2+\beta+f'(V_\xi),
\end{equation}
with $\mathcal L_{\xi}:H^2(\mathbb{R},\mathbb{R}^2)\to L^2(\mathbb{R},\mathbb{R}^2)$, 
linearised near $V_\xi(x):=e^{\ri g}V(x-r)$. Due to the symmetries of the system, we can introduce the two eigenfunctions
\begin{equation}\label{1.phi}
\varphi^r(x):=\partial_rV_\xi(x)\big|_{\xi=0}=-V'(x)\ ~~~ \text{and}\ ~~~ \varphi^g(x):=\partial_gV_\xi(x)\big|_{\xi=0}=\ri V(x),
\end{equation}
of $\mathcal L$ with corresponding zero eigenvalues. The main assumption, which we assume to be satisfied throughout {\ctext the rest} of the paper, is that these are the only zero eigenvalues (no Jordan blocks are allowed at zero) and the rest of the spectrum of $\mathcal L$, $\sigma(\mathcal L)$, belongs to the half-plane with negative real part i.e., 
\begin{equation}\label{2.stable}
\sigma(\mathcal L)\cap\{\Re(z)>0\}=\varnothing,\ ~~~ \sigma(\mathcal L)\cap\{\Re(z)=0\}=\{0\},\ ~~~ \dim(\ker\mathcal L)=2.
\end{equation}
In particular, analyzing the essential spectrum of $\mathcal L$, it is easy to see that (\ref{2.stable}) can be satisfied only if ${\rm Re}(\beta)<0$. Furthermore, we introduce a natural inner product in $L^2(\mathbb C)$  via
\begin{equation}\label{1.inner}
\langle v,w\rangle:= \Re \left
\{\int_{-\infty}^\infty v(x)\overline{w(x)}\,\d x \right \}.
\end{equation}
Then, the adjoint operator $\mathcal L^*=\bar\alpha\partial_x^2+\bar\beta+[f'(V)]^*$ is given by
\begin{equation}
\mathcal L^*z:=\bar\alpha\partial_x^2z+\bar\beta z+\bar\theta(|V|^2)z+|V|^2\bar\theta'(|V|^2)z+V^2\theta'(|V|^2)\bar z.
\end{equation}
It is well-known that $\mathcal L:H^2(\mathbb{R},\mathbb{R}^2)\to L^2(\mathbb{R},\mathbb{R}^2)$ is Fredholm of index zero, so the dimension of the kernel of $\mathcal L^*:H^2(\mathbb{R},\mathbb{R}^2)\to L^2(\mathbb{R},\mathbb{R}^2)$ also equals two and the adjoint eigenfunctions $\psi^r(x)$ and $\psi^g(x)$ can be normalized such that
\begin{equation}\label{adjointcgl2}
\langle \phi^r,\psi^r\rangle=\langle\phi^g,\psi^g\rangle=1,\ ~~~~
\langle \phi^r,\psi^g\rangle=\langle\phi^g,\psi^r\rangle=0.
\end{equation}
Also, since the initial pulse $V(x)$ is even, we also may assume, without loss of generality, that the functions $\phi^g$ and $\psi^g$ are even and the functions $\phi^r$ and $\psi^r$ are odd. In the cases we consider, the essential spectrum of the linear operator is in the left half plane, hence, this restricts us to the case ${\rm Re}(\beta)<0$. Therefore the eigenfunctions $\phi$ and $\psi$ are localised and exponentially decaying as $|x|\to\infty$ or, to be more precise, their asymptotic forms are
\begin{eqnarray}
V(x)&=&p e^{-\lambda|x|}+O(e^{-3\lambda_r|x|}),\nonumber\\
\psi^r(x) &=& q\sgn(x)e^{-\overline{\lambda}|x|}+O(e^{-3\lambda_r|x|}),\label{2.tails}\\
\psi^g(x) &=& s e^{-\overline{\lambda}|x|}+O(e^{-3\lambda_r|x|})\nonumber,
\end{eqnarray}
where $p,q,s$ are complex numbers to be determined numerically and $\lambda=\lambda_r+\ri\lambda_i$ satisfies the dispersion relation \eqref{0.disp} with $\lambda_r>0$; see \cite{ZelMiel09} for more details. Moreover, due to the Fredholm alternative, the equation
\begin{equation*}
\mathcal Lw=\htilde,
\end{equation*}
is solvable if and only if
\begin{equation}\label{e:projection}
P \htilde=0, ~~~ \text{ where } ~~~ P\htilde:=\langle \htilde,\psi^r\rangle\varphi^r+\langle \htilde,\psi^g\rangle\varphi^g.
\end{equation}
Thus, it seems natural to split the solution $w$ of the linearised equation near the pulse $V$ into the sum of the neutral modes $Pw$ and the transversal part $(1-P)w$.
\par
Finally, we introduce the spatial and phase shifted versions of the direct and adjoint eigenfunction associated with the linearization near the shifted pulse $V_\xi$ as
\begin{eqnarray*}
\phi_{\xi}^r(x) &:=&e^{\ri g}\phi^r(x-r),\ ~~~~~\ \phi_{\xi}^g(x) :=e^{\ri g}\phi^g(x-r),\\
\psi_{\xi}^r(x) &:=&e^{\ri g}\psi^r(x-r),\ ~~~~~\ \psi_{\xi}^g(x) :=e^{\ri g}\psi^g(x-r).
\end{eqnarray*}
These shifted functions also satisfy \eqref{adjointcgl2}, namely
\begin{equation}\label{adjointcgl2-shift}
\langle \phi^r_\xi,\psi^r_\xi\rangle=\langle\phi^g_\xi,\psi^g_\xi\rangle=1,\ ~~~~ \
\langle \phi^r_\xi,\psi^g_\xi\rangle=\langle\phi^g_\xi,\psi^r_\xi\rangle=0,
\end{equation}
and equation $\mathcal L_\xi w=\htilde$ is solvable if and only if $P_\xi \htilde:=\langle \htilde,\psi^r_\xi\rangle\varphi^r_\xi(x)+\langle \htilde,\psi^g_\xi\rangle\varphi^g_\xi(x)=0$.

\subsection{The multi-pulse manifold and the projected system}
\label{sec:multimanifold}
We now have the required background to introduce the local coordinates close to the multi-pulse configuration given by
\begin{equation}\label{mult}
V_{\vec\xi}(x):=\sum_{i=1}^nV_{\xi_i}(x),
\end{equation}
where $\vec\xi:=(\xi_1,\cdots,\xi_n)$ and $\xi_i=(r_i,g_i)$ is the position and phase of $i^{\rm th}$ pulse; assuming that the solitons are well-separated i.e
\begin{equation}\label{2.dist}
{\ctext r_1<r_2<\cdots<r_n\ \ \text{and }\ {\ctext d:=\min_{i\ne j}|r_i(t)-r_j(t)|>d_0\gg1}},
\end{equation}
{\ctext where $d_0$ is the minimal separation distance defined in \eqref{0.dist}}. 
Namely, if we define the multi-pulse manifold, {\ctext such that given $d_0>0$ sufficiently large}
\[
{\ctext \mathcal P_n=\mathcal P_n(d):=\{V_{\vec\xi},\ \vec\xi\in(\R\times\mathbb S_1)^n,\ d\geq d_0\}},
\]
then, as proved in \cite{ZelMiel09},  any $u\in L^\infty(\R)$ which is close to $\mathcal P_n$, say, in the $L^\infty$-metric, can be presented in a unique way in the form
\begin{equation}\label{2.ort}
u(x)=V_{\vec\xi}(x)+w(x),\ ~~{\rm such~that}~~\ P_{\xi_i}w=0,\ ~~{\rm for}~~i=1,\cdots, n.
\end{equation}
Our task now is to rewrite the initial system \eqref{e:CGLE} in the new coordinates $\vec\xi=\vec\xi(t)$ and $w=w(t,x)$ in the neighbourhood of $\mathcal P_n$. To this end, we substitute
\begin{equation}\label{ansatzcgl}
u(t,x)=V_{\vec\xi(t)}(x)+w(t,x),\ ~~~~\ P_{\xi_i(t)}w(t)=0,\ ~~{\rm for}~~i=1,\cdots, n,\ \ t\in\R,
\end{equation}
into (\ref{0.init}) and use the definitions
\begin{equation}\label{e:pulseq}
\partial_t V_{\xi_k}=\dot{r}_k\phi^r_{\xi_k} +\dot{g}_k\phi^g_{\xi_k}~~~{\rm and}~~~\alpha\partial_{xx}V_{\xi_k}+\beta V_{\xi_k}+f(V_{\xi_k})=0,
\end{equation}
to yield the equation for the remainder function $w(t,x)$ as
\begin{equation}\label{coup:PDEcgl}
w_t-\alpha w_{xx}-\beta w-f'(V_{\vec\xi})w=-\sum_{k=1}^n \dot{r}_k\phi^r_{\xi_k} - \sum_{k=1}^n \dot{g}_k\phi^g_{\xi_k}
+\Phi(\vec\xi)+G(\vec\xi,w).
\end{equation}
Here the dots denote $\frac{{\rm d}}{{\rm d}t}$, 
\begin{equation} 
\Phi(\vec\xi)=\Phi(\vec\xi,x):=f(V_{\vec\xi}(x))-\sum_{k=1}^n f(V_{\xi_k}(x)),
\end{equation}
is the so called interaction function, and the function 
\begin{equation}
G(\vec\xi,w):=f(V_{\vec\xi}+w)-f(V_{\vec\xi})-f'(V_{\vec\xi})w,
\end{equation}
is of second order in $w$, i.e.
\[
G(\vec\xi,0)=G_w(\vec\xi,0)=0.
\]
We now seek to derive time-evolution equations for the coordinates $\xi_k$ using the orthogonality relations  \eqref{ansatzcgl}. Namely, by the definition of the adjoint eigenfunctions $\psi_{\xi_k}$,
\begin{equation}
\begin{aligned}
\langle w_t-\alpha w_{xx}-\beta w,\psi^r_{\xi_k}\rangle &=& \frac \d{\d t}\langle
w,\psi^r_{\xi_k}\rangle+\langle w,-\partial_t\psi^r_{\xi_k}-\bar\alpha\partial_{xx}\psi^r_{\xi_k}-\bar\beta\psi^r_{\xi_k}\rangle,\\
  &=& \langle f'(V_{\xi_k})w,\psi^r_{\xi_k}\rangle+\dot{r}_k\langle
w,\partial_x\psi^r_{\xi_k}\rangle+\dot{g}_{\xi_k}\langle \ri w,\psi^r_{\xi_k}\rangle, \\
\langle w_t-\alpha w_{xx}-\beta w,\psi^g_{\xi_k}\rangle  &=& \frac \d{\d t}\langle
w,\psi^g_{\xi_k}\rangle+\langle w,-\partial_t\psi^g_{\xi_k}-\bar\alpha\partial_{xx}\psi^g_{\xi_k}-
\bar\beta\psi^g_{\xi_k}\rangle, \\
  &=& \langle f'(V_{\xi_k})w,\psi^g_{\xi_k}\rangle + \dot{r}_k\langle
w,\partial_x\psi^g_{\xi_k}\rangle+\dot{g}_{\xi_k}\langle \ri w,\psi^g_{\xi_k}\rangle,
\end{aligned}
 \end{equation}
using the fact that $w$ is transversal to the neutral modes, i.e that
\begin{equation}\label{transversal}
\langle w, \psi^r_{\xi_k} \rangle=\langle w, \psi^g_{\xi_k} \rangle=0 ~~~~\forall t.
\end{equation}
Thus, taking inner products of \eqref{coup:PDEcgl} with $\psi_{\xi_k}^r$ and $\psi_{\xi_k}^g$, we get the following system of ODEs
\begin{equation}\label{3.ODE}
\begin{aligned}
\sum_{i=1}^n \left[\dot{r}_i\langle \phi_{\xi_i}^r,\psi_{\xi_k}^r \rangle + \dot{g}_i\langle \phi_{\xi_i}^g,\psi_{\xi_k}^r \rangle\right] +\dot{r}_k\langle
w,\partial_x\psi_{\xi_k}^r\rangle+\dot{g}_k\langle \ri w,\psi^r_{\xi_k}\rangle  =  \left\langle \Phi(\vec\xi)+\Psi_k(\vec\xi)w+G(\vec\xi, w),\psi_{\xi_k}^r\right\rangle, \\
\sum_{i=1}^n \left[\dot{r}_i\langle \phi_{\xi_i}^r,\psi_{\xi_k}^g \rangle +
\dot{g}_i\langle \phi_{\xi_i}^g,\psi_{\xi_k}^g \rangle\right] +\dot{r}_k\langle
w,\partial_x\psi_{\xi_k}^g\rangle+\dot{g}_k\langle
\ri w,\psi_{\xi_k}^g\rangle = \left\langle \Phi(\vec\xi)+\Psi_k(\vec\xi)w+G(\vec\xi, w),\psi_{\xi_k}^g\right\rangle,
\end{aligned}
\end{equation}
where $\Psi_k(\vec\xi):=f'(V_{\vec\xi})-f'(V_{\xi_k})$, and $k=1,...,n$.
\par
Note that the system \eqref{3.ODE}  is not explicitly resolved with respect to variables $\dot{r}(t):=(\dot{r}_1(t),\cdots,\dot{r}_n(t)),\;\dot{g}(t):=(\dot{g}_1(t),\cdots,\dot{g}_n(t))$ and, a priori, may be ill posed. However, due to  \eqref{2.tails} and \eqref{adjointcgl2-shift}, if the separation distance $d=d(r)$ between pulses defined by \eqref{2.dist}  is large enough, the inner products
$\langle \phi_{\xi_k}^r,\psi^r_{\xi_i}\rangle$ and $\langle \phi_{\xi_k}^g,\psi_{\xi_i}^g\rangle$
are close to zero for $k\neq i$ (exponentially with respect to the distance
between pulses). Thus, the matrix
\begin{equation}\label{e:cmatr}
 \mathcal C(r,g,w)=\mathcal C(\vec\xi,w) := \begin{pmatrix}
    [\langle  \phi^r_{\xi_i} ,\psi^r_{\xi_j} \rangle + \delta_{ij}\langle w,\partial
_{x} \psi^r_{\xi_j} \rangle]_{i,j=1}^n & [\langle  \phi^g_{\xi_i} ,\psi^r_{\xi_j} \rangle
+ \delta_{ij}\langle \ri w, \psi^r_{\xi_j} \rangle]_{i,j=1}^n  \\
[\langle  \phi^r_{\xi_i} ,\psi^g_{\xi_j}\rangle + \delta_{ij}\langle w,\partial
_{x} \psi^g_{\xi_j}
]_{i,j=1}^n   & [\langle  \phi^g_{\xi_i} ,\psi^g_{\xi_j} \rangle + \delta_{ij}\langle
\ri w, \psi^g_{\xi_j}\rangle]_{i,j=1}^n
    \end{pmatrix},
 \end{equation}
on the left-hand side of \eqref{3.ODE}, is close to the
identity if the pulses are well-separated and the remainder function $w$ is small.  In this case the matrix is invertible and we may write
\begin{equation}\label{3.ODEint}
\frac {\rm d}{{\rm d}t}\(\begin{matrix} r_k\\g_k\end{matrix}\)=\mathcal C^{-1}(\vec \xi,w)\(\begin{matrix}
\left[\left\langle \Phi(\vec\xi)+\Psi_k(\vec\xi)w+G(\vec\xi, w),\psi_{\xi_k}^r\right\rangle\right]_{k=1}^n\\ \left[\left\langle \Phi(\vec\xi)+\Psi_k(\vec\xi)w+G(\vec\xi, w),\psi_{\xi_k}^g\right\rangle\right]_{k=1}^n
\end{matrix}\),
\end{equation}
and denoting the right-hand side of this equation by $H$, we get
\begin{equation}\label{3.ODE1}
\frac {\rm d}{{\rm d}t}\vec\xi={\ctext H(\epsilon,\vec\xi,w)},
\end{equation}
which are the desired equations for the slow components $\vec\xi$. Inserting these equations to the right-hand side of \eqref{coup:PDEcgl}, we remove the terms containing the derivatives of slow variables and arrive at
\begin{equation}\label{3.PDE1}
w_t-\alpha w_{xx}-\beta w-f'(V_{\vec\xi})w=-h(\vec\xi,w)
+\Phi(\vec\xi)+G(\vec\xi,w),
\end{equation}
where 
\begin{eqnarray*}
h(\vec\xi,w) &=& \sum_{k=1}^n \mathcal{C}^{-1}(\vec\xi,w)\left\langle \Phi(\vec\xi)+\Psi_k(\vec\xi)w+G(\vec\xi, w),\psi_{\xi_k}^r\right\rangle\phi^r_{\xi_k}\\ &&+ \sum_{k=1}^n \mathcal{C}^{-1}(\vec\xi,w) \left\langle \Phi(\vec\xi)+\Psi_k(\vec\xi)w+G(\vec\xi, w),\psi_{\xi_k}^g\right\rangle\phi^g_{\xi_k}.
\end{eqnarray*}
Finally, using the orthogonality conditions \eqref{ansatzcgl}, we rewrite this PDE as 
\begin{equation}\label{3.PDE2}
w_t-\mathcal L(\vec\xi)w=-h(\vec\xi,w)
+\Phi(\vec\xi)+G(\vec\xi,w),
\end{equation}
where $\mathcal L(\vec\xi):=\alpha\partial_x^2+\beta+f'(V_{\vec\xi})+P_{\vec\xi}$ and $P_{\vec\xi}\,\theta:=\sum_{i=1}^nP_{\xi_i}\theta$.

Equations \eqref{3.ODE1} and \eqref{3.PDE2} give the desired projected system which consists of a `fast' PDE and a `slow' set of ODEs. The extra term $P_{\vec\xi}$ is added to the operator $\mathcal L(\vec\xi)$ (following \cite{ZelMiel09}) in order to remove the neutral (or almost neutral) modes, which correspond to pulse shifts and rotations, making it invertible and stable. Indeed, on the one hand, this extra term does not destroy the equations since it vanishes due to the orthogonality conditions on $w$ and, on the other hand it makes $\mathcal L(\vec\xi)$ uniformly invertible if the separation distance, $d$, is large. 

The key result on the existence of a stable manifold for the fast-slow system \eqref{3.ODE1}-\eqref{3.PDE2} is given in the following theorem proved in \cite{ZelMiel09}.

\begin{theorem}[\cite{ZelMiel09}]\label{Th.main} Let the pulse $V$ satisfy the stability assumptions \eqref{2.stable}. Then, for all separation distances, {\ctext $d\geq d_0$, for sufficiently large $d_0$}, there exists a (small) neighbourhood of $\mathcal P_n(d)$ (which is independent of $d$), say, in the metric of $L^\infty(\mathbb R)$ such that the projected system \eqref{3.ODE1}-\eqref{3.PDE2} is equivalent to the initial QCGLE \eqref{e:CGLE} {\ctext as long as} the trajectory $u(t,x)$  remains in this neighbourhood.
\par
Moreover, there exists a $C^2$-smooth function $W:\overline{\mathcal P_n(d)}\to L^\infty(\R)$ such that
\begin{equation}\label{e:west}
\|W\|_{C^2}\le C\eb,\ \ \eb:=e^{-\lambda_rd_0},
\end{equation}
where $\lambda_r>0$ is the same as in \eqref{2.tails}, and the manifold defined by the graph of $W$ ( $w=W(\vec\xi)$, $\xi\in \mathcal P_n(d)$) is an invariant centre manifold of the projected system \eqref{3.ODE1}-\eqref{3.PDE2}. In particular, the dynamics on this manifold is determined by the system of ODEs
\begin{equation}\label{3.ODE3}
\frac d{dt}\vec \xi=H(\epsilon,\vec\xi,W(\vec\xi)),\ \ \|H\|_{C^2}\le C\eb,
\end{equation}
and the trajectory $\vec\xi(t)$ starting on this manifold can leave it only through the boundary $\partial P_n(d)$.
\par
Finally, this manifold is exponentially stable and normally hyperbolic, so for any trajectory
$(\vec\xi(t),w(t))$ of the projected system there is a trace trajectory $(\vec\xi_0(t),W(\vec\xi_0(t))$ on the manifold such that
\begin{equation}\label{trace}
\|\vec\xi(t)-\vec\xi_0(t)\|+\|w(t)-W(\vec\xi_0(t))\|_{L^\infty}\le Ce^{-\kappa t},
\end{equation}
for some positive $\kappa$ which is independent of $d$. Estimate \eqref{trace} holds until the trace trajectory $\vec\xi_0(t)$ reaches the boundary $\partial P_n(d)$.
\end{theorem}
{\ctext The constants $C$, $d_0$ and $\kappa$ in the above theorem a priori depend on the parameters of the considered QCGLE and explode as ${\rm Re}(\beta)=\beta_r\to0$. in particular, following the proof given in \cite{ZelMiel09}, we can estimate that the constant $C$ in \eqref{e:west} has the form $C_0\beta_r^{-1}$ for $C_0$ independent of $\beta_r$. On the other hand, the norms of the asymptotic expansions with respect to $\epsilon$ for the function $H$, considered below, do not explode as $\beta_r\to0$, so we may expect much better behaviour of the centre manifold as $\beta_r\to0$. To the best of our knowledge, the dependence of the manifold as $\beta_r\to0$ remains an open problem, and is beyond the scope of this paper.}

We note that the constant $\kappa>0$ in \eqref{trace} is usually not small (it is related to the spectral bound of the operator $\mathcal L$ after excluding the zero eigenvalues), so if we start from initial data sufficiently close to the manifold (e.g., taking $\vec\xi(0)\in \mathcal P_n(d)$ and $w(0)=0$, we will be $\eb$-close to the manifold), after a relatively short transient behaviour, we will be indistinguishably close to it. Moreover, since the dynamics on the manifold is of $O(\epsilon)$, this transient behaviour can only produce a small, $O(\eb)$, shift of the initial data. Cutting off this transient behaviour, we are then able to recover the dynamics on the manifold by solving the projected system with $w(0,x)=0$.

Note also that the system \eqref{3.ODE3} is a coupled system of $2n$ ODEs which can be effectively solved by standard methods if the function $w$ is known. However, in order to recover $w$, we still need to solve the PDE \eqref{3.PDE2}. Below we consider an expansion for the ODE/PDE system consisting of leading order and second order (with respect to $\epsilon$) terms only. We will show that in cell 1 the first order approximation is not sufficient to accurately resolve the motion, but the second order approximation is, and in this case the PDE reduces to solving a quasi-steady PDE, hence simplifying the system and making it numerically tractable.

\subsection{Asymptotic structure of the fast-slow system}
\label{sec:asy}

In this section we derive the leading order {\ctext pulse} correction terms structure of the ODE/PDE system, i.e. up to $O(\epsilon^2)$. The leading order, $O(\epsilon)$, approximation has been found in \cite{ZelMiel09} and we will follow this approach where appropriate.

From \eqref{e:west} it has been argued that $w=O(\epsilon)$ where $\epsilon=e^{-\lambda_rd_0}$. Also the work in \cite{ZelMiel09} shows at leading order the ODE system is \eqref{e:ode2} where the RHS is $O(\epsilon)$. In order to remove the long time scales from the problem we choose to introduce a new time variable. Hence we introduce the new scaled quantities
\begin{equation}\label{scalings}
\that=\epsilon t,~~~\Phi=\epsilon\Phihat,~~~w=\epsilon\what,~~~\Psi_k=\epsilon\Psihat_k,~~~~G=\epsilon^2\Ghat,
\end{equation}
where hatted variables are $O(1)$. Note, both $\Phi$ and $\Psi_k$ are not $O(\epsilon)$ everywhere in $x\in(-\infty,\infty)$, but they are in those regions where they interact with the adjoint eigenfunction. This interaction leads to inner products such as {\ctext $\langle \Phi,\psi^r_{\xi_k}\rangle=O(\epsilon)$} and {\ctext $\langle \Psi_kw,\psi^r_{\xi_k}\rangle=O(\epsilon^2)$} etc. The first of these scales was shown in \cite{ZelMiel09} while the second is confirmed numerically in Appendix \S\ref{appen:innerprods}.

Inserting these quantities into the ODE/PDE system we find that the matrix
\[
\mathcal C(\vec\xi,w)=\mathcal C(\vec\xi,0)+\epsilon\widehat{\mathcal{C}}(\vec\xi,\what),
\]
and, moreover, due to \eqref{2.tails} and the orthogonality relations between direct and adjoint eigenfunctions,
\[
\mathcal C(\vec\xi,w)=\mathbb{I}_{2n}+\epsilon\widehat{\mathcal{C}}(\vec\xi,\what).
\]
Hence its inverse can be easily found and \eqref{3.ODEint} can be written as
\begin{equation}\label{3.ODEasy}
\frac {\rm d}{{\rm d}\that}\(\begin{matrix} r_k\\g_k\end{matrix}\)=\(\begin{matrix}
\left[\left\langle \Phihat(\vec\xi),{\ctext \psi_{\xi_k}^r}\right\rangle\right]_{k=1}^n\\ \left[\left\langle \Phihat(\vec\xi),{\ctext \psi_{\xi_k}^g}\right\rangle\right]_{k=1}^n\end{matrix}\)+\epsilon\left[\(\begin{matrix}
\left[\left\langle \Psihat_k(\vec\xi)\what+\Ghat(\vec\xi, w),\psi_{\xi_k}^r\right\rangle\right]_{k=1}^n\\ \left[\left\langle \Psihat_k(\vec\xi)\what+\Ghat(\vec\xi, w),\psi_{\xi_k}^g\right\rangle\right]_{k=1}^n\end{matrix}\)
-\widehat{\mathcal{C}}\(\begin{matrix}
\left[\left\langle \Phihat(\vec\xi),{\ctext \psi_{\xi_k}^r}\right\rangle\right]_{k=1}^n\\ \left[\left\langle \Phihat(\vec\xi),{\ctext \psi_{\xi_k}^g}\right\rangle\right]_{k=1}^n\end{matrix}\)\right]+O(\epsilon^2),
\end{equation}
or more compactly
\begin{equation}\label{asyODE.compact}
\frac{\d}{\d\that}\vec\xi=\widehat{H}({\ctext \epsilon},\vec\xi,\what)+O(\epsilon^2).
\end{equation}

\subsubsection{The leading order system}

Clearly at $O(1)$ in \eqref{3.ODEasy} the remainder function $\what$ does not feature, hence we need not consider the PDE \eqref{3.PDE2}, and the system reduces to
\begin{equation}\label{3.ODE4}
\frac \d{\d \that}r_k=\langle \Phihat(\vec\xi),\psi_{\xi_k}^r\rangle,\ ~~~ \ \frac \d{\d \that}g_k=\langle \Phihat(\vec\xi),\psi_{\xi_k}^g\rangle,\ \ k=1,\cdots n.
\end{equation}
This system can be further simplified, because based on the tail estimates \eqref{2.tails} and the mean value theorem, one can see that
$$
\Phihat(\vec\xi)={\ctext \frac{1}{\epsilon}}\sum_{i=1}^n f'(V_{\xi_i})(V_{\xi_{i-1}}+V_{\xi_{i+1}})+O(\eb),
$$
where we assume that $V_{\xi_j}=0$ if $j<1$ or $j>n$, see \cite{ZelMiel09} for more details. Furthermore, using \eqref{2.tails} again, we can easily see that
$$
\langle \Phihat(\vec\xi),\psi_{\xi_k}^r\rangle={\ctext \frac{1}{\epsilon}}\langle f'(V_{\xi_k})V_{\xi_{k-1}},\psi_{\xi_k}^r\rangle+{\ctext \frac{1}{\epsilon}}\langle f'(V_{\xi_k})V_{\xi_{k+1}},\psi_{\xi_k}^r\rangle +O(\eb),
$$
and the analogous formulas work also for inner products containing $\psi_{\xi_k}^g$. This allows us to introduce pairwise interaction functions
\begin{equation}
\Theta^r(\xi_k-\xi_{k-1}):={\ctext \frac{1}{\epsilon}\langle f'(V_{\xi_k})(V_{\xi_{k}}-V_{\xi_{k-1}}),\psi^r_{\xi_k}\rangle},\ \ ~~ \Theta^g(\xi_k-\xi_{k-1}):={\ctext \frac{1}{\epsilon}\langle f'(V_k)(V_{\xi_{k}}-V_{\xi_{k-1}}),\psi^g_{\xi_k}\rangle},
\label{eqn:interactionequations}
\end{equation}
and rewrite equations \eqref{3.ODE4} in the form
\begin{equation}\label{3.ODE5}
\frac \d{\d \that} r_k=-\Theta^r(\xi_{k-1}-\xi_k)+\Theta^r(\xi_{k+1}-\xi_k),\ \ ~~ \frac \d{\d \that} g_k=-\Theta^g(\xi_{k-1}-\xi_k)+\Theta^g(\xi_{k+1}-\xi_k),
\end{equation}
where as before, the terms containing $\xi_{j}$ with $j<1$ or $j>n$ are identically zero.
\par
Finally, utilising the structure of the underlying equation \eqref{e:CGLE}, we may  compute functions $\Theta^r$ and $\Theta^g$ in terms of the parameters contained in the asymptotic expansions \eqref{2.tails}. In order to avoid technicalities, we present here only the main result (the details can be found in \cite{ZelMiel09}). Namely,
 \begin{equation}\label{3.Theta}
 \begin{cases}
 \Theta^r(r, g) =M_1 \sin(g - \lambda_i|r| + \kappa_1),\\
 \Theta^g(r,g) =M_2 \operatorname{sgn}({r})\sin(g - \lambda_i|r| + \kappa_2),
 \end{cases}
 \end{equation}
where $\lambda_i$ is the same as in \eqref{0.disp} and the constants $M_1$, $M_2$, $\theta_1$ and $\theta_2$ can be expressed in terms of the coefficients $p$, $q$ and $s$ from \eqref{2.tails} via the following equations
$$
\frac{M_1}{\bar{q}} e^{\ri\kappa_1} = \frac{M_2}{\bar{s}} e^{\ri\kappa_2} = -2p\bar{\alpha}\lambda.
$$
In particular, for the case $n=2$, equations \eqref{3.ODE5} give us the 2D system \eqref{e:ode2} for the differences $\bar r=r_2-r_1>0$ and $\bar g=g_2-g_1$ with $J=-2M_1$ and $K=2M_2$.
\par
We note that the leading order approximation system \eqref{3.ODE5} is {\it reversible}. Indeed, it is invariant with respect to the transformation
$$
\that\to -\that,\ \ \ g_k\to g_k+\pi .
$$
This explains the coincidence of stable and unstable manifolds and formation of cells seen in the case $n=2$, see figure \ref{fig:Sergeyppred}. However, this reversibility does not take place for the QCGLE and is already broken when the second order corrections are included. Moreover, due to this reversibility, the leading order approximation \eqref{3.ODE5} may not be structurally stable even in the simplest case $n=2$ and we need higher order corrections in order to understand the weak interaction of the pulses.

\subsubsection{The second order system}
\label{sec:secondorder}

In \eqref{3.ODEasy}, if we retain the $O(\epsilon)$ terms then it is clear that the remainder function $\what$ now becomes significant and hence needs to be solved for. Substituting the scalings \eqref{scalings} into \eqref{3.PDE2} leads to
\begin{equation}\label{3.PDEasy}
\epsilon^2\what_{\that}-\epsilon\mathcal L(\vec\xi)\what=-\epsilon \hhat(\vec\xi,\what)+\epsilon\Phihat(\vec\xi)+\epsilon^2\Ghat(\vec\xi,\what).
\end{equation}
Hence retaining the leading order components of $\what$ gives the equation
\begin{equation}\label{leadingw}
\mathcal L(\vec\xi)\what=-\Phihat(\vec\xi)+\hhat(\vec\xi,0),
\end{equation}
which is solved together with the fast system \eqref{3.ODEasy}, where
\[
\hhat(\vec\xi,0)=\sum_{k=1}^n \left\langle \Phi(\vec\xi),\psi_{\xi_k}^r\right\rangle\phi^r_{\xi_k}+ \sum_{k=1}^n \left\langle \Phi(\vec\xi),\psi_{\xi_k}^g\right\rangle\phi^g_{\xi_k}.
\]
The numerical evidence for the justification of this system and its efficiency are examined in the next section.

\section{Numerical implementation and justification of projected system}\label{sec:Impl_Just}

\subsection{Computation of the pulse and neutral modes}
To compute the localised pulse solution and the corresponding eigenfunctions and
adjoint eigenfunctions of the QCGLE, we set up the boundary value problem for $V(x)$ on the truncated finite domain $x\in[0,L]$, given by
\begin{equation}\label{e:matlabrhs}
 \alpha V_{xx} + \beta V+ f(V) = 0.
\end{equation}
We split $V~(=u+\ri v)$ into real and imaginary parts, and solve the coupled pair of equations for these functions. Given the even symmetry of the pulse solution sought, we impose the boundary conditions
\begin{equation}\label{e.lboundary}
\left.\frac{\d u}{\d x}\right|_{x=0}=\left.\frac{\d v}{\d x}\right|_{x=0}=0,
\end{equation}
while at $x=L$ we require $u$ and $v$ to have the correct exponential decay at these boundaries, hence we impose
\begin{equation}\label{e.rboundary}
\frac{\d u}{\d x}+\lambda_r u-\lambda_i v=\frac{\d v}{\d x}+\lambda_r v+\lambda_i u=0,
\end{equation}
where $\lambda=\lambda_r+\ri\lambda_i$ satisfies \eqref{0.disp}. 

The four complex parameters $\alpha,~\beta,~\gamma$ and $\delta$ give eight real parameters, of which we fix seven, leaving $\beta_i=\Im(\beta)$ to be found by imposing the additional equation \eqref{0.disp}, which we again split into real and imaginary parts. Equations \eqref{e:matlabrhs} together with \eqref{e.lboundary}, \eqref{e.rboundary} and \eqref{0.disp} are then solved using the \texttt{chebop} routine from chebfun for nonlinear boundary value problems to machine precision \cite{Driscoll2014}. Once the solution is obtained, we extend it to $x\in[-L,0]$ via an even extension.

To calculate the zero eigenfunctions and adjoint eigenfunctions we use the same chebfun structure as for the pulse itself, but we instead solve the linear and adjoint problems $\mathcal{L}(u_1+\ri v_1)=0$ and $\mathcal{L}^*(u_1+\ri v_1)=0$, together with the boundary conditions
\[
u_1(0)=v_1(0)=u_1(L)=v_1(L)=0,
\]
for the odd eigenfunctions, and 
\[
\left.\frac{\d u_1}{\d x}\right|_{x=0}=\left.\frac{\d v_1}{\d x}\right|_{x=0}=u_1(L)=v_1(L)=0,
\]
for the even eigenfunctions. Here $u_1$ and $v_1$ are the real and imaginary parts of the regular or adjoint eigenfunctions, and we extend them to $x\in[-L,0]$ via an even or odd extension as required.

Once we have numerical approximations for the zero eigenfunctions and adjoint eigenfunctions of the linearised operator $\mathcal L$, we need to normalise these for the projected system. To do this we note that eigenfunctions satisfy the conditions
\begin{equation}\label{e:phasecond34}
\Re \left( \int_{-L}^{L} \phi^r \overline{\phi^r} \d x \right) =
\Re \left( \int_{-L}^{L} V_x \overline{V_x} \d x \right), \quad ~~~ \Re \left(
\int_{-L}^{L}\phi^g \overline{\phi^g} \d x \right) = \Re \left( \int_{-L}^{L} V \overline{V} \d x \right).
\end{equation}

Once the eigenfunctions are appropriately normalised we then scale the adjoint eigenfunctions using the following conditions
\begin{equation}\label{e:phasecond56}
\Re \left( \int_{-L}^{L}\phi^g \overline{\psi^g} \d x \right)=\Re \left( \int_{-L}^{L}\phi^r \overline{\psi^r} \d x \right) =  1, \quad  ~~ \Re \left
( \int_{-L}^{L}\phi^r \overline{\psi^g} \d x \right)=\Re \left( \int_{-L}^{L}\phi^g \overline{\psi^r} \d x \right) =  0.
\end{equation}
The final two conditions in \eqref{e:phasecond56} are satisfied exactly due to the odd and even extensions of these functions used.

The projection scheme requires translations and phase rotations of the pulse and the neutral eigenfunctions. But as the functions are constructed via Chebyshev expansions, translations of the these functions are trivial, and thus we can very readily find approximations to the functions $V_{\xi_k},\phi^r_{\xi_k},\phi^g_{\xi_k},\psi^r_{\xi_k}$ and $\psi^g_{\xi_k}$ for any given translation $r_k,$ and phase rotation $g_k$.

The majority of the numerical results presented in this paper {\ctext focus} on the $n=2$ pulse case with $L=10$, and hence we have the 4 variables $(r_1,g_1,r_2,g_2)$ with $\rbar=r_2-r_1>0$ and $\gbar=g_2-g_1$. In order to justify our results and to investigate parameter space we consider the following 3 schemes to solve this two-pulse system:

\begin{itemize}
\item {\bf Standard time-stepping Scheme (SS)} - In this approach we truncate the infinite domain to $x\in[-3,3]$ and discretize the domain using a spatial step size $\Delta x=10^{-2}$. The solution to the QCGLE, $u(t,x)$, is then evaluated on this grid and an $8^{\rm th}$ order finite difference scheme is used to approximate the second derivative. In order to time evolve the stiff system of equations we use \textsc{matlab}'s \texttt{ode15s} solver which has variable order adaptive time-stepping \cite{gaines1997,Shampine94,Shampine97}. The values of $(r_1,g_1,r_2,g_2)$ are found in post-processing by requiring the transversal conditions \eqref{transversal} hold for each of the four adjoint eigenfunctions. These integrals are carried our using Boole's rule. For more background on this numerical approach see \cite{rossides_thesis}.

\item {\bf Projected System (PS)} - This approach solves the second order system documented in \S\ref{sec:secondorder}. Firstly we solve \eqref{leadingw} for $\what$, and use this in \eqref{asyODE.compact} with the $O(\epsilon^2)$ terms neglected. The resulting ODE system is then evolved again using \texttt{ode15s}. In order to solve \eqref{leadingw} we truncate to $x\in[-L,L]$, with $L=10$, and discretize the domain using a spatial step size $\Delta x=10^{-2}$. Again employing an $8^{\rm th}$ order central finite differences for the second derivative term we solve for $\what$ with boundary conditions
\[
\left.\frac{\d\what}{\d x}\right|_{x=-L}=\left.\frac{\d\what}{\d x}\right|_{x=L}=0.
\]
The inner products in \eqref{asyODE.compact} are calculated using Boole's rule, and hence the absolute error of the inner products and spatial differentiation is of  $O(10^{-16})$, i.e. machine precision. {\ctext Note, when calculating $\what$, the linear operator $\mathcal{L}$, defined just after \eqref{3.PDE2}, has the projection $P_{\vec\xi}$ added in order that the matrix formed from the finite difference approximation is invertible, due to the normal modes having been removed.}

\item {\bf Projected ODE System (POS)} - This scheme is the same as the PS scheme above, but rather than {\ctext solving} for $\what$ at each time-step in the \eqref{asyODE.compact} system, $\what$ is calculated prior to the time-integration, and the asymptotic forms of the inner products for $\rbar\gtrsim1.5$ are utilized, see Appendix \ref{appen:innerprods} for full details. The benefit of this system is that it is computationally fast, and as long as we are in the asymptotic regime $\rbar>1.5$ it has the same accuracy as the PS scheme, as we will see in the next session.

\end{itemize}

All of the results presented in this section have been computed on a machine with the following specifications: Intel Core i7-6700U CPU @ 3.40GHz x 8, 16GB RAM, using Ubuntu 18.04 LTS in \textsc{matlab} version R2021b.

\subsection{Numerical validation of Projected Scheme (PS) and Projected ODE Scheme (POS)}\label{subs:justification}

For the justification of the projected scheme (PS), and subsequently the projected ODE scheme (POS), we first confirm the asymptotic scalings used on $w$ and $w_t$ in \S\ref{sec:asy}, and then we compare results from the PS with the full numerical SS results to show they produce quantitatively correct and identical dynamics. Finally, we examine computational timings to show the major speed up of the PS and POS over the SS.

\begin{figure}[htb!]
     \begin{center}
(a){\includegraphics[width=0.45\textwidth]{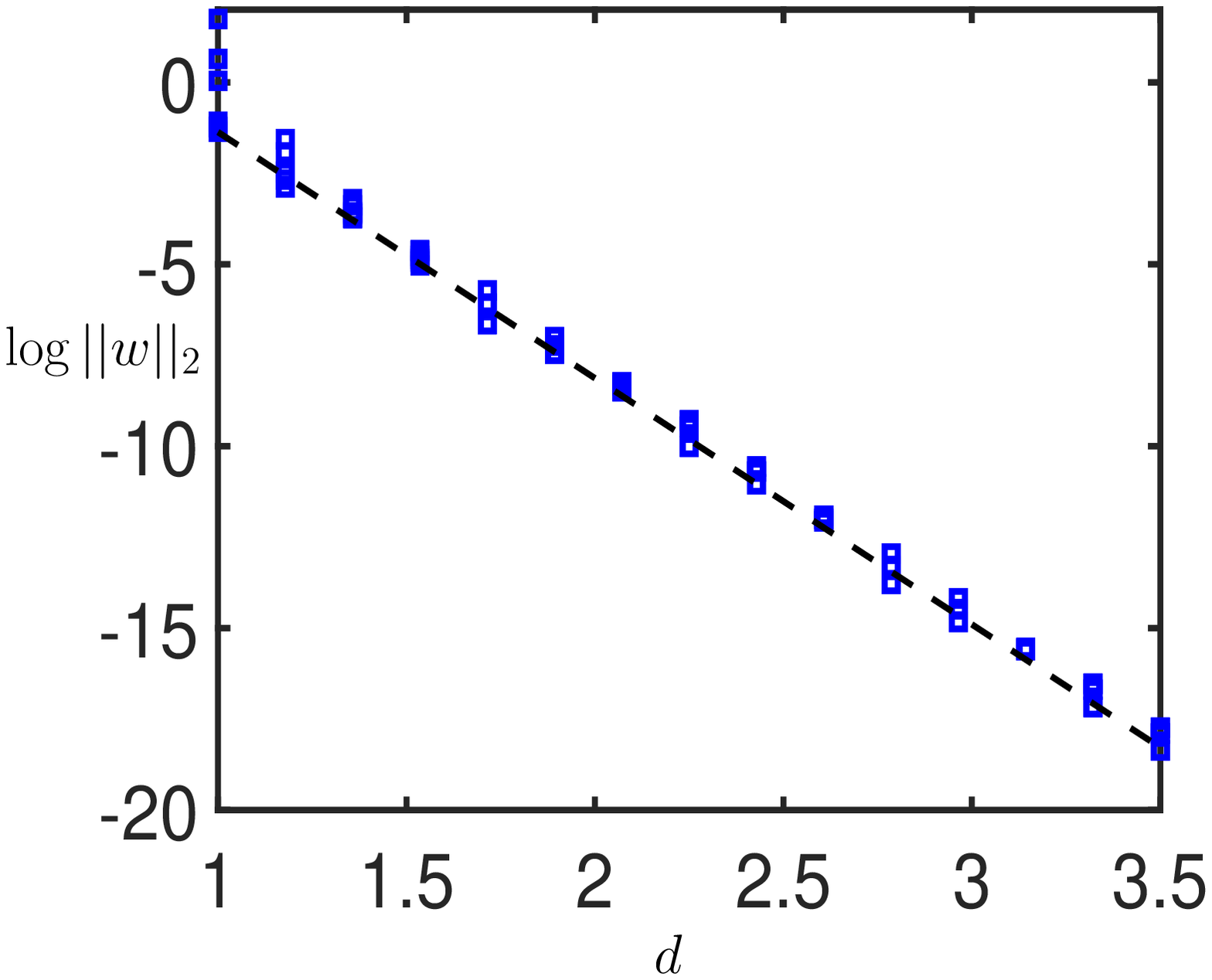}}
(b){\includegraphics[width=0.45\textwidth]{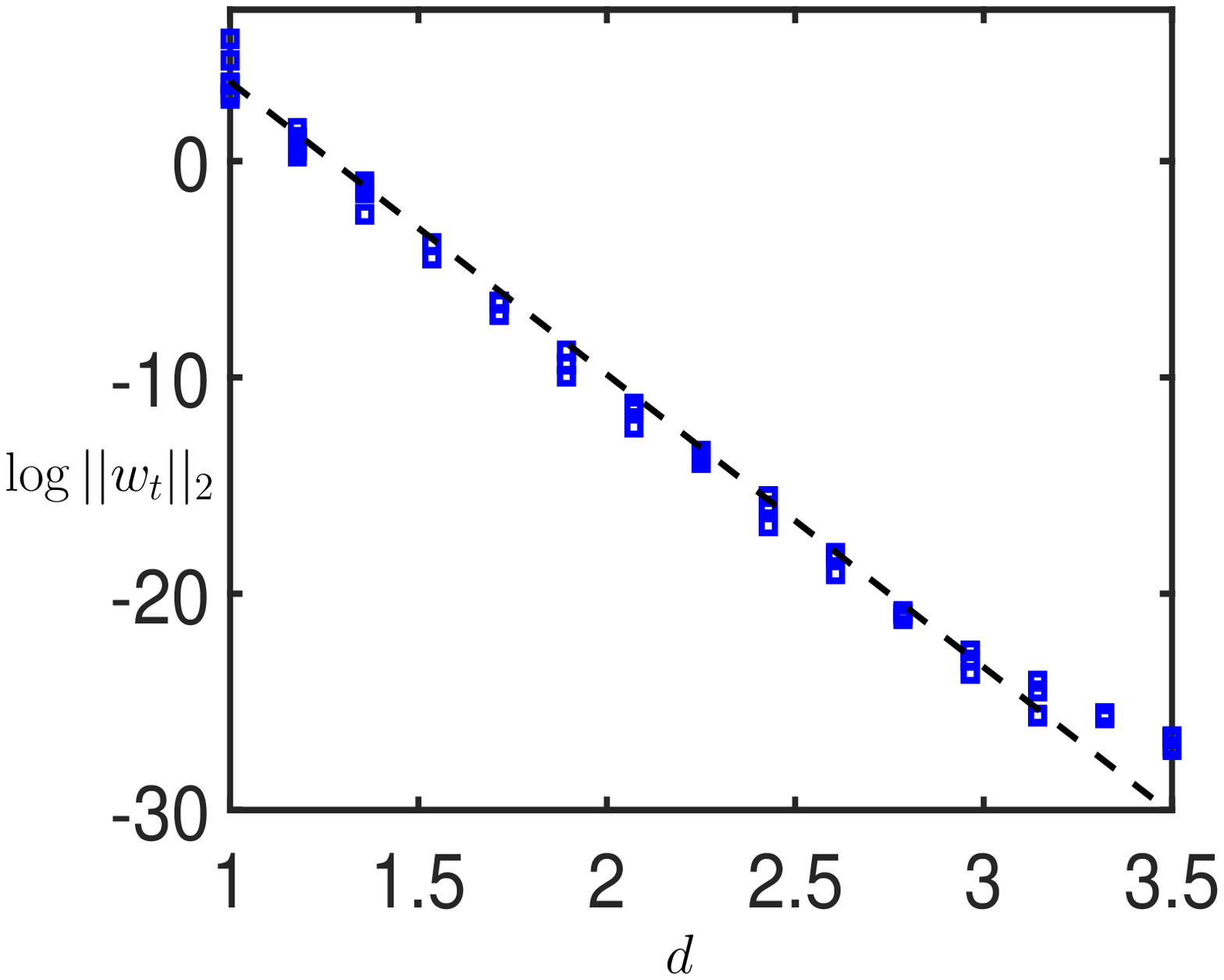}}
    \end{center}
    \caption{(colour online) Comparison of numerical and analytical magnitudes of (a) $\log\|w \|_2$ and (b) $\log\|w_t \|_2$ for separation distance $d=|\rbar| = [1 , 3.5]$ and $\gbar = [0 , \pi]$, for  the system parameters \eqref{eqn:parameters}. The dashed lines in panels (a) and (b) correspond to the analytical estimate for the magnitude decay rates of $-\lambda_r$ and $-2\lambda_r$ respectively.}
   \label{fig:normcomp}
\end{figure}   
For the test cases presented in this section, we consider the parameter values 
\begin{equation}\label{eqn:parameters}
\alpha=0.5+0.5\ri,~~~\beta = -0.02-37.91\ri,~~~\gamma=1.8+\ri,~~~\delta=-0.05+0.05\ri. 
\end{equation}
In figure \ref{fig:normcomp}, we plot the 2-norm of $w=\epsilon\what$ and $w_t=\epsilon^2\what_{\that}$ from solving \eqref{leadingw} and \eqref{3.PDE2} respectively for $d=|\rbar|\in[1,3.5]$ and $\gbar\in[0,\pi]$. The figure shows that the size of these terms agree with the asymptotic prediction in \S\ref{sec:asy}, i.e. that
\begin{equation}\label{wnormest}
\|w\|_2=O(e^{-\lambda_rd_0}),  ~~~\quad \|w_t\|_2=O(e^{-2\lambda_rd_0}).
\end{equation}
Here the dashed line in each panel gives the analytic decay rate $-\lambda_r$ and $-2\lambda_r$ respectively where $\lambda_r=6.77$. This figure thus validates the asymptotic scalings used in in deriving the PS and POS, and highlights that the time derivative term can be neglected from the PDE for the remainder function, leaving a much simpler quasi-steady PDE at leading order.

\begin{figure}[htb!]
     \begin{center}
(a){\includegraphics[width=0.45\textwidth]{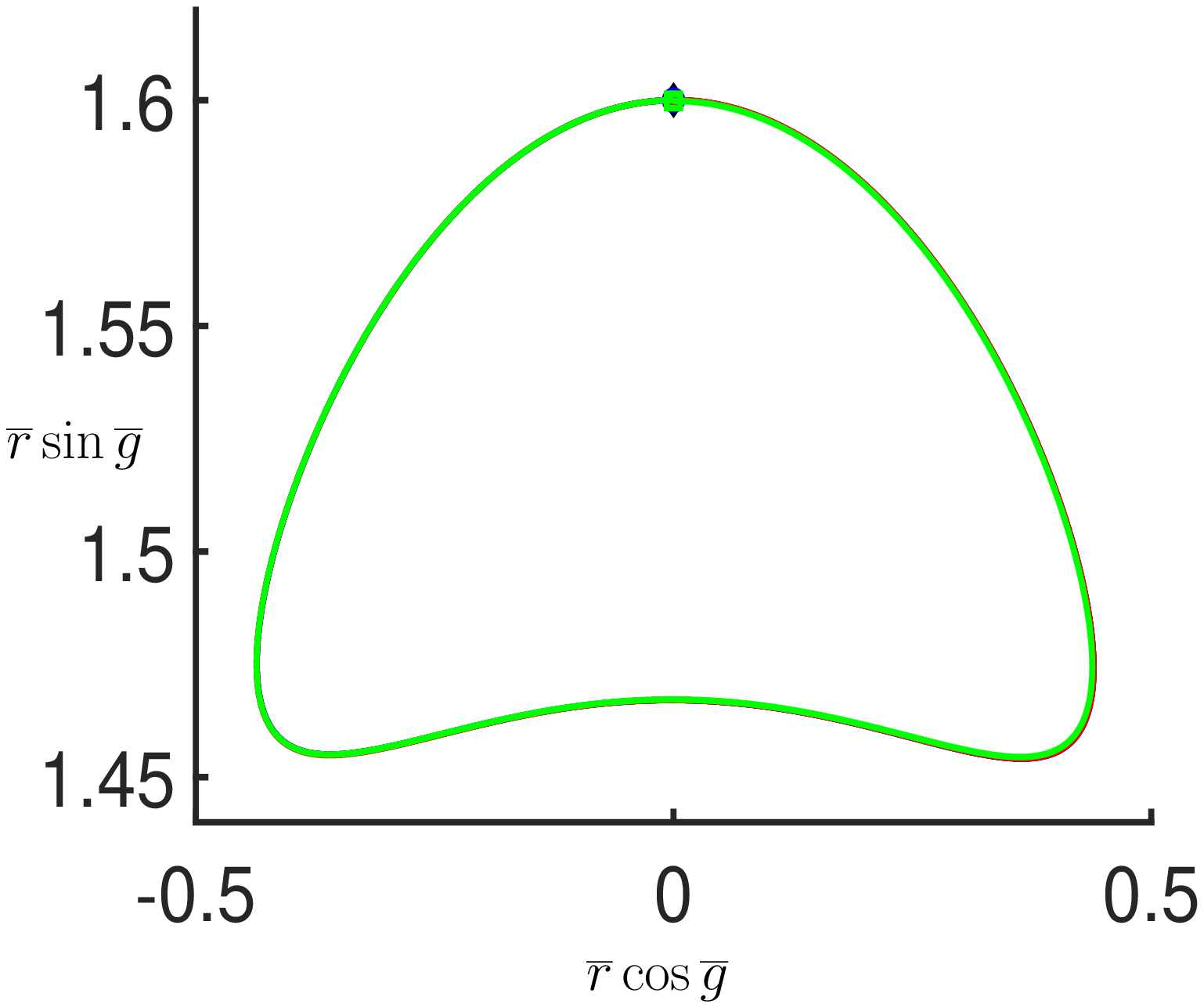}}
(b){\includegraphics[width=0.45\textwidth]{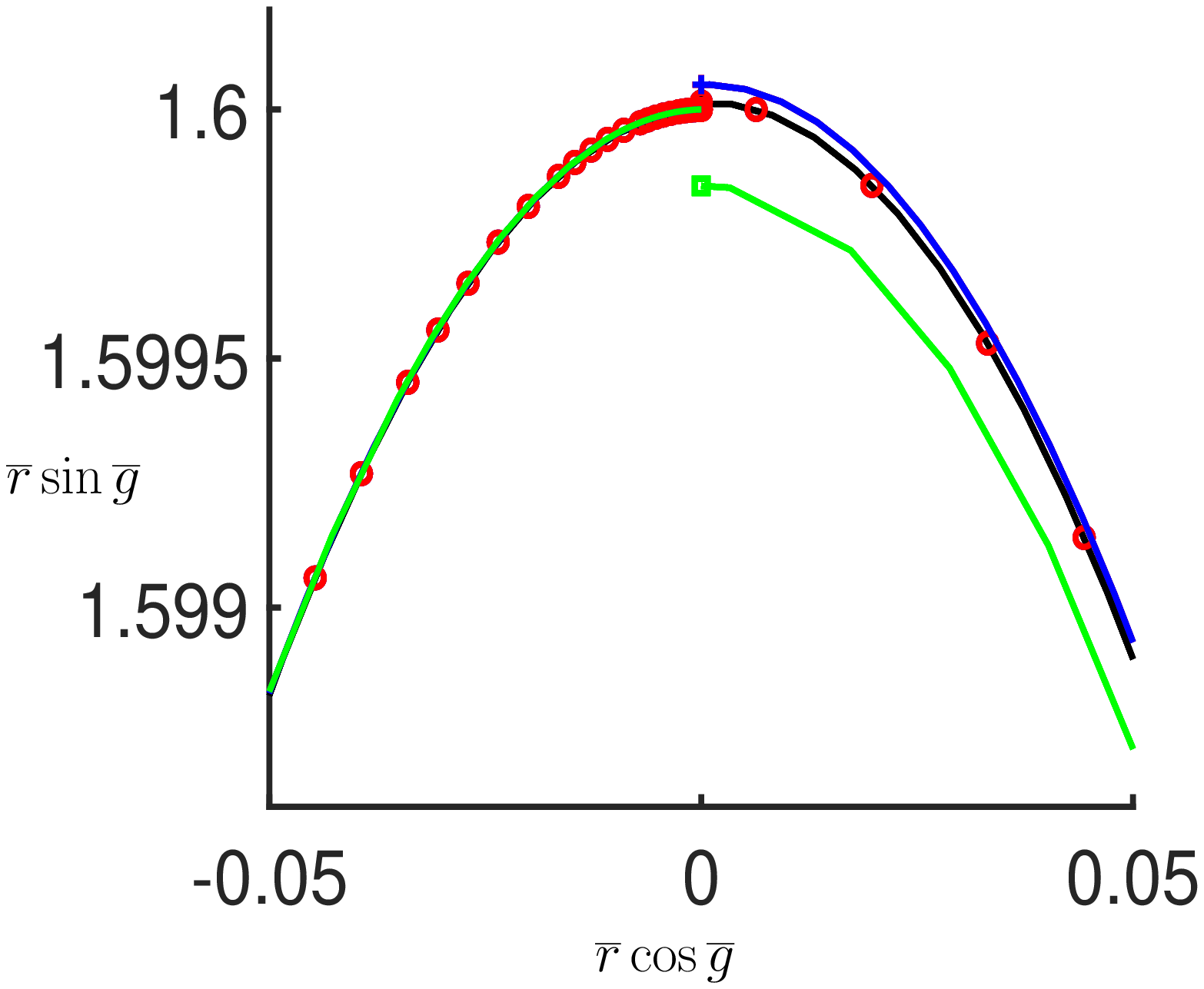}}
(c){\includegraphics[width=0.45\textwidth]{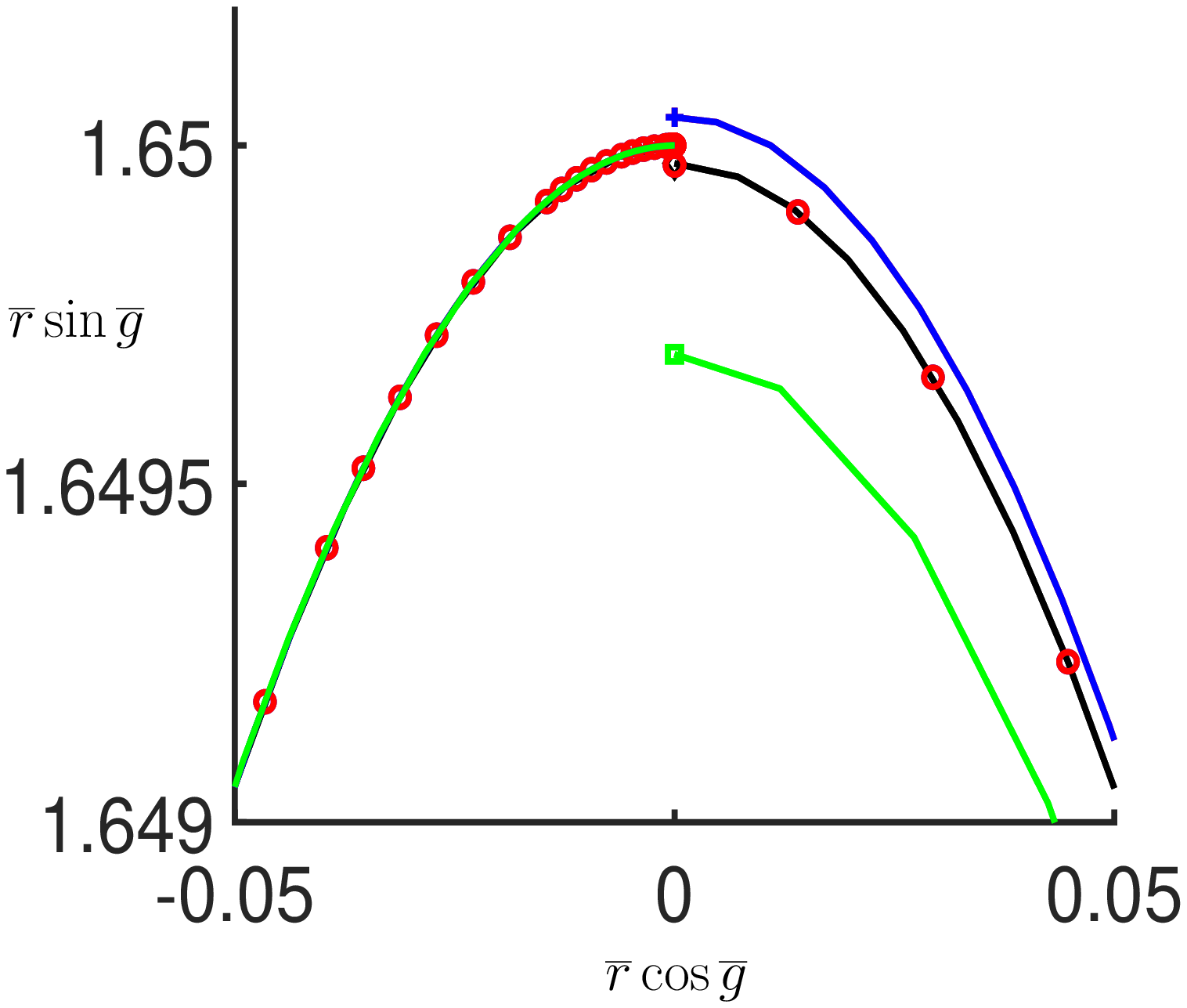}}
    \end{center}
    \caption{(colour online) Poincare section for the phase portrait trajectory with initial condition $\gbar(0)=\frac{\pi}{2}$ and (a,b) $\rbar(0)=1.6$ and (c) $\rbar(0)=1.65$. Each trajectory is integrated until $\that=T$ where, in the leading order system, the trajectory should have returned to $(\rbar(0),\gbar(0))$. Each panel contains 4 results, the SS result (black), the PS result (red) the POS result (blue), and the PS result with $\what\equiv0$ (green). In panels (b) and (c) the red line is represented by the symbols so comparison with the other results is possible. The symbols at $\rbar\cos\gbar=0$ denote the value at $(\rbar(T),\gbar(T))$.}
   \label{fig:strob_sec}
\end{figure} 
In figure \ref{fig:strob_sec} we make a comparison between the phase plane characteristics of the PS and POS with the full numerical SS, specifically, by considering solution trajectories in the first cell of the phase plane. We can make this comparison by introducing a Poincare (stroboscopic) section of the dynamics. To do this we fix the initial phase difference between pulses as $\gbar(0)= \frac{\pi}{2}$, and vary the initial separation distance, $\rbar(0)$. By doing this we create a phase plane where we can study the stability of each trajectory with initial starting point $(\rbar(0),\gbar(0))$ by examining the separation distance $\rbar(T)$ at the next time when trajectory intersects the Poincare section. The time $T$, if it exists, is detected using the event detection feature in \textsc{matlab}'s \verb2ode15s2 subroutine. The solver's event detection will stop the solver at some time $T$ when
\[
\gbar(T)-\gbar(0)=0.
\]
However, the value $T$ is not unique since there exist 2 values for which the above condition is satisfied,  therefore we also impose one extra condition regarding the sign of the derivative $\gbar_{\that}(T)$. As the trajectory directions change in each successive cell, the sign of $\gbar_{\that}(T)$ needs to be made positive or negative depending on whether the trajectory direction is clockwise or counter-clockwise.
\begin{table}[htb!]
\centering
\begin{tabular}{|c|l|l|}
\hline
      Numerical & $\rbar(T)\sin\gbar(T)$  & $\rbar(T)\sin\gbar(T)$\\
			scheme    & figure \ref{fig:strob_sec}(b) & figure \ref{fig:strob_sec}(c) \\
	\hline
	\hline
 SS    & 1.60001223 & 1.64997342 \\
	 PS    & 1.60001462    & 1.64997099 \\
	 POS   & 1.60004979    & 1.65004149 \\
   PS ($\what\equiv0$) & 1.59984641 & 1.64969150 \\
\hline
	\end{tabular}
\caption{Value of $\rbar(T)\sin\gbar(T)$ for the 4 numerical results plotted in figure \ref{fig:strob_sec}(b) and figure \ref{fig:strob_sec}(c).}
	\label{tab:orbit_compare}
	\end{table} 

There are four separate trajectories in the panels in figure \ref{fig:strob_sec}, namely results of the SS (black), the PS (red) the POS (blue), and the PS with $\what\equiv0$ (green), with the markers in figure \ref{fig:strob_sec}(a) denoting the point $(\rbar(T),\gbar(T))$. The results show that for $\rbar(0)=1.6$ (panels (a,b)) and $\rbar(0)=1.65$ (panel (c)) the SS and PS are in excellent quantitative agreement, and are indistinguishable from one another. See table \ref{tab:orbit_compare} for a comparison of the value of $\rbar(T)\sin\gbar(T)$ for each method. Also, when examining panels (b) and (c) we also see that these results have the same dynamics, namely that in panel (b) $\rbar(T)-\rbar(0)>0$, and in (c) $\rbar(T)-\rbar(0)<0$. The significance of this will be discussed shortly. The results of the POS (red) are also in very good agreement with the SS and PS results, and in panel (a) the dynamics are also the same. In panel (b) the POS dynamics appear different to the SS and PS results, but as we discuss in \S\ref{sec:Results}, this is a region where the POS and PS begin to differ in behaviour. In figure \ref{fig:strob_sec}, the green trajectory, given by the PS with $\what\equiv0$, has been included in order to highlight the significance of including the remainder function $\what$ in the numerical method and not just the other $O(\epsilon)$ terms in the matrix $\mathcal{C}$ in \eqref{3.ODEint}. An approach neglecting $\what$ would be desirable, because it would mean the PDE \eqref{leadingw} does not need to be solved at each time-step, but clearly this approach not only produces qualitatively inaccurate results in the first cell, but also incorrect trajectory dynamics. Note that all the trajectories displayed in figure \ref{fig:strob_sec} are robust with regard to space step $\Delta x$, and changing this will not change the picture qualitatively. {\ctext Note, we also performed a calculation where we solved \eqref{3.PDE1} fully, as was the approach in \cite{rossides2014}, but did not include this result in figure \ref{fig:strob_sec}. We found that the results agreed exactly with the SS, as we would expect, and the run time was of the same order of magnitude as the SS runs. Hence any advantage the PS and the POS have over the SS, are the same advantage as they have over the approach taken in \cite{rossides2014}.}

By examining phase portrait trajectories using these Poincare sections, we can deduce information about the dynamics of the system by comparing $\left(\rbar(0),\gbar(0)\right)$ with the {\ctext first return point} $\left(\rbar(T),\gbar(T)=\gbar(0)\right)$. To make this comparison easier we introduce the function $\Pi(\rbar(0)) = \rbar(T)-\rbar(0)$, which we will use to classify the dynamics of the two pulse interaction. There are three possible scenarios to consider:
\begin{enumerate}
  \item $\Pi(\rbar(0))> 0$ : unstable dynamics, the trajectory is repelled from the equilibrium point $(\rbar_{\eq},\gbar_{\eq})$,
  \item $\Pi(\rbar(0))= 0$ : neutral dynamics, existence of a periodic orbit/limit cycle,
  \item $\Pi(\rbar(0))< 0$ : stable dynamics, the trajectory is attracted towards the equilibrium point $(\rbar_{\eq},\gbar_{\eq})$,
\end{enumerate}
where we choose $\rbar(0)>\rbar_{\eq}$ and $(\rbar_{\eq},\gbar_{\eq})$ are the coordinates of the equilibrium points labeled $A,~B,~C~$ in figure \ref{fig:Sergeyppred}. 

By considering the Poincare section in figure \ref{fig:strob_sec} we note that $\Pi(1.6)>0$ highlighting that these trajectories produce unstable dynamics, while $\Pi(1.65)<0$ and so this initial condition leads to stable dynamics. This suggests that for $1.6<\rbar(0)<1.65$ there exists an initial condition which produces a periodic orbit. This will be examined in detail in \S\ref{sec:results_two}.

\begin{figure}[htb!]
     \begin{center}
(a){\includegraphics[width=0.45\textwidth]{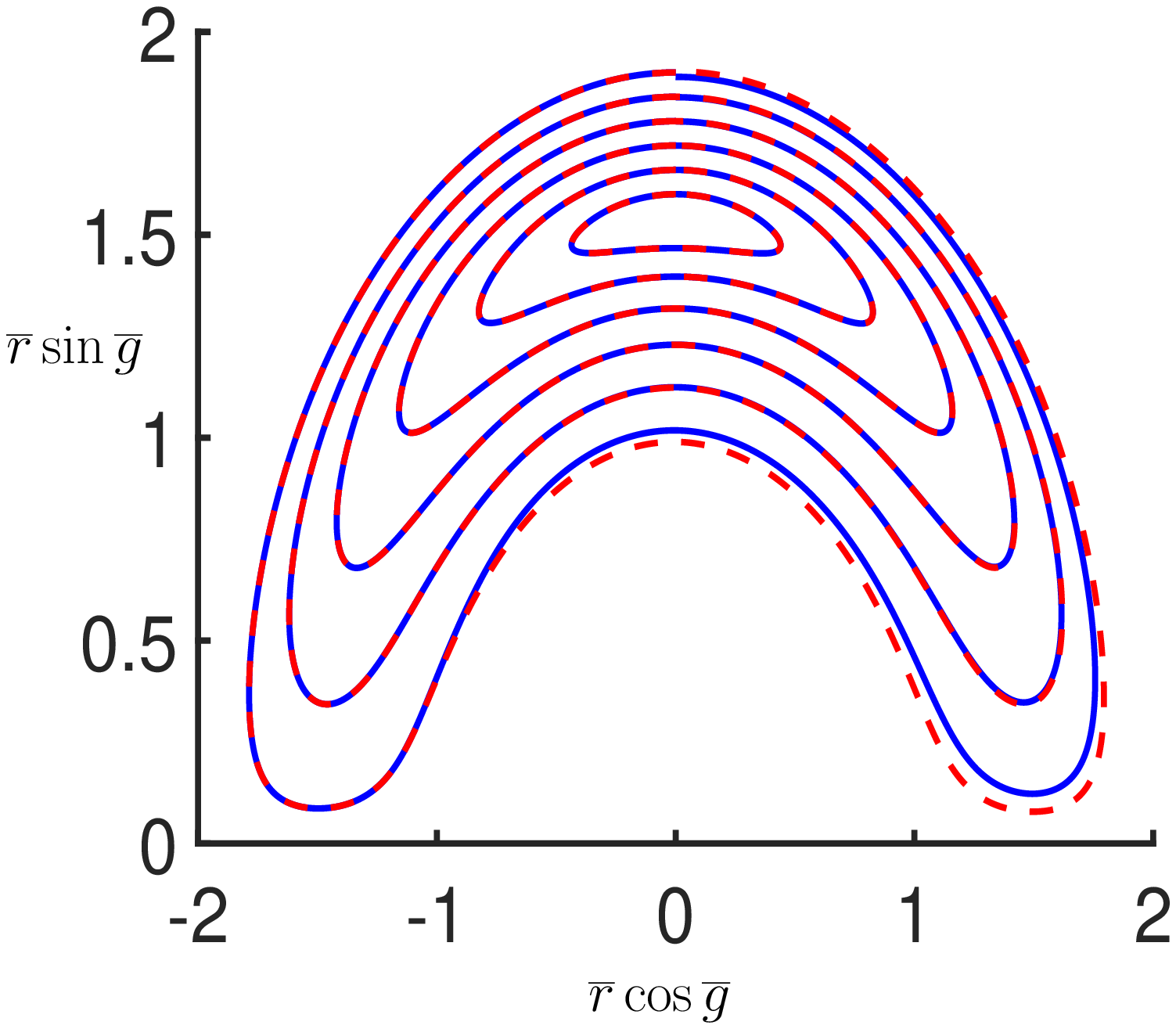}}
(b){\includegraphics[width=0.45\textwidth]{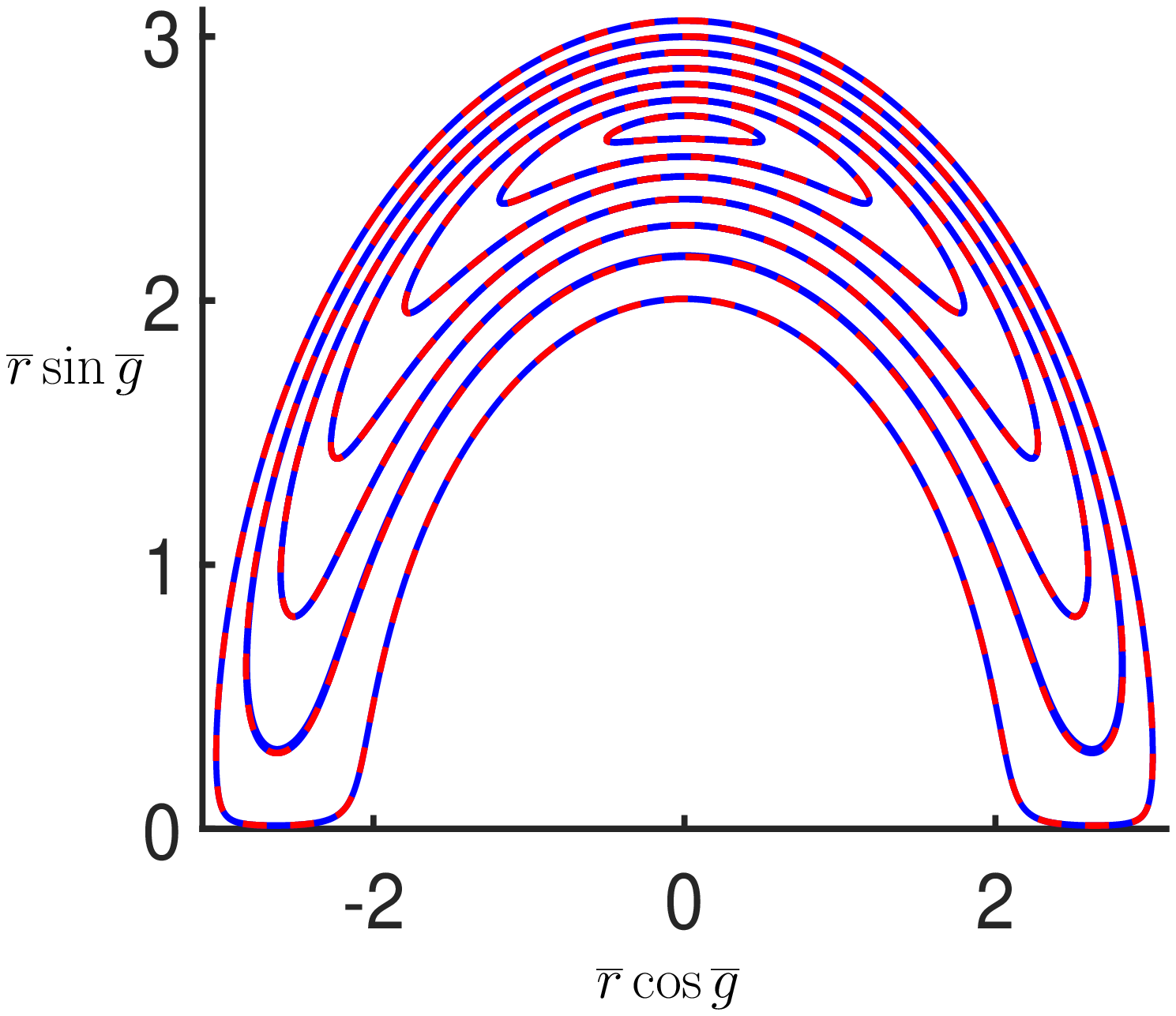}}
    \end{center}
    \caption{(colour online) Numerical Comparison of the Projected Scheme (PS) (solid red line) and the Projection ODE Scheme (POS) (dashed blue line) for the parameters (\ref{eqn:parameters}). Panel (a) plots the dynamics in cell 1, while (b) plots the dynamics in cell 2.}
   \label{fig:case1_test}
\end{figure}   
In figures \ref{fig:case1_test}(a) and \ref{fig:case1_test}(b) we consider a selection of solution trajectories for various initial points $(\rbar(0),\pi/2)$ for both the PS and POS in cell 1 and cell 2 respectively. The results here show that the POS provides an excellent approximation to the PS in cell 2, as well as in the majority of cell 1. On trajectories where $\rbar(t)\lesssim1.5$ the POS and PS agreement begins to diverge and produce different dynamics in the system. This becomes apparent when we consider plots of $\Pi(\rbar)$ in \S\ref{sec:results_two}, but as long as $\rbar(t)\gtrsim1.5$ during the simulation then the POS offers a suitable method to simulate the system.

\subsection{Computational time advantages}

Here we investigate a significant reason why the projected scheme (PS) and projected ODE scheme (POS) are advantageous in investigating the dynamics of well-separated pulses.
\begin{table}[htb!]
\centering
\begin{tabular}{|c|c|c|c|c|}
\hline
      Scheme & $\rbar(0)$ & $\gbar(0)$ & T & Comp. Time(s)\\
	\hline
	\hline
   SS    & 1.60 & $\dfrac{\pi}{2}$   & 500  & 3.978e+05 \\
	 PS    &     &         &       & 9.6419e+02 \\
	 POS   &     &         &       & 25.7409 \\
	\hline
	 SS    & 1.65 & $\dfrac{\pi}{2}$  & 500  & 3.235e+05 \\
	 PS    &      &        &       & 1.2387e+03 \\
	 POS   &      &        &       & 16.2092 \\
	\hline
	 SS    & 1.70 & $\dfrac{\pi}{2}$  & 500  & 3.0859e+05 \\
	 PS    &      &        &       & 1.9991e+03 \\
	 POS   &      &        &       & 15.9248 \\
	\hline
	\end{tabular}
\caption{Computational wall time comparison between SS, PS and POS for two pulse interaction in the QCGLE for 3 initial conditions. Here $\Delta x = 10^{-2}$ in both the SS and PS and the system parameters are given by \eqref{eqn:parameters}.}
	\label{tab:CGL}
	\end{table} 
In table \ref{tab:CGL} we examine how the three schemes SS, PS and POS compare with regard to computational wall time for interacting pulses. In all cases, the solution is evolved for the interval $t = [0,500]$ and we consider the initial phase difference between two pulses to be $\gbar(0)=\frac{\pi}{2}$. As can be seen from table \ref{tab:CGL}, the PS and POS are significantly faster than the SS, with the PS being more than 150 times faster, and the POS being about 20000 times faster. Varying the initial distance between the two pulses $r(0)$, we observe that the PS and POS excel in cases where $\rbar(0)$ is large. Hence, it is clear that splitting the solution of the QCGLE into a sum of pulses and a remainder term \eqref{0.mult}, and solving a boundary-value problem for $w$ allows the PS to take larger time steps than the SS and still maintain maximum accuracy for weakly-interacting pulses. The speed of the POS is useful for well-separated pulses as it allows us to readily examine wide areas of parameter space in order to identify the dynamics of the system and any new phenomena that arise. This would not be possible using the SS, and even using the PS this is very time consuming and unwieldy.

\section{Results}\label{sec:Results}

\subsection{Two-pulse interaction in QCGLE}
\label{sec:results_two}

In this section we examine the dynamics of the two-pulse interaction for the QCGLE using the PS and the POS. The projected system has been derived for the two-pulse solution of the form
\[
 u(t,x) = V_1 + V_2 + w(t,x),
\]
in \cite{Turaev07}, who considers only the leading order terms, with remainder function $w(t,x) \equiv 0$. Even for this simple two-pulse case, it is unknown as to what kind of dynamics appear, and so we examine part of parameter space in order to determine an answer to this question. Hence we provide a full schematic phase portrait for the two-pulse interaction in QCGLE for a range of system parameters. 

\subsubsection{Two-pulse interaction dynamics: $\beta_r=\Re(\beta)=-0.02$}\label{sec:res1}

In this section we use the same parameter values as in \eqref{eqn:parameters} to investigate the dynamics of the two-pulse solutions in the QCGLE. {\ctext Note that, according to theorem 2.1, the centre-manifold reduction exists for all $\beta_r<0$ (if our assumptions on the initial pulse are satisfied) and likely does not exist if $\beta_r\geq 0$ (since assumption \eqref{2.stable} is not satisfied in this case). However, we do not a priori know at which pulse separation distance  this reduction starts to work. The case where it fails corresponds to the strong interaction of pulses and the dynamics there may be more complicated, including pulse collision etc. In other words, for fixed values of parameters, we do not know what analytical quantity we should compute to guarantee this fact. Instead, we do computations utilizing the formal asymptotic expansions and related projected schemes and compare the obtained results with the standard scheme involving direct numerical simulation of the PDE. After that we a posteriori conclude that we are in the regime of a weak pulse interaction and the reduction indeed works.} 

\begin{figure}[ht!]
\centering
(a) \includegraphics[width=0.5\textwidth]{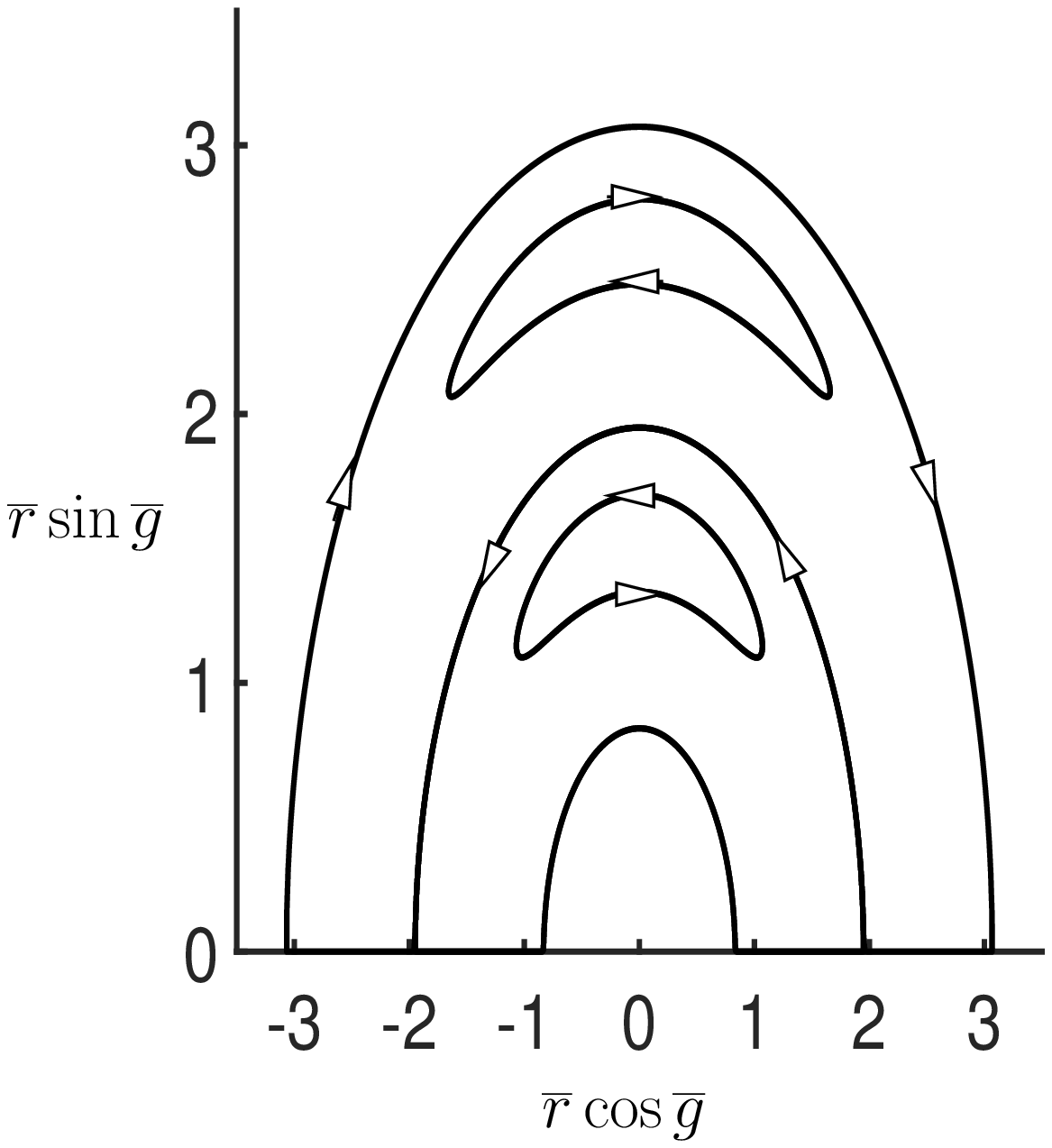}
\\(b) \includegraphics[width=0.45\textwidth,height=0.35\textwidth]{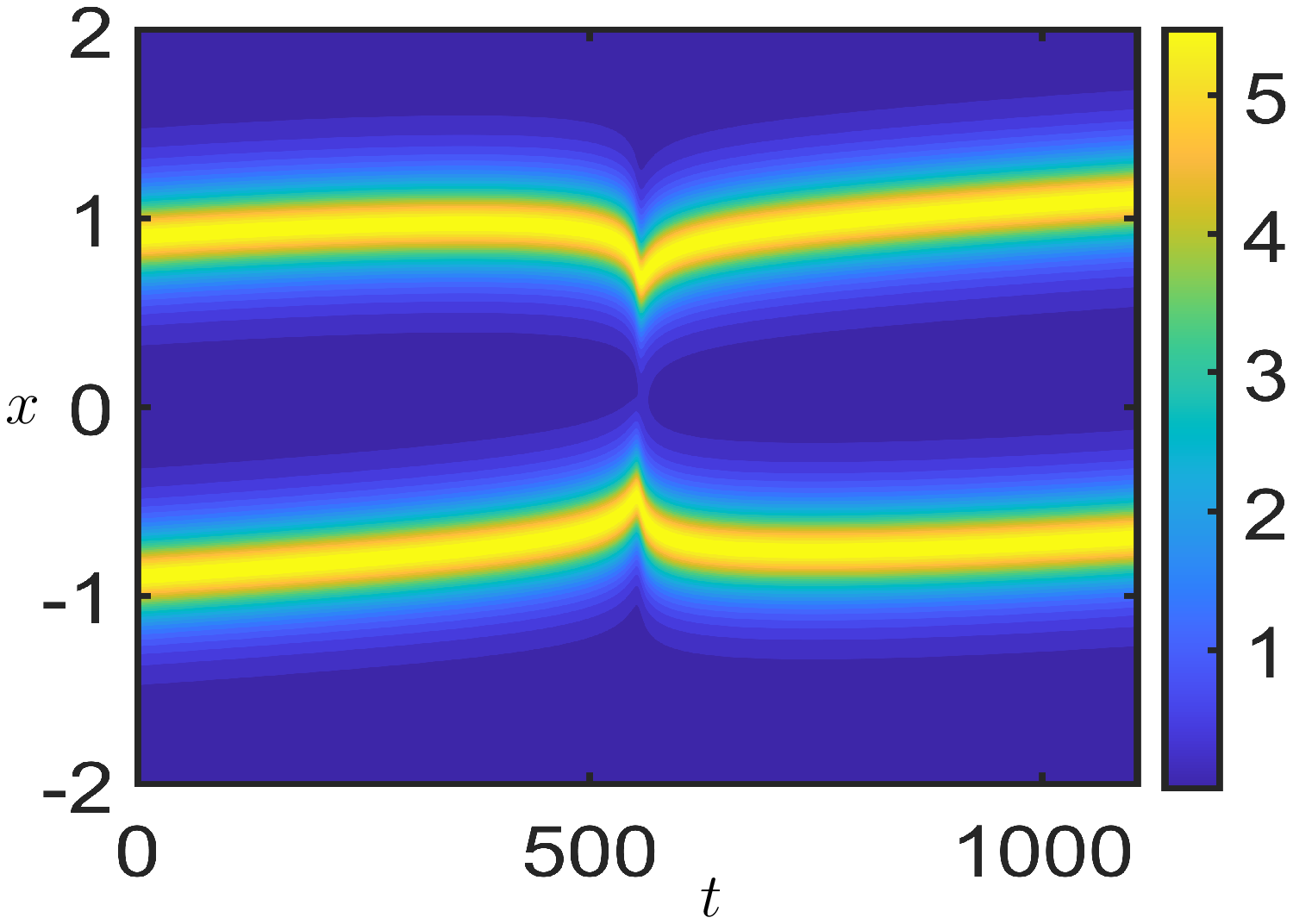}
(c) \includegraphics[width=0.45\textwidth,height=0.35\textwidth]{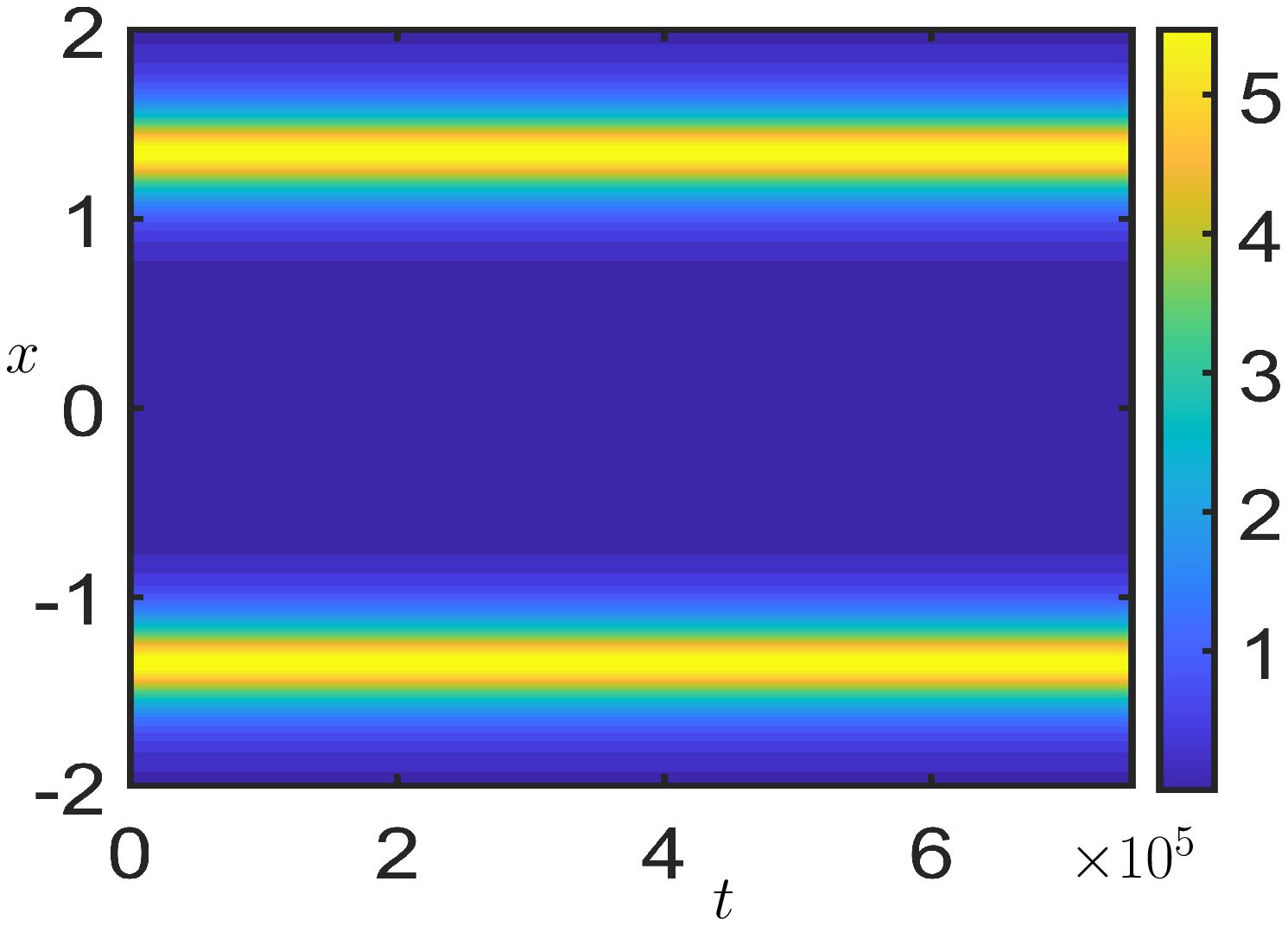}
\caption{(color online) (a) Phase portrait for two-pulse interaction in QCGLE, including trajectories for weakly (WI) and very-weakly-interacting (VWI) pulses. (b) WI: Profile of weakly-interacting pulses located in cell 1 with initial condition $r_2(0)-r_1(0) = 1.8,\;g_2(0)-g_1(0)=\frac{\pi}{2}$, \textbf{(c)} VWI: Cell 2 very-weakly-interacting dynamics with initial condition $r_2(0)-r_1(0) = 2.7,\; g_2(0)-g_1(0)=\frac{\pi}{2}$.}
\label{fig:two_cell_traj}
\end{figure}

We begin by showcasing some typical trajectories in cell 1 (which we call weakly-interacting (WI)) and cell 2 (which we call very-weakly-interacting (VWI)), see figure \ref{fig:two_cell_traj}(a). 
 Note that although the trajectories appear to be closed, numerically we observe that this is not the case. In figures \ref{fig:two_cell_traj}(b) \& (c) we plot $|u(t,x)|$ for both the WI and VWI solution trajectories in figure \ref{fig:two_cell_traj}(a). In both cases the pulses have an approximately periodic movement and the weaker interaction in the cell 2 results in much larger periods for very-well-separated pulses. Moreover, in figure \ref{fig:two_cell_traj} we see that depending on the initial distance between pulses $|\rbar(0)|$, the pulses can travel in space in different directions. In particular, for very-well-separated pulses, that belong in cell 2, pulses travel to the left, whereas for weakly-interacting pulses, that belong in cell 1, pulses travel to the right. In order to investigate the dynamics within each of these cells further, we will use the stroboscopic/Poincare section as introduced in \S\ref{subs:justification}, and in particular we examine the behaviour of the function $\Pi(\rbar(0))=\Pi(\rbar)$. 
  
\vspace*{0.5cm}
\textbf{Cell 1 stability}
\vspace*{0.5cm}

\begin{figure}[htb!]
  \centering
\hspace*{-1cm}
 \includegraphics[width=4in]{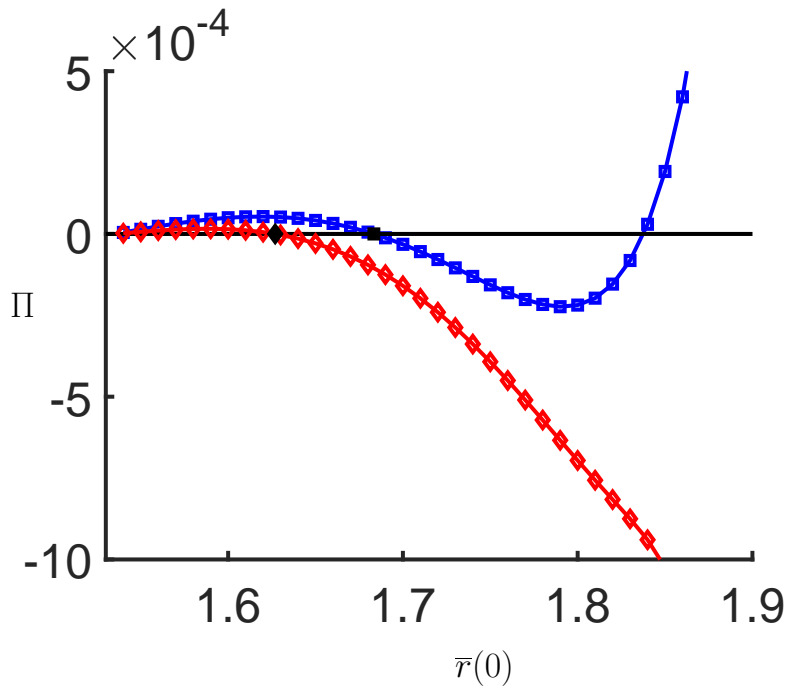}
\caption{(colour online) The function $\Pi(\rbar):=\rbar(T)-\rbar(0)$ as a function of the initial difference $\rbar(0)$. The red diamonds represent results of the PS while the blue squares are results of the POS. In each case the black symbols represent the case $\Pi(\rbar)=0$ which denotes a periodic orbit.}
\label{fig:case1_Poin}
  \end{figure}

For the parameter values \eqref{eqn:parameters} the position of the equilibrium point labeled $A$ in figure \ref{fig:Sergeyppred} is $(\rbar_{\eq},\gbar_{\eq})=(1.5353,\pi/2)$. We fix $\gbar(0) = \pi/2$ and vary values of $\rbar(0)>\rbar_{\eq}$ in order to calculate $\Pi(\rbar(0))$ and produce the results in figure \ref{fig:case1_Poin} showing the system dynamics in the cell 1.
 In this figure the red diamonds represent results of the PS while the blue squares represent results of the POS. We observe that for trajectories close to $\rbar_{\eq}$, $\Pi(\rbar(0)) > 0$. This demonstrates that these trajectories are spiraling out from the equilibrium point, albeit very slowly due to the small magnitude of $\Pi$. As we increase $\rbar(0)$ we observe that $\Pi$ passes through zero at $\rbar(0)=1.6269$ (denoted by the black diamond) and is then negative for the remainder of the $\rbar(0)$ values considered. For $\rbar(0)>1.6269$ this indicates that the trajectories are spiraling towards the equilibrium point, thus at $\rbar(0)=1.6269$ we must have a stable periodic orbit. Note that the periodic orbit position, denoted by the black diamond, is calculated using the bisection method to within a tolerance of $10^{-10}$. The blue squares of the POS are in qualitative agreement with the {\ctext PS} for $\rbar_{\eq}\leq\rbar(0)\leq1.8$ (here the periodic orbit is found at $\rbar(0)=1.6834$), but for $\rbar(0)>1.8$ the POS suggests there is another periodic orbit. However, in this region the pulses move towards each other during their motion and have value of $\rbar(t)<1.5$ for some time values (see figure \ref{fig:two_cell_traj}(b)). Hence, the asymptotic forms of the inner products given in Appendix \ref{appen:innerprods} are no longer valid, and thus the POS no longer gives qualitatively correct dynamics. Note that the results displayed in figure \ref{fig:case1_Poin} are robust in regard to space step $\Delta x$.

As the dynamics of the PS predicts that the trajectories spiral outwards from the equilibrium point for $\rbar(0)>1.6269$, a natural question to ask is, what happens to the cell boundary between cells 1 and 2? This boundary, as seen in figure \ref{fig:Sergeyppred}(b), consists of two heteroclinic orbits connecting two saddles, but the inclusion of the higher order terms are expected to destroy the heteroclinic orbits and cause the orbit to converge into either cell 1 or 2.
\begin{figure}[htb!]
     \begin{center}
(a)\includegraphics[width=0.45\textwidth]{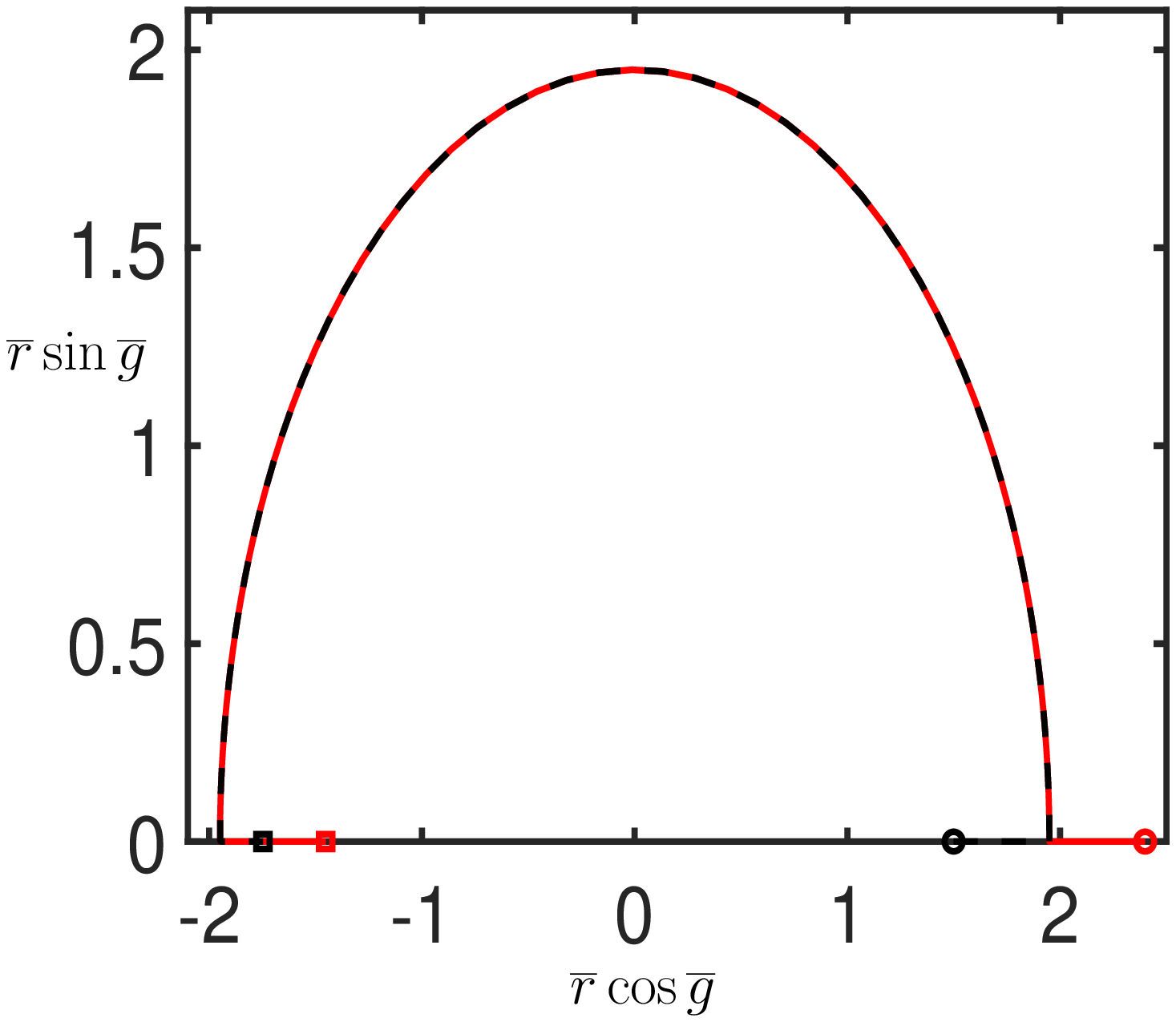}
(b)\includegraphics[width=0.45\textwidth]{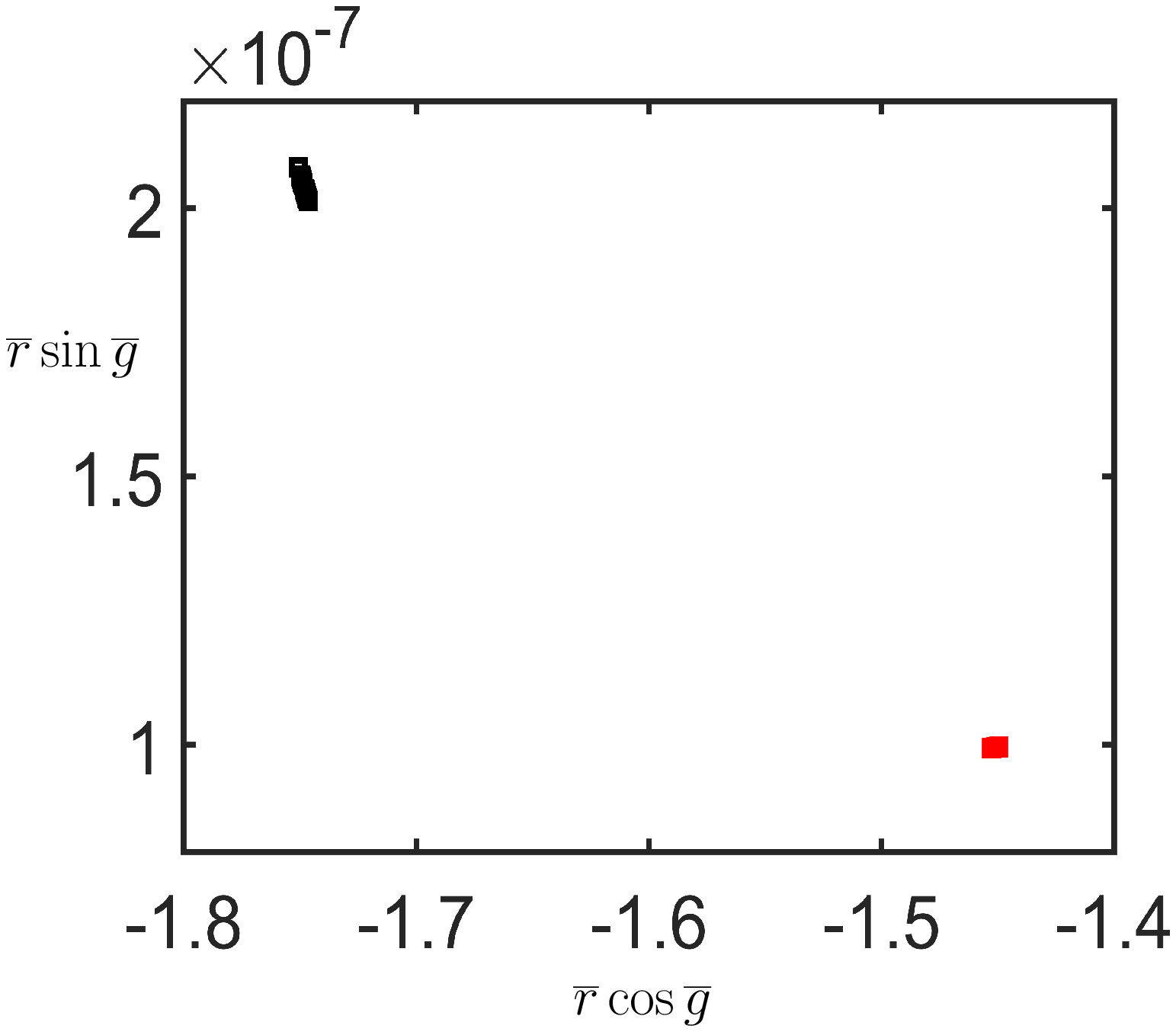}
(c)\includegraphics[width=0.45\textwidth]{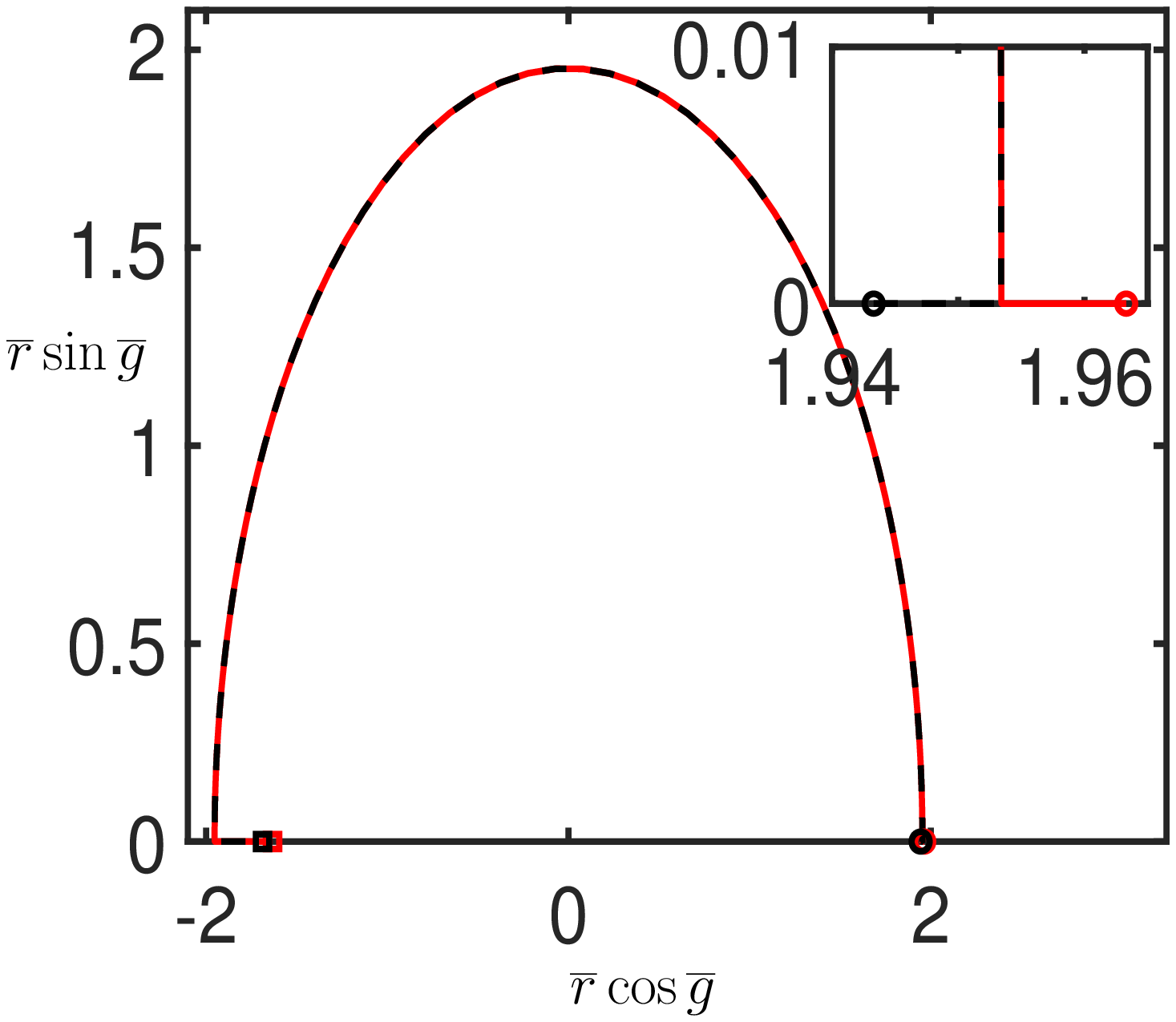}
        \end{center}
\caption{(colour online) (a) Trajectories which pass close to the cell 1/cell 2 saddle points with the solid red line originating from cell 2 and the black dashed line originating from cell 1, calculated via the POS. The circles give their starting positions and the squares give the computation end positions. (b) The end position of the ensemble of $2^{20}$ parameters in the POS, all of which are located in cell 1. (c) The equivalent to panel (a) except calculated via the PS. The insert gives a close-up of the starting point of the trajectories.}
\label{fig:1cellsplit}
\end{figure}
In figure \ref{fig:1cellsplit}(a) we plot two trajectories which pass close to the saddle point labeled $2'$ in figure \ref{fig:Sergeyppred}(b). The initial points of the two trajectories are $(\rbar_{\eq}\pm0.45,10^{-15})$, where $\rbar_{\eq}=1.9533$, i.e. the initial points lie close to the $x$-axis with one (red solid line) in cell 2 and the other (black dashed line) in cell 1. As these trajectories pass close to two saddle points, the dynamics of the system slows and integration times increase, hence we present these results for the POS in order to maintain a high level of numerical accuracy. The results in figure \ref{fig:1cellsplit}(a) show that the two trajectories starting in cells 1 and 2 both have end points (square symbols) in cell 1, this demonstrates that the two heteroclinic orbits in the leading order problem actually split, and on the edge of this cell the trajectories move into cell 1. Once the trajectory from cell 2 enters cell 1, then figure \ref{fig:case1_Poin} shows that it will spiral in and converge to the periodic orbit in that cell. 

To investigate the robustness of the results in figure \ref{fig:1cellsplit}(a), in particular the trajectory which starts in cell 2 and ends in cell 1, we perform an ensemble calculation of the 20 $B_i$ and $\mu_i$ coefficients noted in Appendix \ref{appen:innerprods}. The take home result here is that in every one of the $2^{20}$ simulations, both trajectories starting in cells 1 and 2 of panel (a) end up in cell 1. Hence the splitting of the cell boundaries is a robust phenomena. For completeness, in figure \ref{fig:1cellsplit}(c) we consider the same case, except this time calculated by the PS. Note, we begin much closer to the saddle point, but the qualitative dynamics are still identical to the POS result, with both trajectories ending in cell 1.

\vspace*{0.5cm}
\textbf{Cell 2 and 3 stability}
\vspace*{0.5cm}

We now investigate the second equilibrium point of the phase portrait, labeled $B$ in figure \ref{fig:Sergeyppred}(b), using the same Poincare section method as we used to investigate the stability of the cell 1 dynamics. This case is considered as a very-weakly-interacting pulse case and the dynamics in this cell have not been computed up to now due to the computational expense of the calculations. A plot of what the phase-plane trajectories look like around this point can be seen in figure \ref{fig:case1_test}(b). Here it appears that all the trajectories are closed orbits, but this is not the case.

\begin{figure}[htb!]
  \centering
\hspace*{-1cm}
(a)\includegraphics[width=0.45\textwidth]{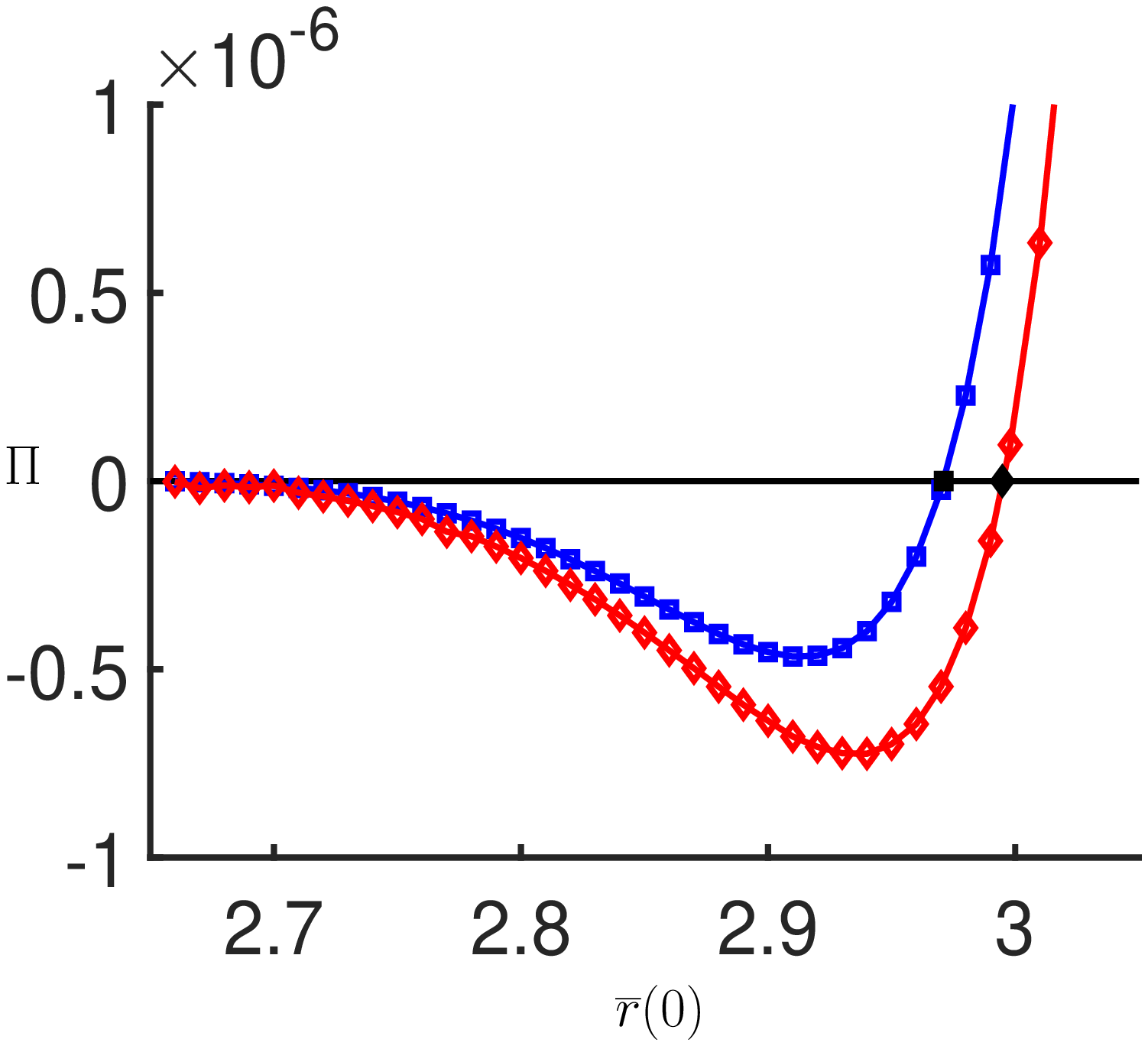}
(b)\includegraphics[width=0.45\textwidth]{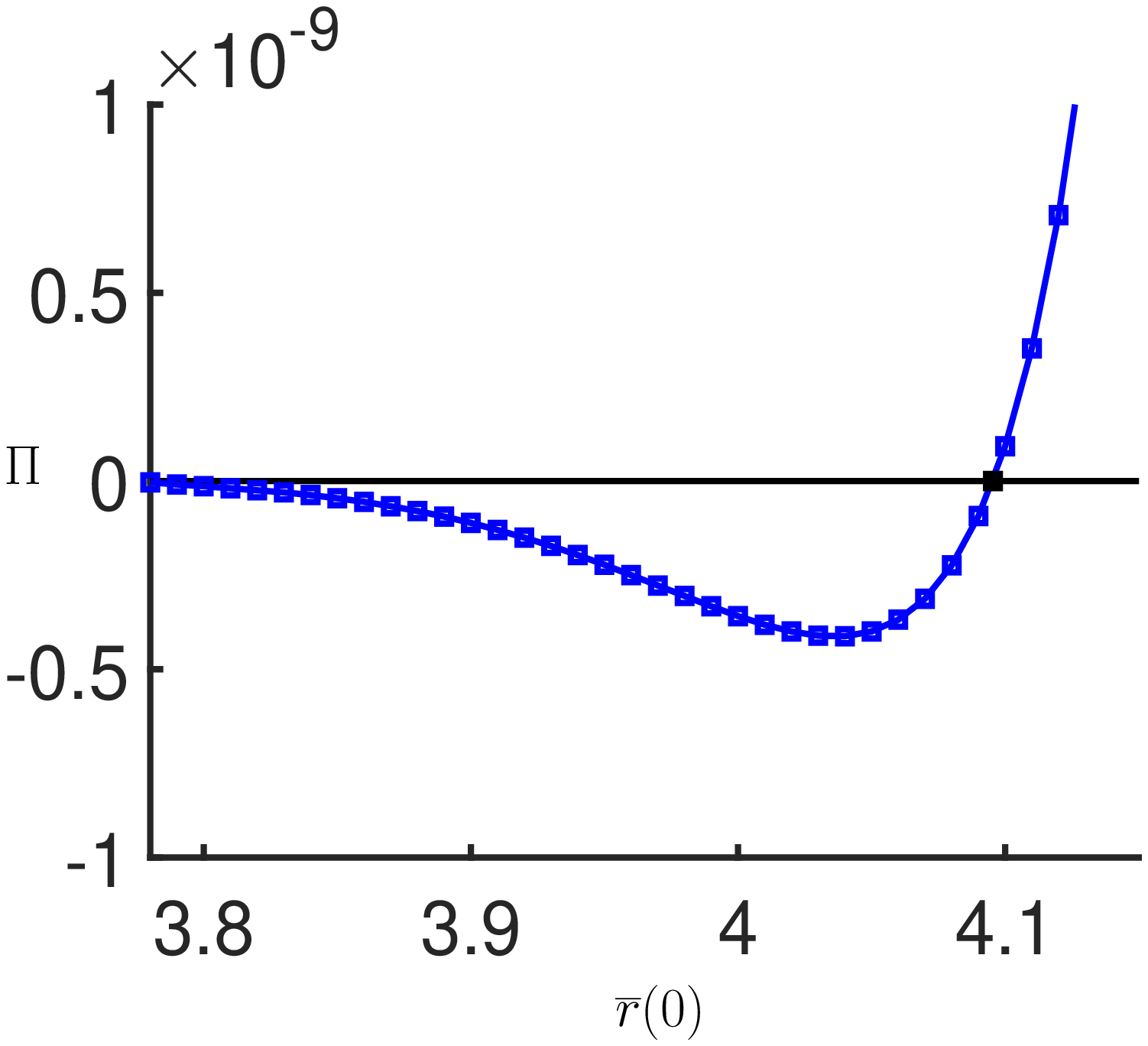}
\caption{(colour online) The function $\Pi(\rbar):=\rbar(T)-\rbar(0)$ as a function of the initial difference $\rbar(0)$ for trajectories in (a) cell 2 and (b) cell 3. The red diamonds represent results of the PS while the blue squares are results of the POS. Only results for the POS were possible in cell 3 given the very weak interaction of the pulses. In each case the black symbols represent the case $\Pi(\rbar)=0$ which denotes a periodic orbit.}
\label{fig:2cellstab}
\end{figure}
In figure \ref{fig:2cellstab}(a) we consider the function $\Pi(\rbar)$ for a range of phase-plane trajectories. These results show that, as in cell 1, we again have a periodic orbit in this cell ($\Pi(\rbar(0))=0$) at $\rbar(0)=2.9948$ for the PS, and $\rbar(0)=2.9711$ for the POS. However, the dynamics of this periodic orbit are the opposite to those seen in cell 1. For $\rbar_{\eq}<\rbar(0)<2.9948$ (for the PS), the trajectories spiral in to the equilibrium point $(2.6571,\pi/2)$, while for $\rbar(0)>2.9948$ the trajectories spiral away from the periodic orbit. Thus this periodic orbit is an unstable limit cycle. These dynamics are consistent with those we observed in cell 1. For example for $\rbar(0)>2.9948$ in cell 2, the trajectory will spiral outwards towards the edge of the cell, and then, as it passes close to the saddle point at the boundary of cell 1 and 2, as in figure \ref{fig:1cellsplit}, it will move into cell 1 and then spiral into the periodic orbit in that cell.

\begin{figure}[ht!]
     \begin{center}
(a)\includegraphics[width=0.45\textwidth]{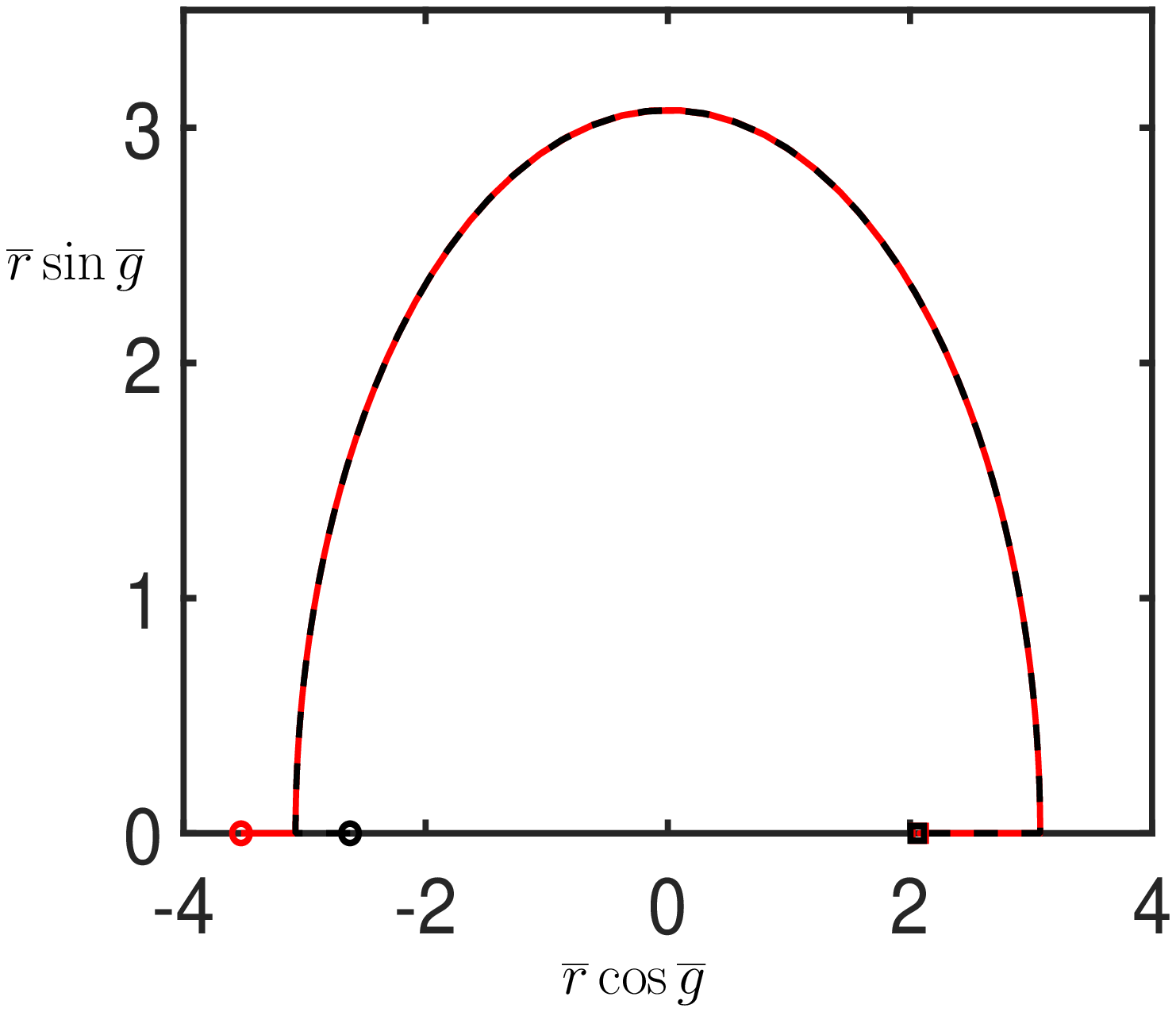}
(b)\includegraphics[width=0.45\textwidth]{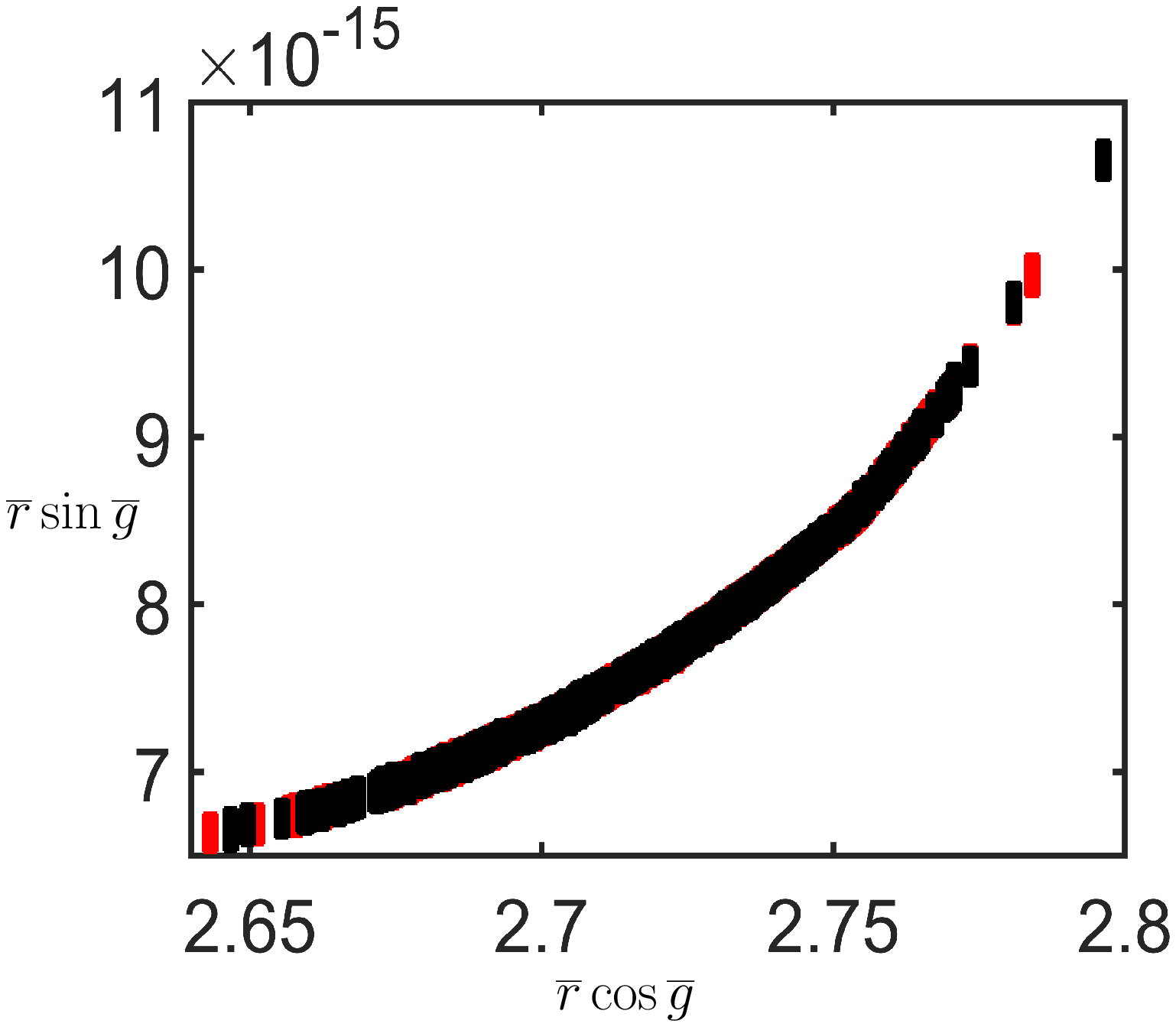}
    \end{center}
\caption{(colour online) (a) Trajectories which pass close to the cell 2/cell 3 saddle points with the solid red line originating from cell 3 and the black dashed line originating from cell 2 calculated via the POS. The circles give their starting positions and the squares give the computation end positions. (b) The end position of the ensemble of $2^{20}$ parameters in the POS, all of which are located in cell 2.}
\label{fig:2cellsplit}
\end{figure}

For cell 3, in figure \ref{fig:2cellstab}(b), the dynamics are the same as in cell 2, i.e. the equilibrium point at $(3.7789,\pi/2)$ is stable, and trajectories with $\rbar_{\eq}<\rbar(0)<4.0954$ spiral into this equilibrium point. On the other hand, trajectories with $\rbar(0)>4.0954$ spiral outwards. This outward spiraling behaviour is consistent with the dynamics in cell 2 as long as, at the cell 2 and 3 boundary, the heteroclinic orbits split such that trajectories close to the edge of cell 3 end up in cell 2. This is precisely what we show in figure \ref{fig:2cellsplit}. Again the ensemble calculation of $2^{20}$ results confirms the robustness of this calculation.

\vspace*{0.5cm}
\textbf{Full schematic projected plane for two-pulse interactions for parameter set \eqref{eqn:parameters}}
\vspace*{0.5cm}

At double numerical precision we are unable to use the POS to investigate the dynamics in cells beyond cell 3, but we can formulate a good idea of the system dynamics by considering only these 3 cells. 
\begin{figure}[ht!]
  \centering
(a)\includegraphics[width=0.6\textwidth]{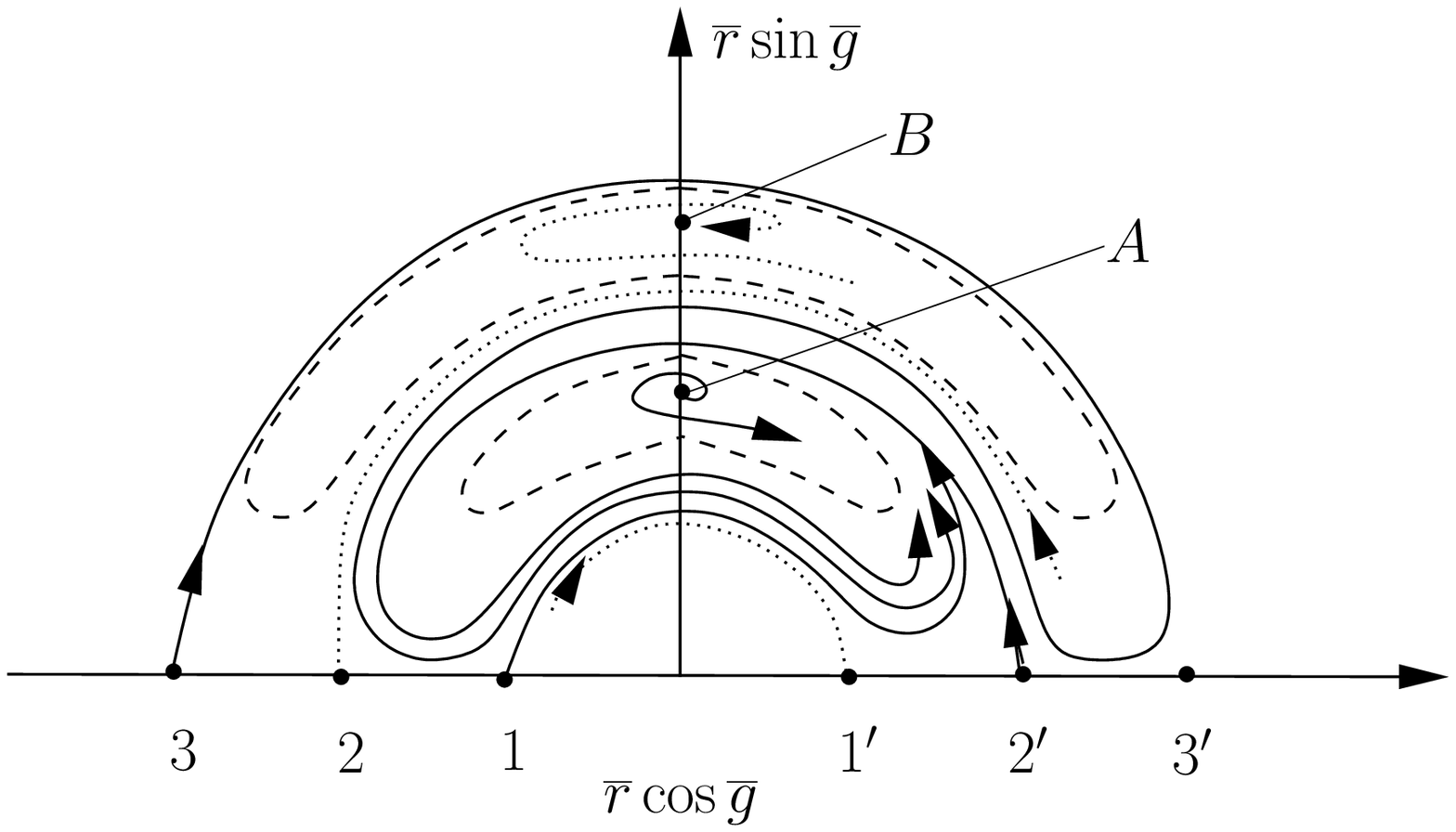}\\
(b)\includegraphics[width=0.45\textwidth]{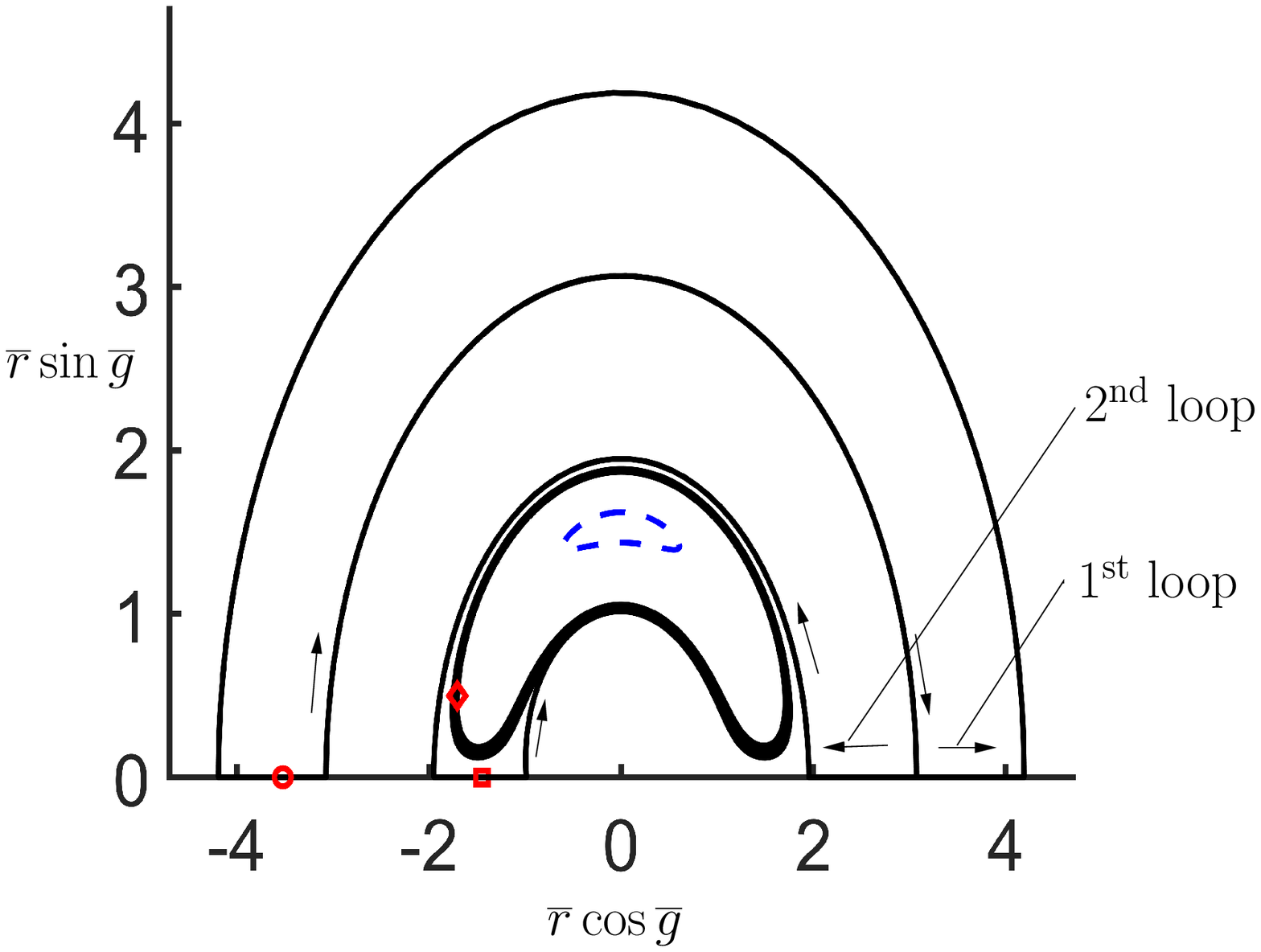}
(c)\includegraphics[width=0.45\textwidth]{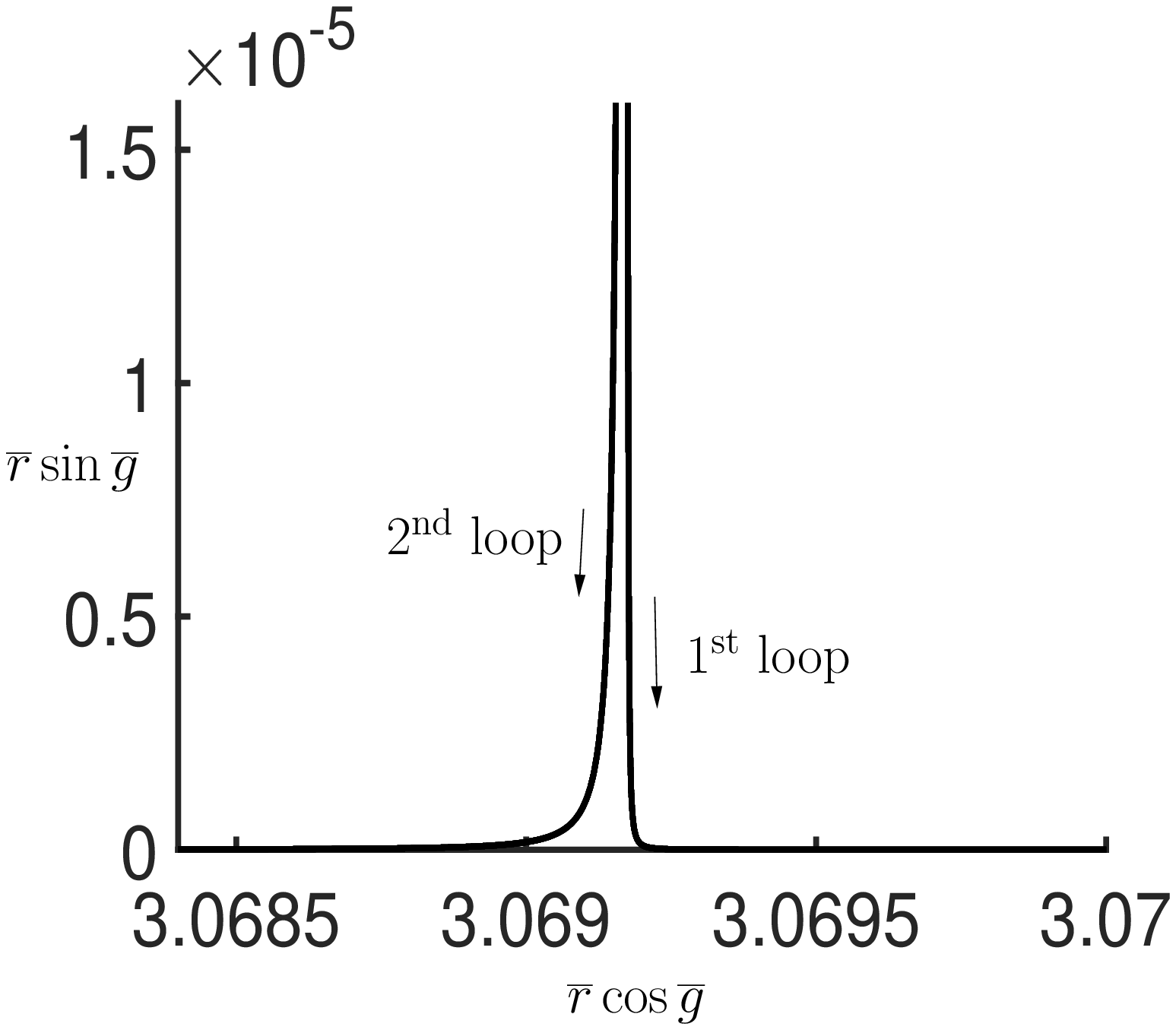}
(d)\includegraphics[width=0.45\textwidth]{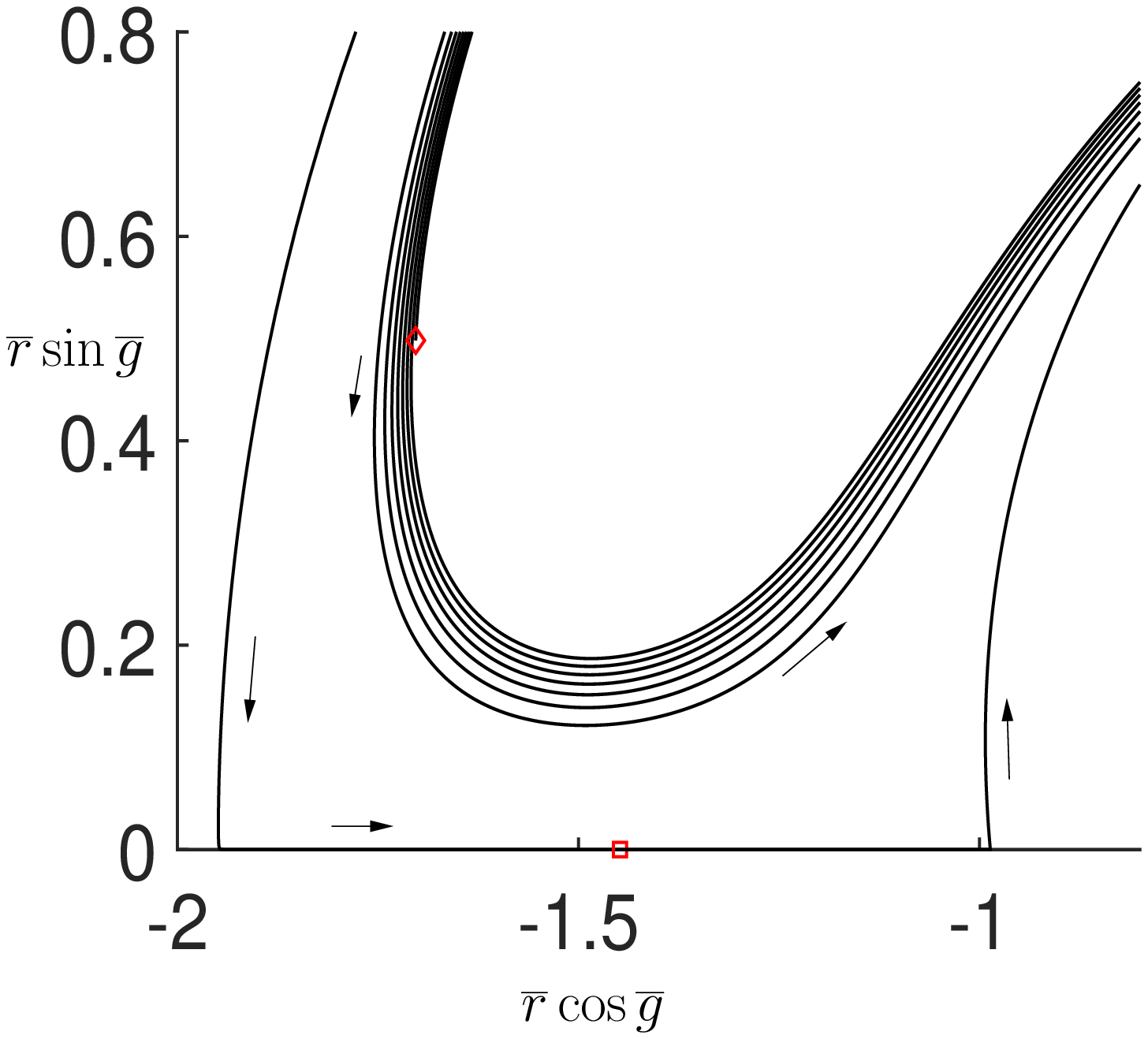}
  \caption{ (a) Schematic phase plane for two-pulse interaction in QCGLE, from very-weakly-interacting to weakly-interacting pulses. (b) A computed trajectory which starts close to the cell 3 boundary (circle). At the square symbol the POS calculation converts to a PS calculation and the trajectories spiral in towards the periodic orbit (dashed line). The diamond gives the end point of the trajectory. (c) A close up of the cell 2 and 3 boundary where the trajectory splits. (d) A close up of the trajectories in cell 1 showing the spiraling behaviour.}
 \label{fig:ppfullschem}
 \end{figure}
In figure \ref{fig:ppfullschem}(a) we give a schematic diagram of this qualitative dynamics for the two-pulse interaction in the QCGLE in cells 1 and 2. Here the dashed lines represent the two periodic orbits, the solid lines give example trajectories emanating from the equilibrium points $1$, $3$, $2'$ and $A$, while the dotted lines give example trajectories which converge at the equilibrium points $1'$, $2$ and $B$. 

In figure \ref{fig:ppfullschem}(b) we plot a true numerical trajectory for a trajectory starting at the circle symbol close to the equilibrium point $3$ in cell 3. The trajectory then moves round to $3'$, but the initial point is carefully chosen such that the trajectory does another circuit of cell 3 and, on its second visit to $3'$ the trajectory this time splits into cell 2. It then travels around to equilibrium point $2$ and into the square symbol between points $2$ and $1$. This trajectory, up to this point, is calculated using the POS as the pulses are very well separated, but at the square symbol they become close enough together that the POS is no longer valid. Hence we use this point as an initial condition in the PS and follow the trajectory into cell 1, where it spirals in towards the dashed line, representing the periodic orbit. We terminate the computation at the diamond symbol.

\subsubsection{Robustness of the two-pulse interaction dynamics}

In \S\ref{sec:res1} we investigated the dynamics of the two-pulse system for the specific set of parameters given in (\ref{eqn:parameters}). In this section we ask the question as to whether the dynamics seen in \S\ref{sec:res1} are robust and are seen for other parameter values, whether we can identify different dynamics, or whether we can identify the origin of the periodic orbits observed?

To perform the parameter study we vary the parameter $\beta_r=\Re(\beta)$ to follow the position of the periodic orbit, {\ctext keeping all other parameters in \eqref{eqn:parameters} fixed, except $\beta_i$ which also varies to ensure the dispersion relation \eqref{0.disp} is satisfied.} Note that while we only vary $\beta_r$ in this work, we observe qualitatively similar results when varying other parameters. 
\begin{figure}[ht!]
  \centering
(a)\includegraphics[width=0.45\textwidth]{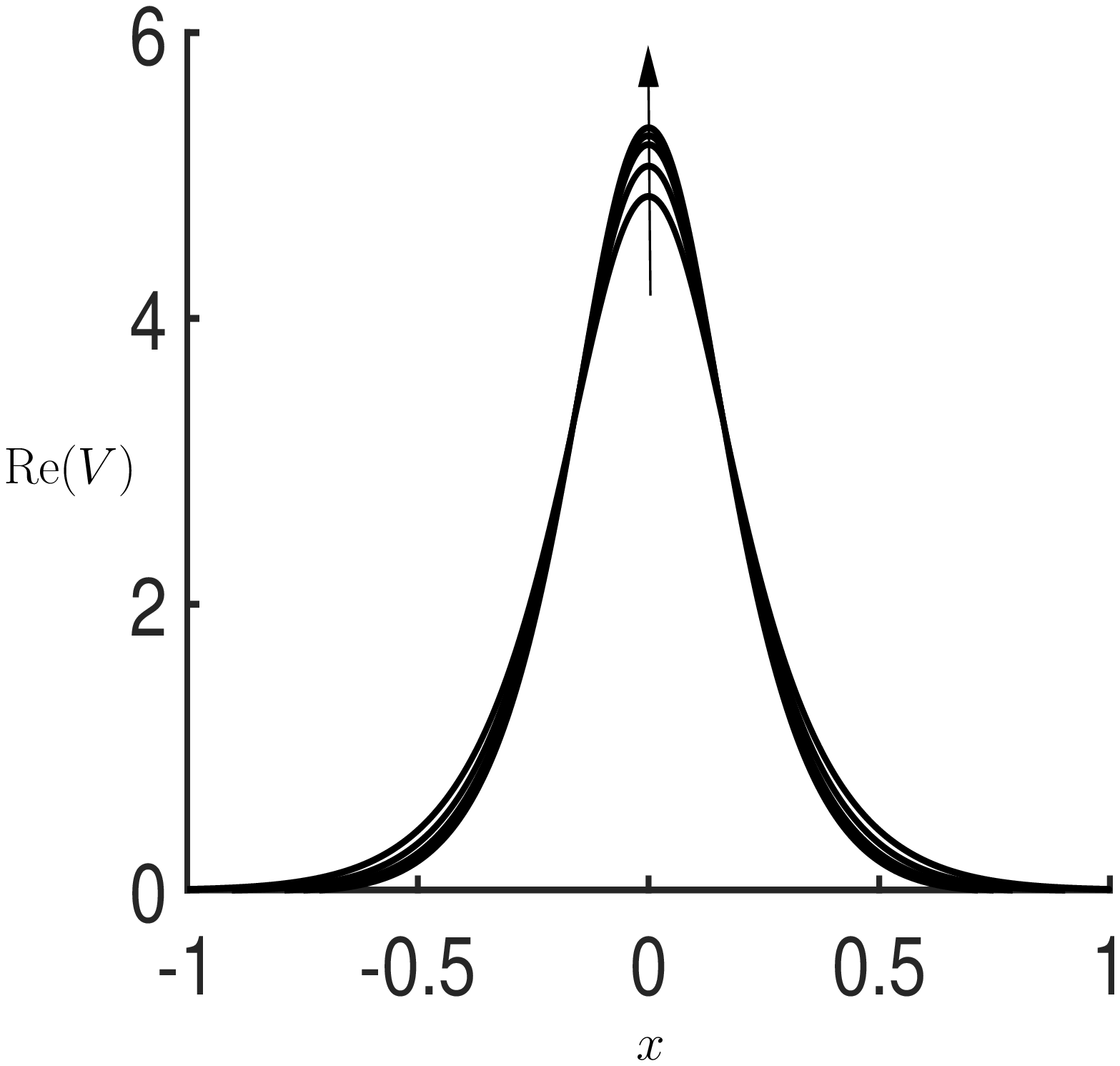}
(b)\includegraphics[width=0.45\textwidth]{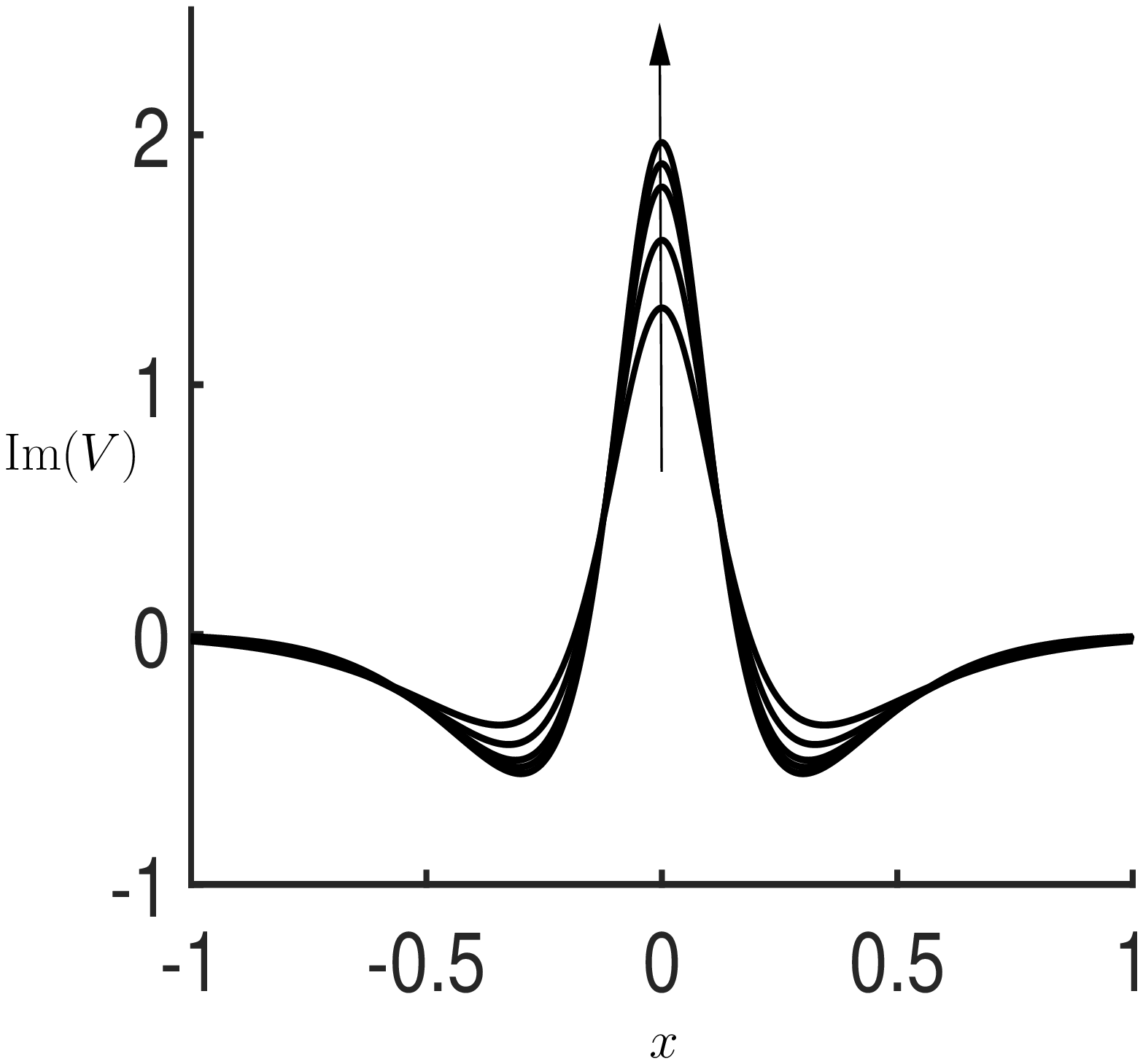}
  \caption{(a) The real part and (b) the imaginary part of the pulse $V(x)$ for $\Re(\beta)=\beta_r=-6,-3,0.02,3$ and $6$, with the arrow indicating the direction of increasing $\beta_r$.}
 \label{fig:pulse_variation}
 \end{figure}
 Figure \ref{fig:pulse_variation} shows how the real and imaginary parts of the pulse solutions $V(x)$ vary as a function of $\beta_r$. Essentially, as $\beta_r$ decreases, the pulses reduce in magnitude and become more broad.

To calculate the position of the periodic orbit we use the bisection method, and when considering the behaviour in cell 1 we use the PS, while in cells 2 and 3 we use the POS given its computational time advantage over the PS. Once the value of $\rbar(0)$ which denotes the periodic orbit is computed, the sign of the gradient $\frac{\d\Pi}{\d\rbar(0)}$ determines whether the periodic orbit is stable ($\frac{\d\Pi}{\d\rbar(0)}<0$) or unstable ($\frac{\d\Pi}{\d\rbar(0)}>0$). The dynamics of the system can then be identified for that parameter set.

\begin{figure}[ht!]
  \centering
(a)\includegraphics[width=0.6\textwidth]{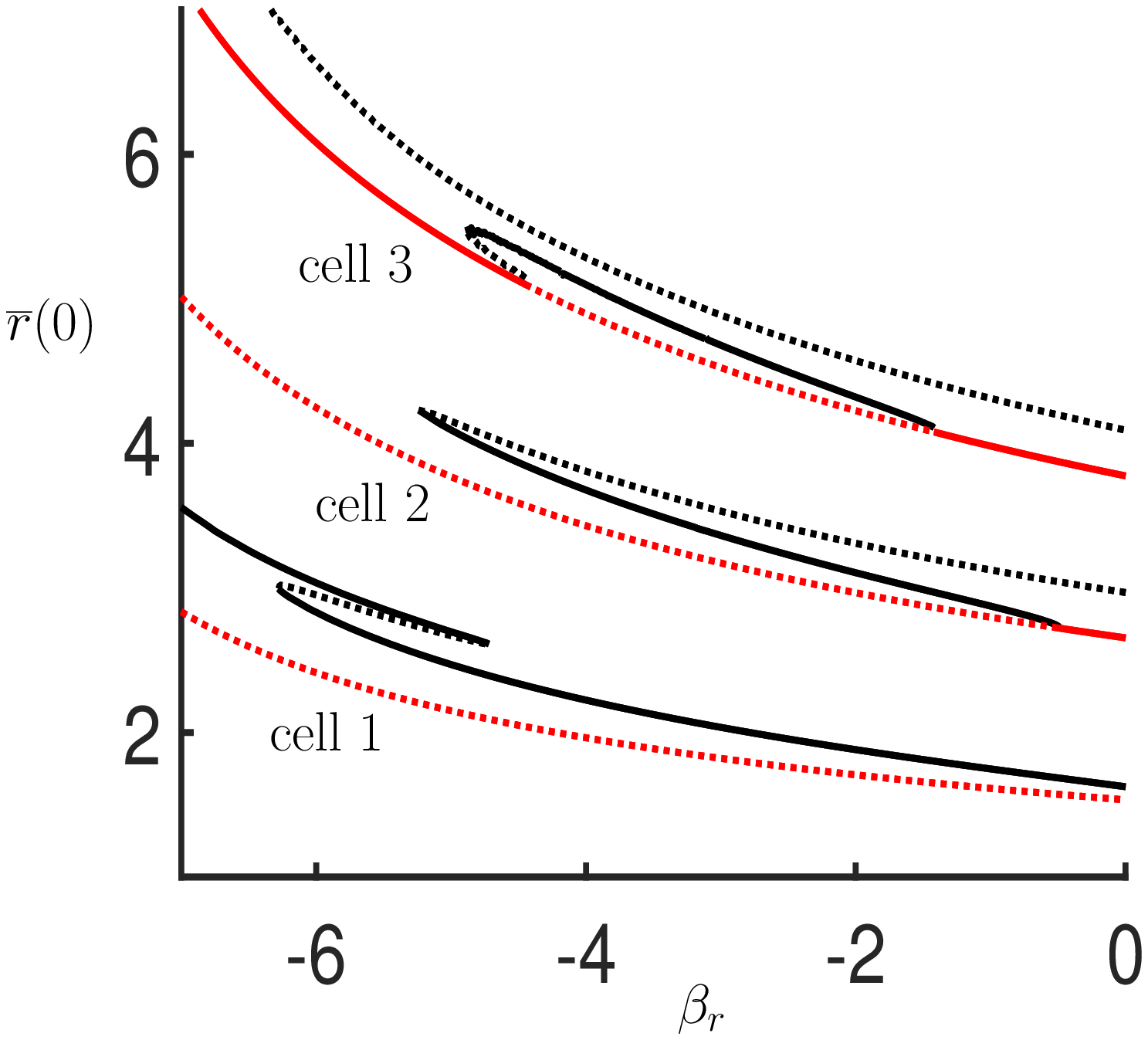}	\\
(b)\includegraphics[width=0.45\textwidth]{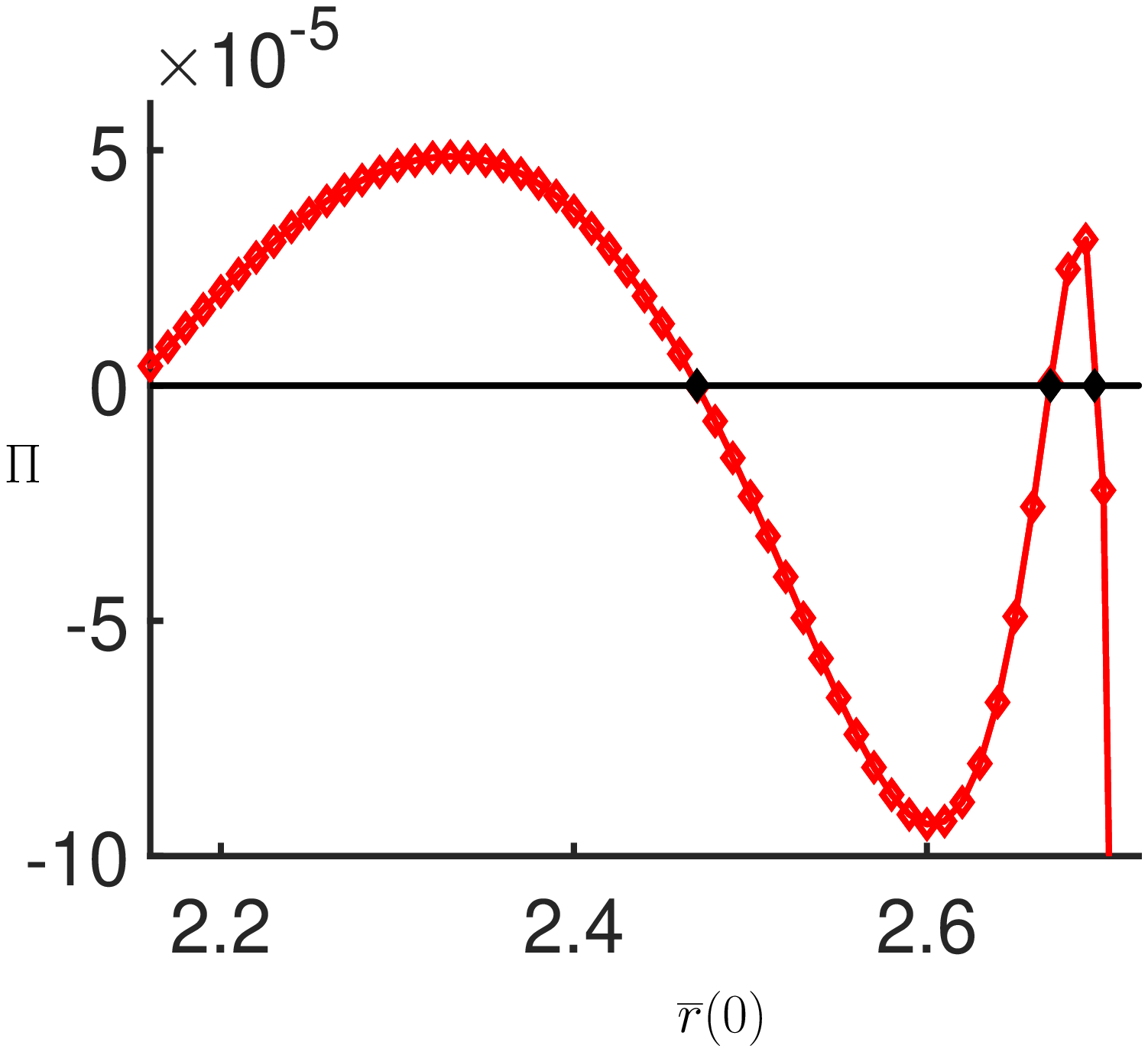}
(c)\includegraphics[width=0.45\textwidth]{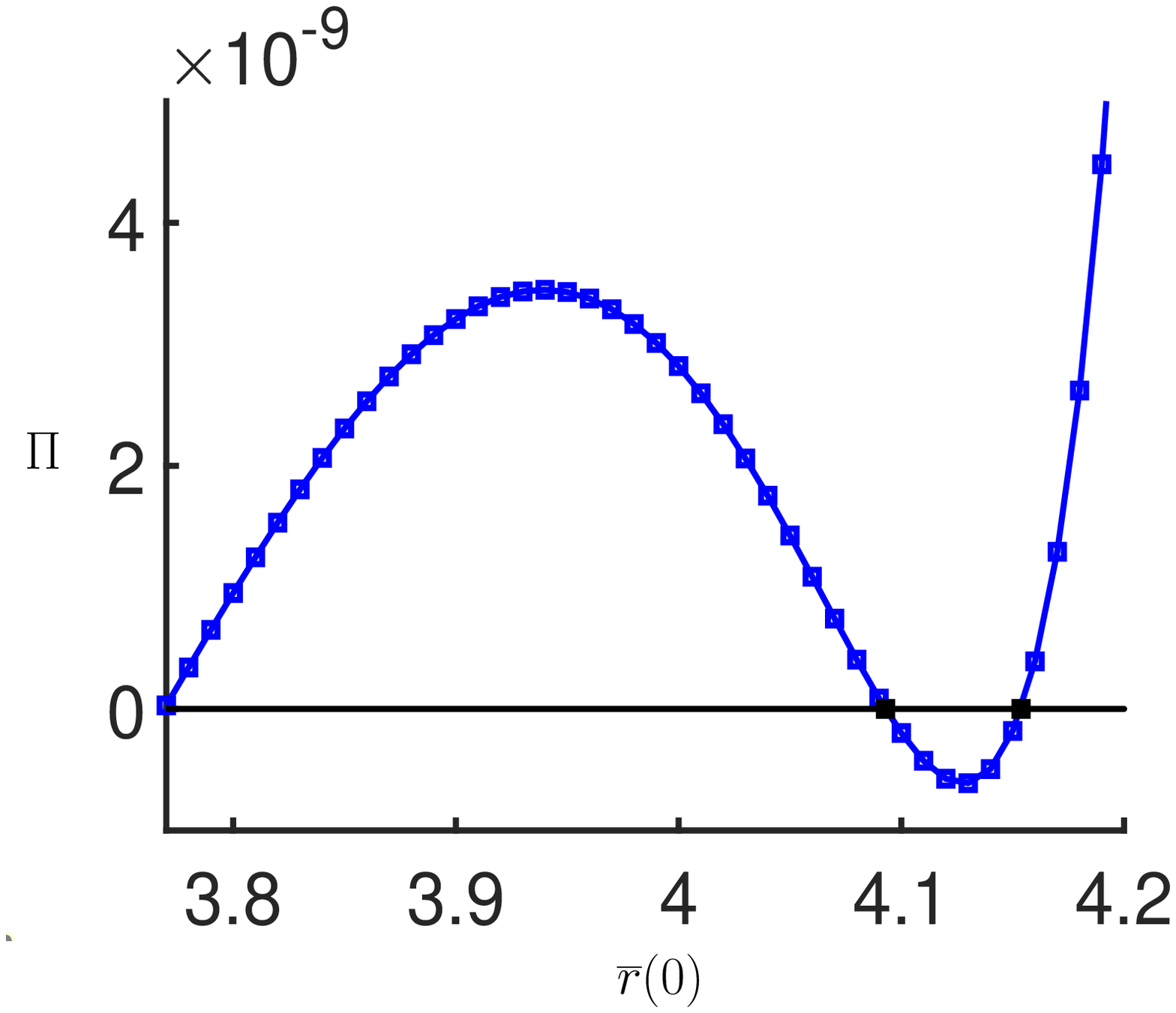}
  \caption{(color online) (a) The stability properties of the first 3 cells as a function of $\beta_r$. The black lines signify the position of the limit cycle (for $\rbar(0)>\rbar_{\eq}$) and the red lines give the position of the equilibrium points $\rbar_{\eq}$ (foci). Solid lines depict a stable focus/limit cycle and dashed lines signify an unstable focus/limit cycle. 
In (b) and (c) we plot the cell 1 and 2 stability respectively for the parameters in (\ref{eqn:parameters}) except with $\beta=-5.0-27.84\ri$. In (b) the PS is used to calculate the result, while in (c) we use the POS.}
 \label{fig:robust}
 \end{figure}
In figure \ref{fig:robust}(a) we plot the stability properties of cells 1-3 as a function of $\beta_r$. The figure shows a range of complex dynamics occurring as $\beta_r$ varies, including the generation of limit cycles through, what are likely to be, saddle-node bifurcations, and the merger of limit cycles with the equilibrium foci. From our initial assumption that the essential spectrum of the linear operator is in the left half plane, we are restricted to $\beta_r<0$. This is because there is no theory developed for the centre manifold reduction scheme {\ctext for $\beta_r\geq 0$} and hence it is not expected to work. However, we did extend this figure to $\beta_r=1$ and found that the scheme did appear to work. Hence we found that for $\beta_r\gtrsim0.95$ there is no limit cycle in cell 1 and an unstable limit cycle in both cell 2 and 3. In this case the dynamics are very similar to those documented for $\beta_r=-0.02$ above. Any solution trajectory starting in cell 1 converges to the stable focus in that cell, while in cells 2 and 3, trajectories starting inside the limit cycles converge to the foci in those cells, and those starting outside the limit cycles cascade down to cell 1 and converge at the focus in that cell. As $\beta_r$ is reduced, the dynamics remain similar, except there is now convergence to a limit cycle in cells 1-3 as $\beta_r$ passes approximately $0.95,~-0.49$ and $-1.42$ respectively.

In cell 2 the dynamics are relatively straight forward, and the two limit cycles move towards each other and eventually merge and annihilate at $\beta_r\approx-5.24$. In cell 3 the dynamics are different and here an unstable limit cycle persists for the whole range of $\beta_r$ studied, and there is a stable limit cycle which merges with an unstable limit cycle, which connects to the equilibrium point. Thus for $-4.88\lesssim\beta_r\lesssim-4.45$ there exists two stable bound states in cell 3 for the trajectories to converge to.

The dynamics in cell 1 are different again. We again find the stable limit cycle merges with an unstable limit cycle at $\beta_r\approx-6.28$, but this time the unstable limit cycle does not persist to positive $\beta_r$, as we found in cell 2, and instead merges with a second stable limit cycle at $\beta_r\approx-4.72$. Therefore, for $-6.28\lesssim\beta_r\lesssim-4.72$ there are two stable limit cycles for trajectories to converge to (and an unstable one {\ctext in between}). This also means that any solution trajectories starting outside the unstable limit cycles in cells 2 and 3 cascade down and converge to the outer of the two stable limit cycles in this parameter region. A plot showing $\Pi(\rbar(0))$ depicting the 3 limit cycles in cell 1 and the 2 limit cycles in cell 2 for $\beta_r=-5$ are given in figures \ref{fig:robust}(b) and \ref{fig:robust}(c) respectively.  

For $\beta_r<-6.28$ only a stable limit cycle persists in cell 1, and this is located at the very edge of the cell, while in cell 3 there is also a stable focal equilibrium point. Thus in this region all trajectories starting inside the unstable limit cycle in cell 3 converge to the equilibrium point, and all other trajectories starting outside this limit cycle, or anywhere in cells 1 and 2, converge to the limit cycle in cell 1. As $\beta_r$ is reduced further we ultimately expect the stable and unstable limit cycles in cells 1 and 3 respectively to move further towards the outside of their respective cells and vanish, possibly via a homoclinic bifurcation with the saddle points at the cell boundary. However, we do not see this bifurcation for the values of $\beta_r$ we were able to calculate in sensible time-frames.

The question, from where, and by which mechanism these periodic orbits bifurcate is left as an open problem. We attempted to resolve this problem by varying the other system parameters, namely $\alpha,~\gamma,~\delta$, to the point where the periodic orbits in cell 2 move to the focus, for all values of $\beta_r$. We considered cell 2 so that we could utilize the POS for computational speed. However, when varying these other parameters we identified similar figures to that in figure \ref{fig:robust}(a), with the stable periodic orbit region translated to larger or smaller $\beta_r$, depending on whether the varied parameter is increased or decreased. Ultimately, as these parameters are varied further, we reached one of two problems: either the position of the cell moves to larger radii values, hence the POS becomes inaccurate due to the values of $\Pi(r(0))$ becoming of the order of machine precision; or the position of the cell moves to smaller radii, meaning the POS is no longer valid and the PS is required to calculate the dynamics. But this vastly increases the computational running time per parameter making it unwieldy.

\subsection{Three-pulse interaction dynamics, $\beta_r=-0.02$}

Finally, to demonstrate the full power of the PS, we investigate the interaction of three pulses in the QCGLE for two sets of initial conditions. Firstly we consider a simple case of three equidistant, in-phase, pulses, hence we expect to have a solution with symmetry about the middle pulse, and secondly we consider a case where this symmetry is broken by examining equidistant pulses with differing phases. Here the six time dependent positions and phases of the pulses are $(r_1,r_2,r_3)$ and $(g_1,g_2,g_3)$, with $r_1<r_2<r_3$.

\subsubsection{Case 1: Symmetric pulse distribution}

\begin{figure}[ht!]
(a)\includegraphics[width=0.45\textwidth]{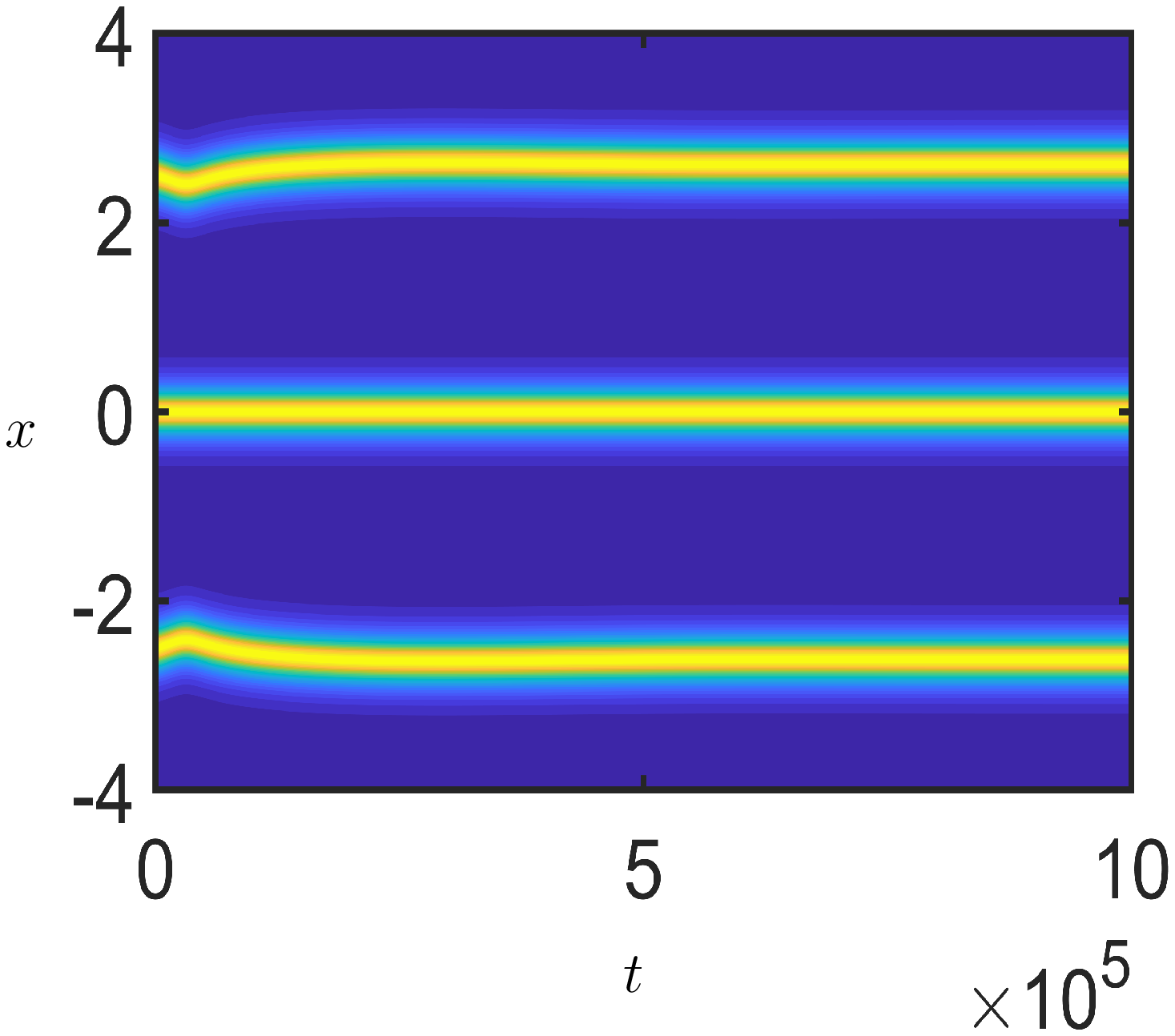}
(b)\includegraphics[width=0.45\textwidth]{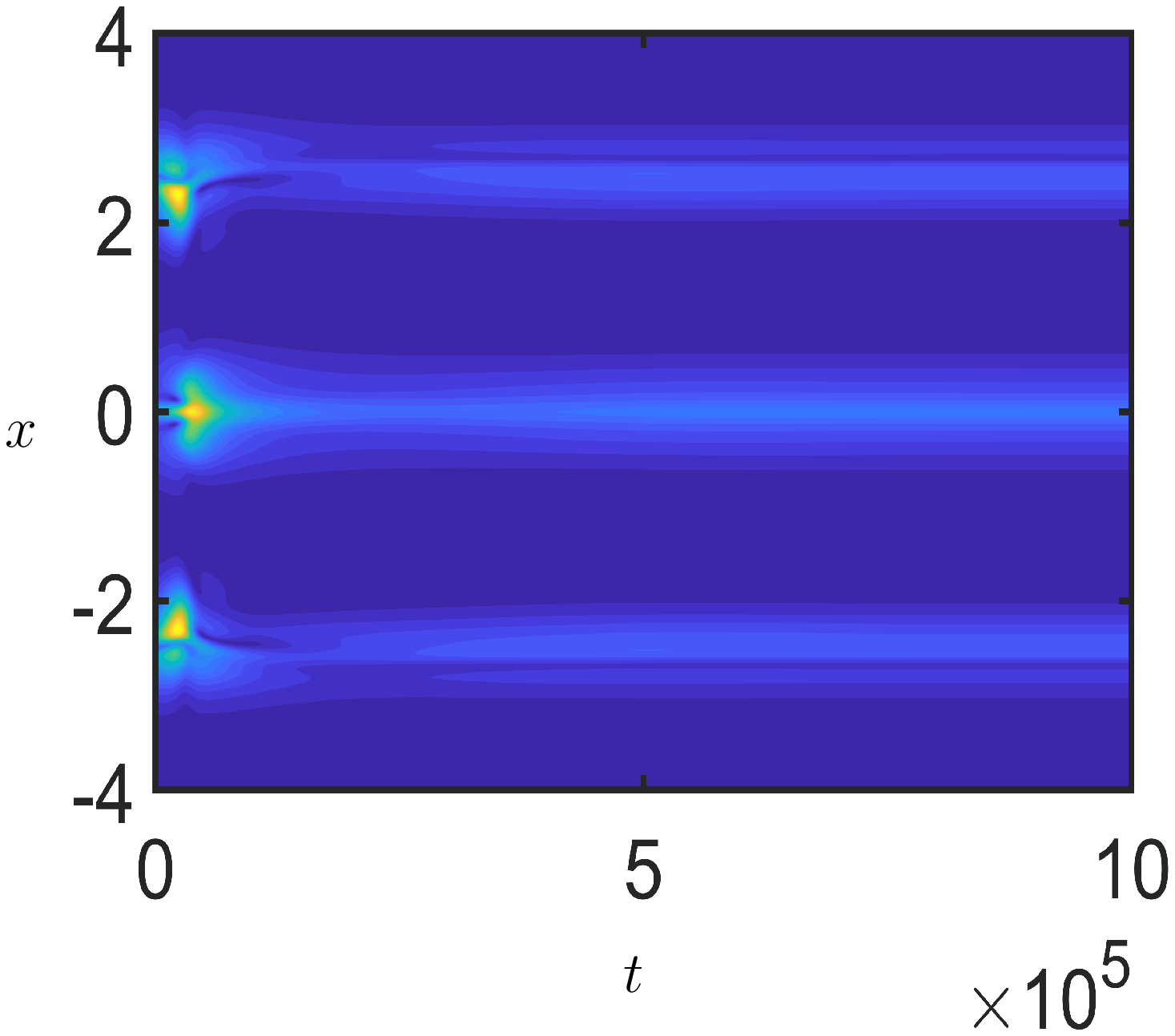}
(c)\includegraphics[width=0.45\textwidth]{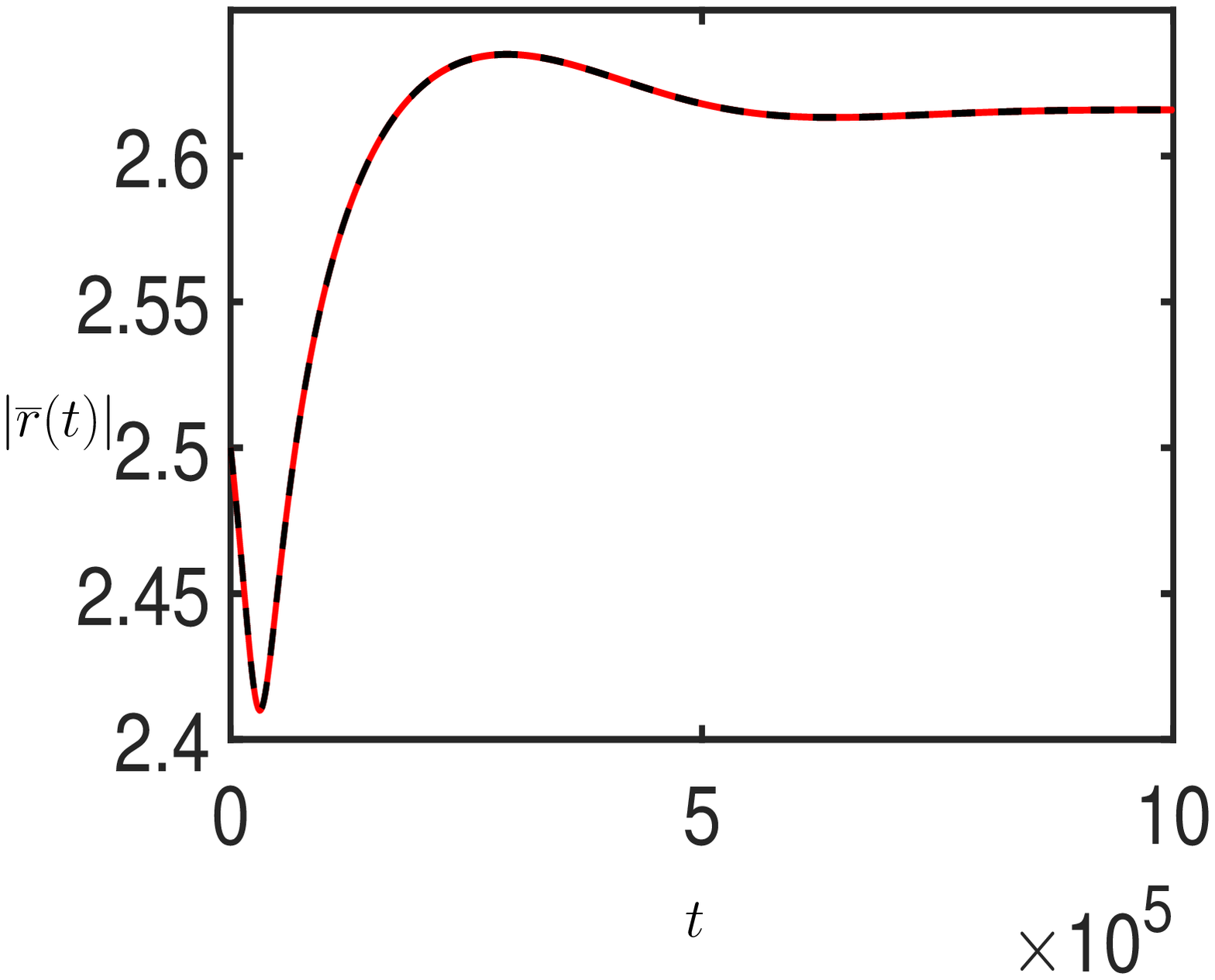}
(d)\includegraphics[width=0.45\textwidth]{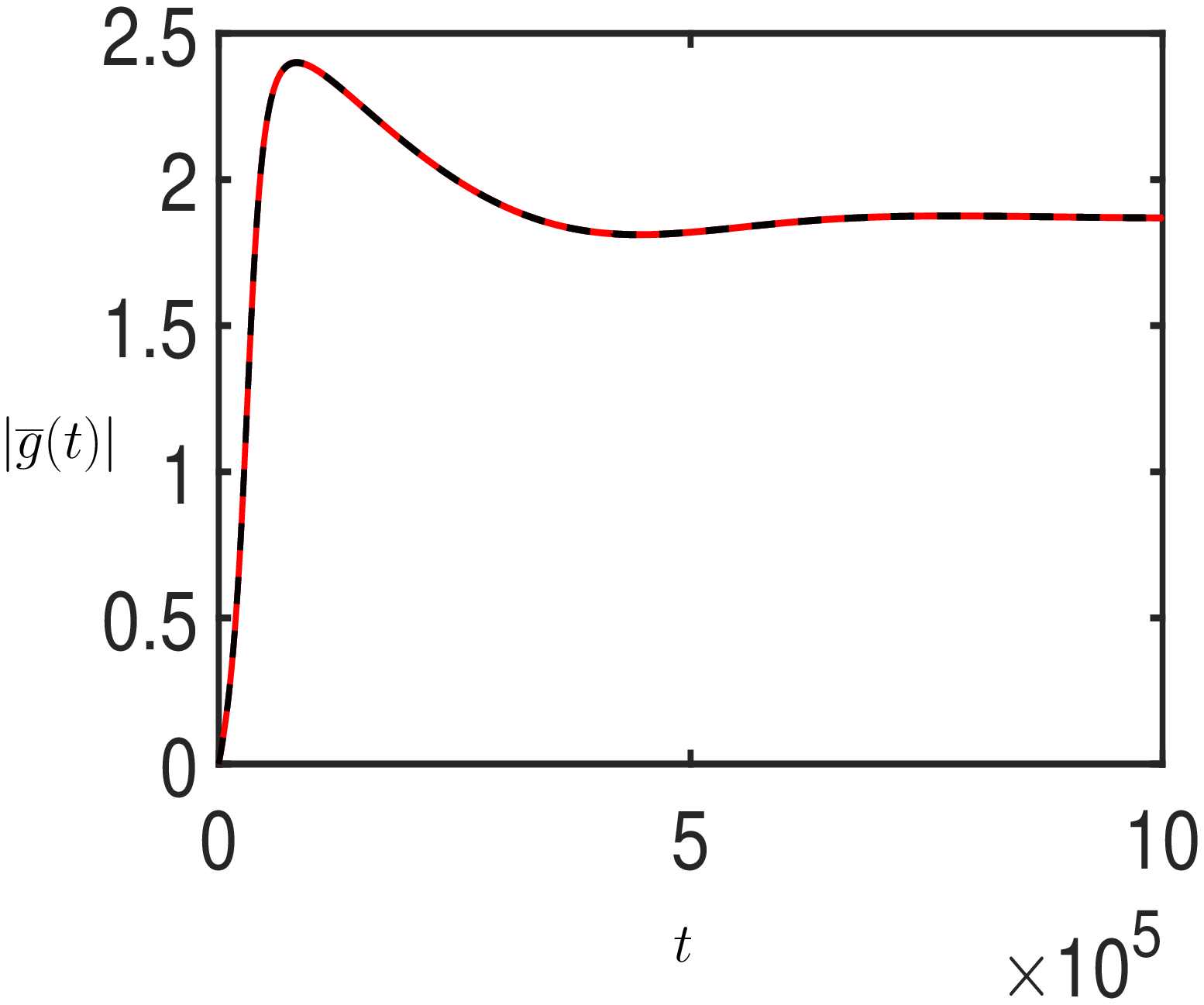}
\caption{(colour online) Three equidistant in-phase pulse interaction with initial conditions $r_1 = -2.5,\; r_2 = 0.0,\; r_3 = 2.5,\quad g_1 = g_2 = g_3 = 0$. (a) Profile of pulse evolution $|u(t,x)| = |V_1(x)+V_2(x)+V_3(x)+w(t,x)|$, (b) remainder function $|w(t,x)|$, (c) $|\rbar(t)|$  and (d) $|\gbar(t)|$. In panels (c) and (d) the black dashed line gives $|\rbar|=|r_1-r_2|$ and $|\gbar|=|g_1-g_2|$, while the red solid line gives $|\rbar|=|r_3-r_2|$ and $|\gbar|=|g_3-g_2|$.}
\label{fig:col_w_case1}
\end{figure}

In this case we consider the initial positions of the pulses to be at
\begin{equation}
r_1 = -2.5,~~~ r_2 = 0.0,~~~r_3 = 2.5~~~{\rm with}~~~g_1 = g_2 = g_3 = 0.
\label{eqn:3case1_parameters}
\end{equation}
The motion of the three pulses is given in figure \ref{fig:col_w_case1}(a) which plots $|u(t,x)|=|V_1+V_2+V_3+w(t,x)|$. Here we see that both outer pulses are initially attracted to the centre pulse, however after some initial transient region, the interaction between the three pulses appears to cancel each other out, resulting in a three-pulse bound state. In figure \ref{fig:col_w_case1}(b), we consider the remainder function $|w(t,x)|$, which is small initially, but it increases in magnitude to reach its maximum at the point where the two outer pulses are closest to the central pulse. Once the pulses converge to the bound state, the form of $w(t,x)$ becomes a constant function in time.

In figures \ref{fig:col_w_case1}(c) and (d), we consider the differences $|r_3-r_2|,|r_2-r_1|$ and $|g_3-g_2|,|g_2-g_1|$ respectively, which show the movement and phase change in the outer pulses compared to the central pulse. As seen in figure \ref{fig:col_w_case1}(a), the outer pulses translate and change
their phase identically with relation to the central pulse, before converging to their bound state. In particular, in figure \ref{fig:col_w_case1}(c) we observe that the final separation distance of the pulses in the bound state is increased in comparison to the initial separation. Similarly, the phase difference, converges to $|\gbar|\approx 1.8$ for each pair.

We observe identical results to this when we change the initial phase variables to $\gbar = \pi,$ or $\frac{\pi}{2}$. In particular, we converge to the same bound state as the one shown in figure \ref{fig:col_w_case1}.
 The same bound state is also reached for three in-phase pulses with an initial separation distance $\rbar\lesssim1.94$, except in this case the outer two pulses initially move away from the central pulse before converging to the bound state.

\subsubsection{Case 2: Asymmetric pulse distribution}

\begin{figure}[ht!]
\begin{center}
(a)\includegraphics[width=0.45\textwidth]{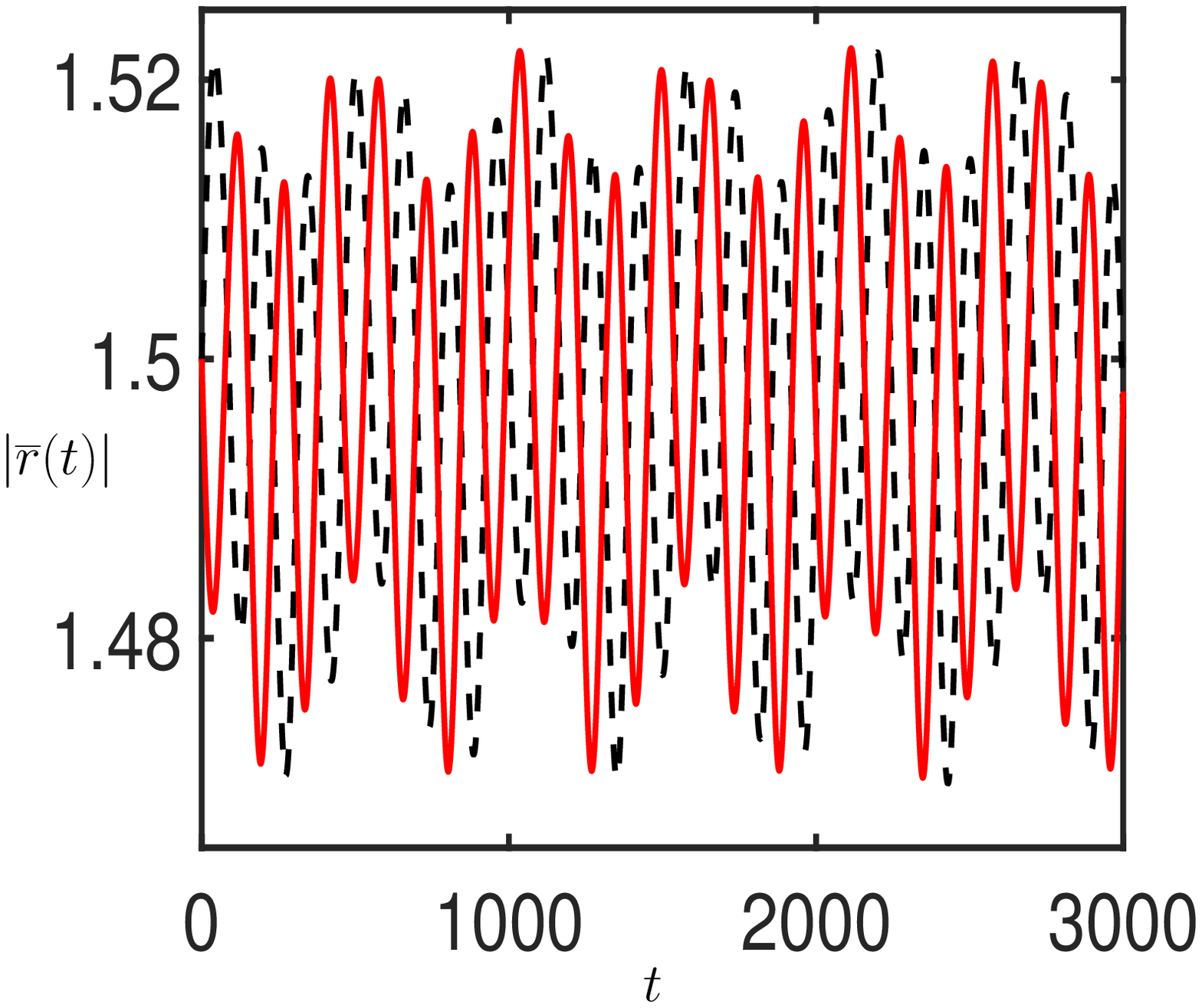}
(b)\includegraphics[width=0.45\textwidth]{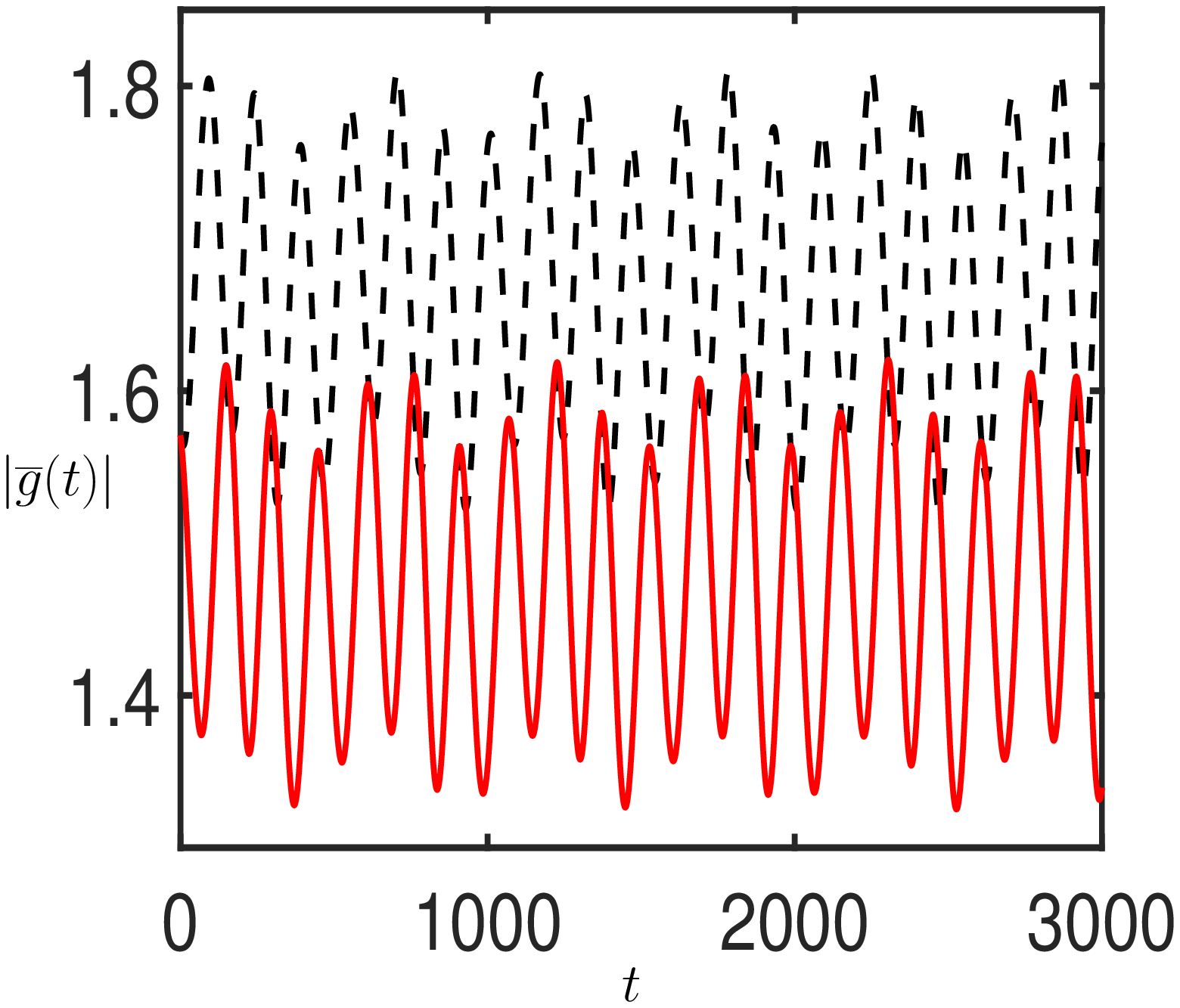}
(c)\includegraphics[width=0.45\textwidth]{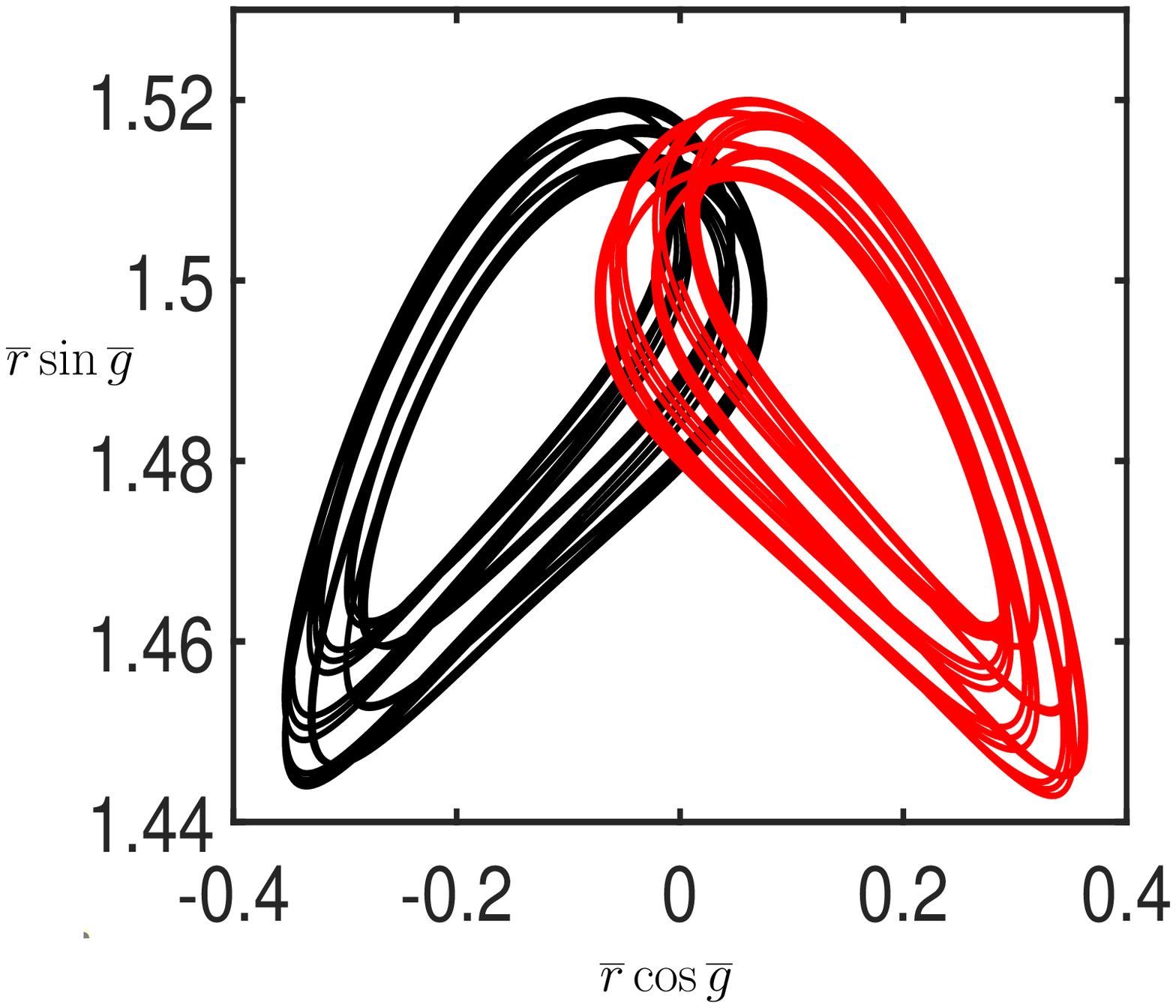}
\end{center}
\caption{(colour online) Three equidistant, asymmetric phase, pulse interaction with initial conditions $r_1 = -1.5,~r_2 = 0.0,~r_3 = 1.5,\quad g_1 =\pi,~g_2 =\pi/2,~g_3 = 0$. (a) $|\rbar(t)|$, (b) $|\gbar(t)|$ and (c) the phase plane where the black line gives $|\rbar|=|r_1-r_2|$ and $|\gbar|=|g_1-g_2|$, while the red line gives $|\rbar|=|r_3-r_2|$ and $|\gbar|=|g_3-g_2|$.}
\label{fig:case22_1c}
\end{figure}

In this case we consider three pulses with the initial parameters
\begin{equation}
r_1 = -1.5,~~~ r_2 = 0.0,~~~r_3 = 1.5~~~{\rm with}~~~g_1 = \pi,~~~~g_2 =\frac{\pi}{2},~~~~ g_3 = 0,
\label{eqn:3case2_parameters}
\end{equation}
and follow the time interaction of the pulses in figure \ref{fig:case22_1c}. In this case the distance and phase difference between the two pulse pairs varies in time, but without a clearly defined period. When looking at the phase-plane in figure \ref{fig:case22_1c}(c) we see that in this case the dynamics appear to be chaotic. 
 We have not performed a calculation here to determine whether or not the system is in fact chaotic, but this result shows that the PS presented in this paper allows us to readily study such systems of multiple pulses in reasonable computation times. 

{\ctext We present results in this section using the PS only, but equally we could have used the POS too. However, in order to use the POS we would need to calculate a large number of coefficients which appear in the inner products of \eqref{3.ODEasy}. As there is currently no systematic way to do this, the PS becomes the more effective approach for higher numbers of pulses.}

\section{Conclusions and discussion}\label{sec:Conclusion}

In this paper, we investigated the dynamics of multi-pulse solutions in the QCGLE. We highlighted a Projected Scheme (PS), which by writing the solution as a sum of the pulses plus a remainder term, $u(t,x)=V_{\vec\xi}(x)+w(t,x)$, allows for the dynamics to be investigated via a slow ODE system, coupled to a fast PDE for the remainder function. The PS has a computational speed advantage compared to more traditional numerical schemes, which allowed for a more comprehensive investigation of multi-pulse interactions in the QCGLE than ever before. We also introduced a method to approximate the contribution of the remainder function $w(t,x)$ in the very-well-separated pulse limit, the Projected ODE Scheme (POS), which gave an even greater computational speed up, and allowed us to investigate the dynamics in the very-weakly-interacting region of the phase plane. The speed up of these two systems is of the order of 150 times and 20000 times faster respectively than a traditional Standard time-stepping Scheme (SS).

Having formulated the centre-manifold reduction projection scheme for the case of multiple pulses, we focused initially on the two-pulse system, and identified the system dynamics from weakly-interacting to very-weakly-interacting pulses. We performed a parameter study for this system as a function of $\beta_r=\Re(\beta)$. For a fairly typical parameter set ({\ctext $\beta_r=-0.02$}) we found that each cell of the phase plane $(\rbar\cos\gbar,\rbar\sin\gbar)$, contained a periodic limit cycle. In cell 1 this limit cycle was found to be stable, while in cells 2 and 3 it was found to be unstable. Thus for all initial conditions such that the phase plane trajectory starts in cell 1, the dynamics led to the trajectory converging to the limit cycle. For initial conditions where the trajectory starts inside the cell 2 or 3 limit cycle, then the solution converged to the equilibrium fixed point within the limit cycle, with constant separation and a phase difference of $\gbar=\pi/2$. However, for conditions such that the trajectory begins outside, either of these limit cycles, but inside cell 2 or 3, then the deformation of the heteroclinic orbits at the edges of these cells meant that the trajectory spiraled out to the edge of the cell and then cascaded down and into cell 1, where the solution converged to the cell 1 limit cycle.

By varying the parameter $\beta_r$ we were able to investigate the robustness of the solution structure found for $\beta_r=-0.02$ for cells 1-3. What was discovered were apparent saddle-node bifurcations leading to the generation and annihilation of the limit cycles observed. For large positive $\beta_r$, the limit cycles in cells 2 and 3 persist and are unstable,  while there is no limit cycle in cell 1. In this case any trajectories beginning inside these limit cycles converged to the foci in those cells, while trajectories beginning outside these limit cycles, or beginning in cell 1, converged to the focus in cell 1. As $\beta_r$ is reduced, stable limit cycles bifurcate from the foci and move towards the edge of the cell and eventually annihilate with an unstable limit cycle. The most interesting dynamics occurred in cell 1, where there was a fold in the bifurcation parameter-space leading to the existence of three limit cycles for the range $-6.28\lesssim\beta_r\lesssim-4.72$. Two of these limit cycles were stable, with an unstable limit cycle located inbetween, leading to an extra bound state to which the system could converge. 

Finally, we investigated the interaction of three pulses for two separate initial conditions. The first showed the three pulses converging to a steady bound state with constant separation and phase differences, while the second showed the pulses experiencing chaotic behaviour. Both these behaviours were found for the same set of QCGLE parameters, hence showing the rich variety of dynamics in the three-pulse solution of the QCGLE.

From the results presented, it is very clear that the POS has significant speed up advantages over the PS, but because we fix the inner product constants in the very-weakly-interacting region of parameter space, its quantitative results become less accurate as the two-pulses come close together during a cycle of their dynamics, as in this region the very-weakly-interaction approximation is not valid. However, it is still good to know that the POS gives excellent qualitative information about the system. Hence if a solution for $w(t,x)$, from \eqref{leadingw}, in terms of $\rbar$ and $\gbar$ could be written down explicitly then this would allow the POS to extend to the same range of validity as the PS. However, as was shown in the inner product calculation in Appendix \ref{appen:innerprods}, a large number of higher order exponential terms are required in the solution for $w$ in order to accurately calculate the relevant inner products to $O(e^{-\lambda_r|\rbar|})$ and $O(e^{-2\lambda_r|\rbar|})$.

One of the main advantages of the PS and POS, over more traditional numerical schemes, is that the explicit computation of the small interaction terms means we are able to keep close control on all the numerical error tolerances. However, we are still limited by machine precision arithmetic, i.e. for double precision considered here, we are limited to separation distances of less than $\log(\nu\times 10^{-16})$ where $\nu$ is the linear decay rate to the base state. In this paper this meant that we could use the POS scheme to examine the dynamics in cell 3, but not beyond. 
To examine beyond cell 3 we could use a multiple precision package to perform the calculation of the initial pulse solution and the eigenfunctions. This however, is a non-trivial exercise, as any built-in chebfun subroutines from \textsc{matlab} used, need to be re-written by the user for the arbitrary precision software.

The global centre-manifold theory used here and in Zelik \& Mielke~\cite{ZelMiel09}, has been proven for a general parabolic PDE system with a strongly elliptic differential operator, which involves general weakly-interacting localised structures in several space dimensions. Hence, the numerical Projected Scheme presented in this paper has the potential to extend to a large class of problems. In particular this scheme should allow one to study the interaction of localised pulses in the two-dimensional plane, which form either spots or hexagon patches. 

\vspace{0.5cm}

{\Large{\bf Acknowledgments}}

\vspace{0.5cm}

MRT and DJBL would like to thank Jamie Chavez Malacara whose work on an EPSRC funded Summer Internship helped in the completion of this work. TR would like to thank the EPSRC for funding his PhD studentship. For the purpose of open access, the authors have applied a Creative Commons Attribution (CC BY) license to any Author Accepted Manuscript version arising.

\vspace{0.5cm}

{\Large{\bf Data Access Statement}}

\vspace{0.5cm}

The data underlying this article will be shared, on request to the corresponding author.

\appendix

\renewcommand{\theequation}{\thesection.\arabic{equation}}

\section{The second order system formulated as ODEs for $n=2$ (POS formulation)}
\label{appen:innerprods}

For the leading order form of the slow system \eqref{3.ODE4}, \cite{ZelMiel09} evaluated the inner products explicitly to give closed form solutions. In this appendix we show how to determine, and numerically evaluate the form of the inner products in the Projected System, in the very-weakly-interacting regime, such that the coupled ODE/PDE system can be reduced to just an ODE system, which leads to a rapid decrease in computation time. From \eqref{leadingw}, the quasi-steady equation for $\what$ is a function of $\vec\xi$ only, i.e. it has no explicit $\that$ dependence. Hence inner products involving this function should also have an explicit form in terms of $r_1(t),r_2(t),g_1(t)$ and $g_2(t)$ only, with no explicit $t$ dependence. This appendix focuses on the two-pulse case, but this approach could be used to identify a POS for $n>2$ pulses too.

\begin{figure}
\begin{center}
\includegraphics[width=0.55\textwidth]{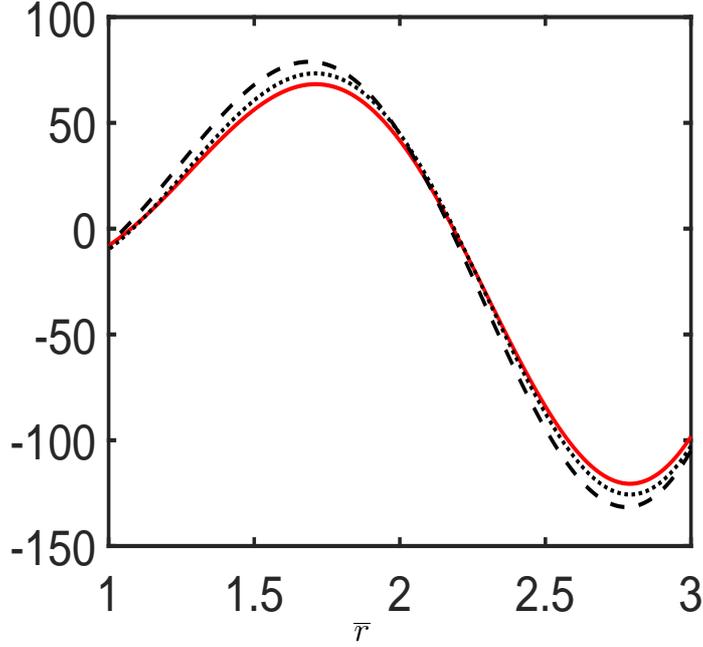}
\end{center}
\caption{Numerical form (solid line) and well-separated asymptotic forms given by (\ref{eqn:innerprod1}) (dashed line) and (\ref{eqn:innerprod2}) (dotted line) of the inner product $e^{\lambda_r|\rbar|}\langle \varphi^r_{\xi_2}, \psi^r_{\xi_1} \rangle$ with $\gbar=\sqrt{2}$.}
\label{fig:appen_inner}
\end{figure}

Before considering the inner products containing $\what(t,x)$ we first consider the inner products which occur in the matrix $\widehat{\mathcal{C}}$, which have the form $\langle \varphi^r_{\xi_2},\psi^r_{\xi_1}\rangle$. These inner products can be evaluated by using the far-field asymptotic forms of the eigenfunctions \eqref{2.tails} as 
\begin{eqnarray}
\langle \varphi^r_{\xi_2},\psi^r_{\xi_1}\rangle&=&\Re\left[\int_{-\infty}^\infty\left[-e^{\ri g_2} V'(x-r_2)\right]\overline{\left[e^{\ri g_1}\psi^r(x-r_1)\right]}\,\d x\right],\nonumber\\
&=&\Re\left[\int_{r_1}^{r_2}-\lambda p\overline{q}e^{\ri\gbar}e^{-\lambda\rbar}\,\d x\right],\nonumber\\
&=&A_1\rbar e^{-\lambda_r\rbar}\sin(-\lambda_i\rbar+\gbar+\kappa_{1}),
\label{eqn:innerprod1}
\end{eqnarray}
where
\[
A_1e^{\ri\kappa_{1}}=-\ri p \overline{q}\lambda.
\]
This approximation to the inner product $\langle \varphi^r_{\xi_2},\psi^r_{\xi_1}\rangle$ is given by the dashed line in figure \ref{fig:appen_inner}, and this appears to show that an additional $O(e^{-\lambda_r\rbar})$ term is missing.
 The resolution of this missing term comes from the consideration of higher order terms in the asymptotic expansions \eqref{2.tails}, where in fact we can show that
\begin{eqnarray}
V(x)&=&pe^{-\lambda|x|}+\phat e^{-2\lambda_r|x|}e^{-\lambda|x|}+O(e^{-5\lambda_r|x|}),\nonumber\\
\psi^r(x)&=&q\sgn(x)e^{-\lambda|x|}+\sgn(x)e^{-3\lambda_r|x|}\left(\qhat_1e^{\ri\lambda_i|x|}+\qhat_2e^{-3\ri\lambda_i|x|}\right)+O(e^{-5\lambda_r|x|}),\label{2.tails2}\\
\psi^g(x)&=&se^{-\lambda|x|}+e^{-3\lambda_r|x|}\left(\shat_1e^{\ri\lambda_i|x|}+\shat_2e^{-3\ri\lambda_i|x|}\right)+O(e^{-5\lambda_r|x|}),\nonumber
\end{eqnarray}
with
\begin{eqnarray*}
\phat&=&-\frac{\gamma|p|^2p}{\alpha(9\lambda_r^2+6\ri\lambda_r\lambda_i-\lambda_i^2)+\beta},\\
\qhat_1&=&-\frac{2|p|^2\overline{\gamma}q}{\overline{\alpha}(9\lambda_r^2-6\ri\lambda_r\lambda_i-\lambda_i^2)+\overline{\beta}},~~~~\qhat_2=-\frac{p^2\gamma\overline{q}}{\overline{\alpha}(9\lambda_r^2+18\ri\lambda_r\lambda_i-9\lambda_i^2)+\overline{\beta}},\\
\shat_1&=&-\frac{2|p|^2\overline{\gamma}s}{\overline{\alpha}(9\lambda_r^2-6\ri\lambda_r\lambda_i-\lambda_i^2)+\overline{\beta}},~~~~\shat_2=-\frac{p^2\gamma\overline{s}}{\overline{\alpha}(9\lambda_r^2+18\ri\lambda_r\lambda_i-9\lambda_i^2)+\overline{\beta}}.
\end{eqnarray*}
This time when evaluating the inner product $\langle \varphi^r_{\xi_2},\psi^r_{\xi_1}\rangle$ the solution becomes
\begin{equation}\label{eqn:innerprod2}
\langle \varphi^r_{\xi_2},\psi^r_{\xi_1}\rangle=A_1\rbar e^{-\lambda_r\rbar}\sin(-\lambda_i\rbar+\gbar+\kappa_{1})+A_2e^{-\lambda_r\rbar}\sin(-\lambda_i\rbar+\gbar+\kappa_{2}),
\end{equation}
with 
\[
A_2e^{\ri\kappa_{2}}=-\ri\left[\frac{\qbar\widehat{p}(\lambda+2\lambda_r)}{2\lambda_r}+\frac{\lambda p\overline{\widehat{q}_1}}{2\lambda_r}+\frac{\lambda p \overline{\widehat{q}_2}}{2(\lambda_r-2\ri\lambda_i)}\right],
\]
and this approximation is plotted as the dotted line in figure \ref{fig:appen_inner}. This approximation is clearly better than the approximation in \eqref{eqn:innerprod1}, but as this calculation shows, higher order asymptotic terms in the expansions \eqref{2.tails2} are required to give the correct, and accurate, form of the inner product at $O(e^{-\lambda_r\rbar})$. Hence, in order to get accurate forms of the coefficients in \eqref{eqn:innerprod2} we would require the inclusion of the $O(e^{-5\lambda_r|x|})$, the $O(e^{-7\lambda_r|x|})$ terms etc, in the expansions \eqref{2.tails2}. While this may be tangible for the inner products of the form \eqref{eqn:innerprod2}, this is not practical for those inner products which involve $\what$.

We can overcome this technicality and calculate asymptotic forms for all the inner products in \eqref{asyODE.compact} by numerically calculating the very-weakly-interacting form of the inner products, $\rbar\gtrsim1.5$, for various values of $\rbar$ and $\gbar$ and then approximating a best fit surface to these values using \textsc{matlab}'s \verb2fit2 subroutine. This routine gives a 95\% confidence interval for the coefficients of a pre-specified function. To make sure we are in this very-weakly-interacting regime, we choose values in the range
\begin{equation}
\rbar\in[-\log(10^{-12})/(2\lambda_r),-\log(10^{-15})/(2\lambda_r)],~~~{\rm and}~~~\gbar\in[-\pi,\pi]. 
\label{eqn:fittingregion}
\end{equation}
Examining the numerical form of the inner products we are able to determine that they have the following very-weakly-interacting forms
\begin{eqnarray*}
\langle \Phihat,\psi^{r,g}_{\xi_1}\rangle&=&A_1\sin(-\lambda_i\rbar+\gbar+\kappa_1),\\
\langle \Phihat,\psi^{r,g}_{\xi_2}\rangle&=&A_1\sin(-\lambda_i\rbar-\gbar+\kappa_1),\\
\langle \varphi^{r,g}_{\xi_2},\psi^r_{\xi_1}\rangle&=&e^{-\lambda_r\rbar}\left(-A_1\rbar\sin(-\lambda_i\rbar+\gbar+\kappa_1)-A_2\sin(-\lambda_i\rbar+\gbar+\kappa_2)\right),\\
\langle \varphi^{r,g}_{\xi_1},\psi^r_{\xi_2}\rangle&=&e^{-\lambda_r\rbar}\left(-A_1\rbar\sin(-\lambda_i\rbar-\gbar+\kappa_1)-A_2\sin(-\lambda_i\rbar-\gbar+\kappa_2)\right),\\
\langle \what,\partial_x\psi^{r,g}_{\xi_1}\rangle=\langle \ri\what,\partial_x\psi^{r,g}_{\xi_1}\rangle&=&A_1\cos(-\gbar+\lambda_i\rbar+\kappa_1),\\
\langle \what,\partial_x\psi^{r,g}_{\xi_2}\rangle=\langle \ri\what,\partial_x\psi^{r,g}_{\xi_2}\rangle&=&A_1\cos(\gbar+\lambda_i\rbar+\kappa_1),\\
\langle \Psihat_1\what,\partial_x\psi^{r,g}_{\xi_1}\rangle=\langle \Ghat,\partial_x\psi^{r,g}_{\xi_1}\rangle&=&
B_0+B_1\cos(-2\gbar+2\lambda_i\rbar+\mu_1)+B_2\cos(-2\lambda_i\rbar+\mu_2)+B_3\cos(-2\gbar+\mu_3),\\
\langle \Psihat_2\what,\partial_x\psi^{r,g}_{\xi_2}\rangle=\langle \Ghat,\partial_x\psi^{r,g}_{\xi_2}\rangle&=&
B_0+B_1\cos(2\gbar+2\lambda_i\rbar+\mu_1)+B_2\cos(-2\lambda_i\rbar+\mu_2)+B_3\cos(2\gbar+\mu_3),\\
\end{eqnarray*}
where $A_i,\kappa_i,B_i$ and $\mu_i$ have different values in each inner product. Therefore the ODEs integrated forward in time are then given by
\begin{eqnarray}
\frac{\d\overline{r}}{\d\that}&=&\langle \Phihat,\psi^r_{\xi_2}\rangle\left[1-e^{-\lambda_r\rbar}\langle \what,\partial_x\psi^r_{\xi_2}\rangle+\langle \varphi^r_{\xi_2},\psi^r_{\xi_1}\rangle\right]-\langle \Phihat,\psi^r_{\xi_1}\rangle\left[1-e^{-\lambda_r\rbar}\langle \what,\partial_x\psi^r_{\xi_1}\rangle+\langle \varphi^r_{\xi_1},\psi^r_{\xi_2}\rangle\right]\nonumber\\
&&+\langle \Phihat,\psi^g_{\xi_1}\rangle\left[e^{-\lambda_r\rbar}\langle \ri \what,\psi^r_{\xi_1}\rangle-\langle \varphi^g_{\xi_1},\psi^r_{\xi_2}\rangle\right]+\langle \Phihat,\psi^g_{\xi_2}\rangle\left[-e^{-\lambda_r\rbar}\langle \ri \what,\psi^r_{\xi_2}\rangle+\langle \varphi^g_{\xi_2},\psi^r_{\xi_1}\rangle\right]\nonumber\\
&&+e^{-\lambda_r\rbar}\langle \Psihat_2,\psi^r_{\xi_2}\rangle-e^{-\lambda_r\rbar}\langle \Psihat_1,\psi^r_{\xi_1}\rangle+e^{-\lambda_r\rbar}\langle \Ghat,\psi^r_{\xi_2}\rangle-e^{-\lambda_r\rbar}\langle \Ghat,\psi^r_{\xi_1}\rangle,\label{eqn:ode1}\\
\frac{\d\overline{g}}{\d\that}&=&\langle \Phihat,\psi^g_{\xi_2}\rangle\left[1-e^{-\lambda_r\rbar}\langle \ri \what,\psi^g_{\xi_2}\rangle+\langle \varphi^g_{\xi_2},\psi^g_{\xi_1}\rangle\right]-\langle \Phihat,\psi^g_{\xi_1}\rangle\left[1-e^{-\lambda_r\rbar}\langle \ri \what,\psi^g_{\xi_1}\rangle+\langle \varphi^g_{\xi_1},\psi^g_{\xi_2}\rangle\right]\nonumber\\
&&+\langle \Phihat,\psi^r_{\xi_1}\rangle\left[e^{-\lambda_r\rbar}\langle \what,\partial_x\psi^g_{\xi_1}\rangle-\langle \varphi^r_{\xi_1},\psi^g_{\xi_2}\rangle\right]+\langle \Phihat,\psi^r_{\xi_2}\rangle\left[-e^{-\lambda_r\rbar}\langle \what,\partial_x\psi^g_{\xi_2}\rangle+\langle \varphi^r_{\xi_2},\psi^g_{\xi_1}\rangle\right]\nonumber\\
&&+e^{-\lambda_r\rbar}\langle \Psihat_2,\psi^g_{\xi_2}\rangle-e^{-\lambda_r\rbar}\langle \Psihat_1,\psi^g_{\xi_1}\rangle+e^{-\lambda_r\rbar}\langle \Ghat,\psi^g_{\xi_2}\rangle-e^{-\lambda_r\rbar}\langle \Ghat,\psi^g_{\xi_1}\rangle.\label{eqn:ode2}
\end{eqnarray}

We can further constrain the parameters in the inner products by noting that we still require the fixed saddle point equilibrium points $(\rbar_{\rm eq},\gbar_{\rm eq})$ to lie at points with $\gbar=n\pi$ where $n$ is an integer, hence this leads restrictions on the coefficients such that
\[
\langle \what,\partial_x\psi^{r}_{\xi_1}\rangle=\langle \what,\partial_x\psi^{r}_{\xi_2}\rangle~~~{\rm and}~~~\langle \what,\partial_x\psi^{g}_{\xi_1}\rangle=-\langle \what,\partial_x\psi^{g}_{\xi_2}\rangle,
\]
at $\gbar=n\pi$ etc. Hence in total there are $48$ independent coefficients needed for the POS. For the parameters in \eqref{eqn:parameters} the 95\% confidence intervals for these coefficients to 2 decimal places are shown in table \ref{table:wcoeffs3}.
\begin{table}[!htb]
\begin{center}
\resizebox{\columnwidth}{!}{%
\begin{tabular}{|c|c|c|c|c|c|c|c|}

\hline

 & $A_1$ & $\kappa_1$ & $A_2$ & $\kappa_2$ & & & \\

\hline
$\langle \Phihat, \psi^r_{\xi_1} \rangle$  & (67.48,67.48)   & (2.33,2.33) & & & & &\\
$\langle \Phihat, \psi^g_{\xi_1} \rangle$  & (423.29,423.30) & (2.73,2.73) & & & & &\\
$\langle \varphi^r_{\xi_2}, \psi^r_{\xi_1} \rangle$  & (47.72,47.72)  & (-1.60,-1.60) & (11.79,11.79) & (1.33,1.33) & & &\\
$\langle \varphi^g_{\xi_2}, \psi^r_{\xi_1} \rangle$  & (6.52,6.52) & (2.72,2.72) & (0.75,0.75) & (-0.43,-0.43) & & &\\
$\langle \varphi^r_{\xi_2}, \psi^g_{\xi_1} \rangle$  & (299.32,299.32)  &(-1.20,-1.20) & (-58.46,-58.46) & (-1.36,-1.36) & & &\\
$\langle \varphi^g_{\xi_2}, \psi^g_{\xi_1} \rangle$  & (40.88,40.88)  & (3.12,3.12) & (2.85,2.85) & (0.28,0.28) & & &\\
$\langle \what, \partial_x\psi^r_{\xi_1} \rangle$  & (5.94,5.94)  & (-1.71,-1.71) &  & & & &\\
$\langle \what, \partial_x\psi^g_{\xi_1} \rangle$  & (10.99,10.99)  & (-0.19,-0.19) & & & & &\\
$\langle \ri\what, \psi^r_{\xi_1} \rangle$  & (1.91,1.91)  & (2.91,2.91) & & & & &\\
$\langle \ri\what, \psi^g_{\xi_1} \rangle$  & (5.99,6.00)  & (-3.00,-3.00) & & & & &\\

\hline
\hline

 & $B_0$ & $B_1$ & $\mu_1$ & $B_2$ & $\mu_2$ & $B_3$ & $\mu_3$ \\

\hline

$\langle \Psihat_1\what, \psi^r_{\xi_1} \rangle$  & (-860.89,-860.50)  & (1111.08,1111.39) & (2.26,2.26) & (699.22,699.75)& (-2.48,-2.48)& (575.40,575.72)& (2.14,2.14)\\
$\langle \Psihat_1\what, \psi^g_{\xi_1} \rangle$  & (-679.49,-677.29)  & (5212.81,5214.75) & (1.71,1.71) & (4386.23,4388.67) & (-2.08,-2.08) & (3609.25,3611.04) & (1.74,1.74)\\
$\langle \Ghat_1\what, \psi^r_{\xi_1} \rangle$  & (443.28,443.46)  & (630.03,630.27) & (0.22,0.22) &  -- & -- & -- & --\\
$\langle \Ghat_2\what, \psi^g_{\xi_1} \rangle$  & (1264.10,1264.79)  & (2459.95,2460.98) & (-0.42,-0.41) & -- & -- & -- & --\\

\hline
\end{tabular}
}
\end{center}
\caption{The 95\% confidence intervals for the coefficients of the inner products found in \eqref{eqn:ode1} and \eqref{eqn:ode2} for the system parameters in \eqref{eqn:parameters}.}
\label{table:wcoeffs3}
\end{table}
 For the POS in this paper we take the central value of these confidence intervals to denote the coefficient used in the ODE system, but we note that it is in the final four sets of inner products that these coefficients have the most variability. Hence, when we present the key results in the main paper we also present an ensemble calculation  of the most variable coefficients in order to validate the robustness of the dynamics observed. Thus, this ensemble is an ensemble of $2^{20}$ simulations using the upper and lower limits of the 20 $B_i$ and $\mu_i$ coefficients. 

\begin{figure}[htb!]
\begin{center}
(a)\includegraphics[width=0.45\textwidth,height=0.45\textwidth]{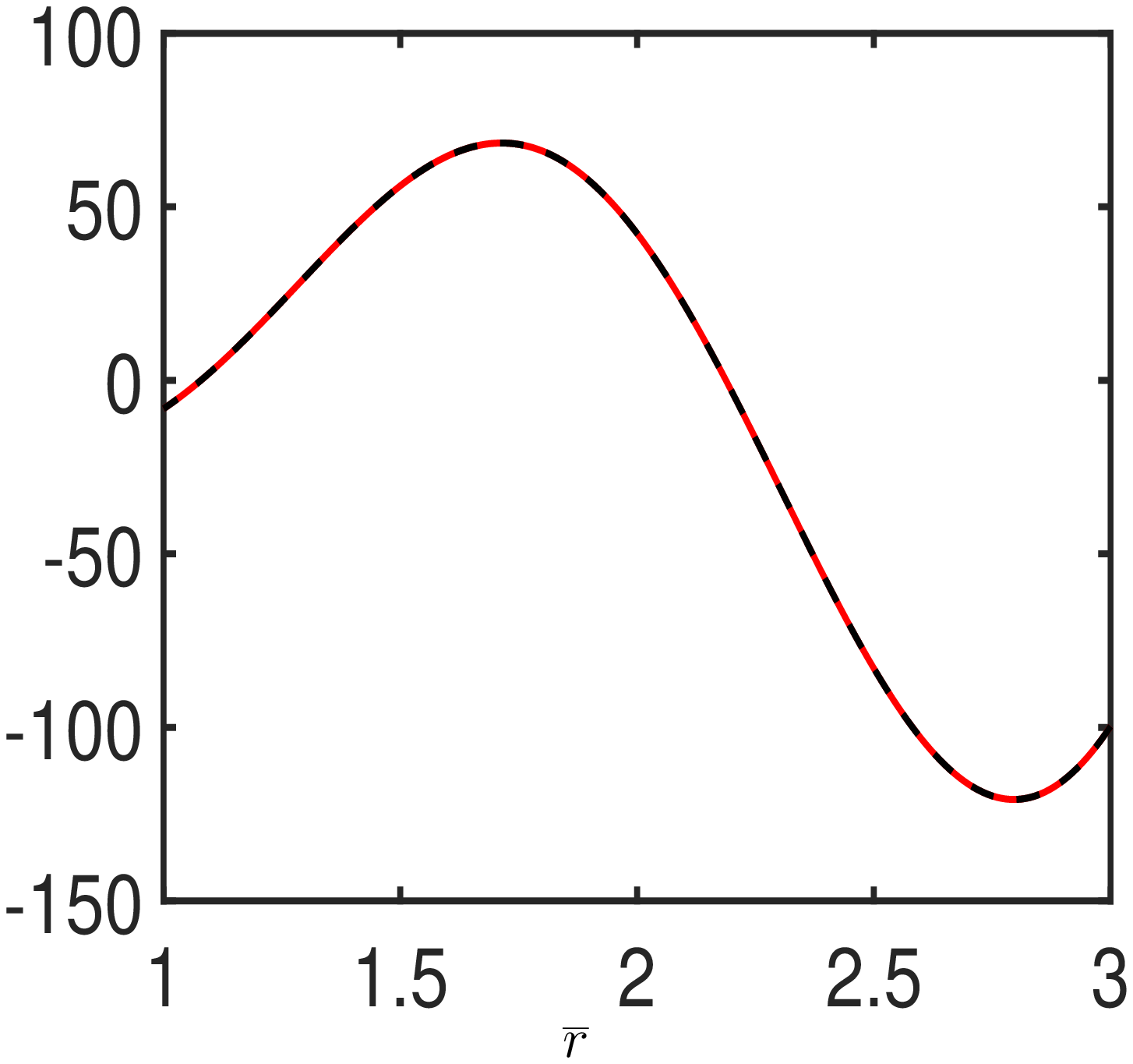}
(b)\includegraphics[width=0.45\textwidth,height=0.45\textwidth]{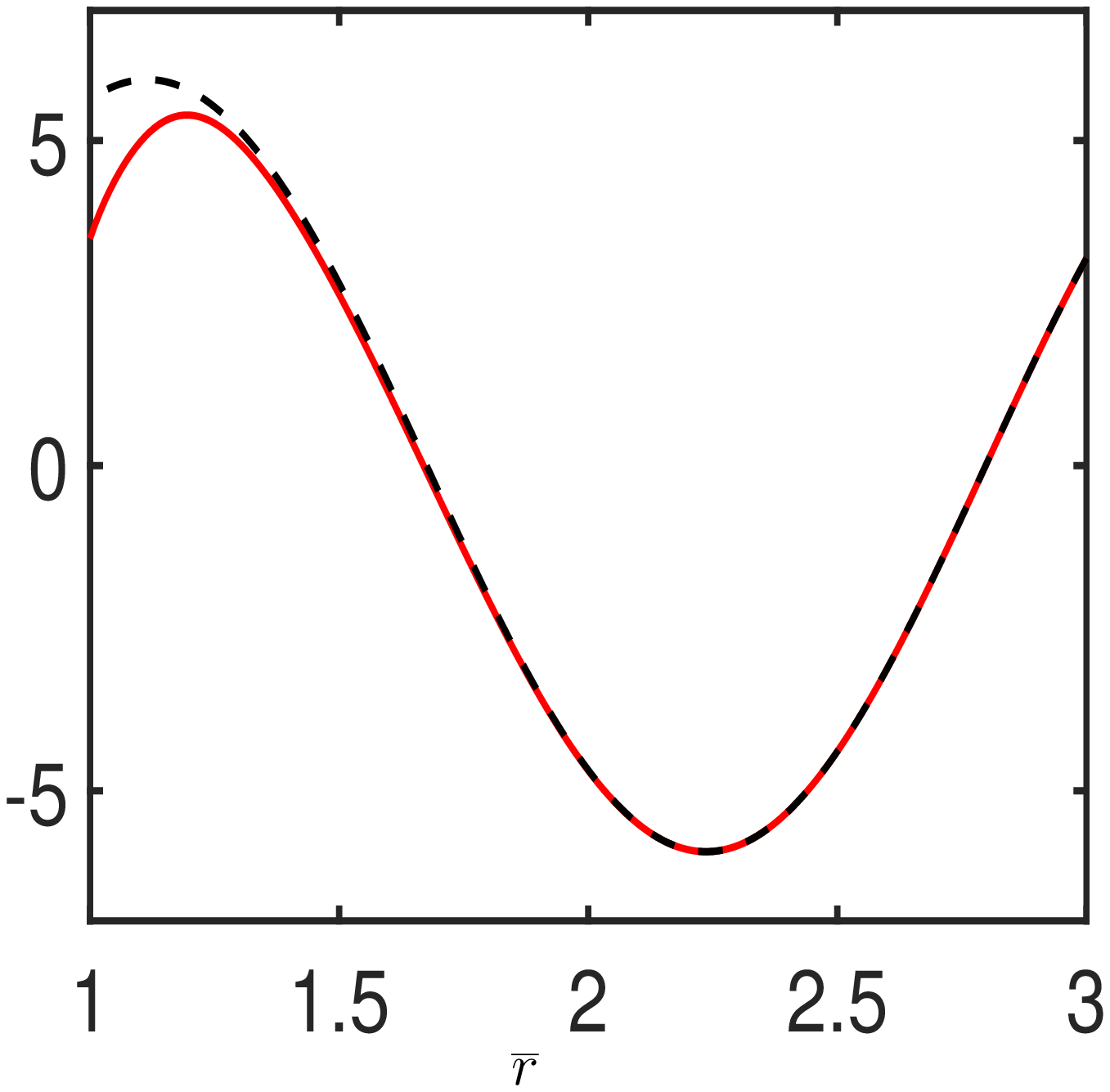}
(c)\includegraphics[width=0.45\textwidth,height=0.45\textwidth]{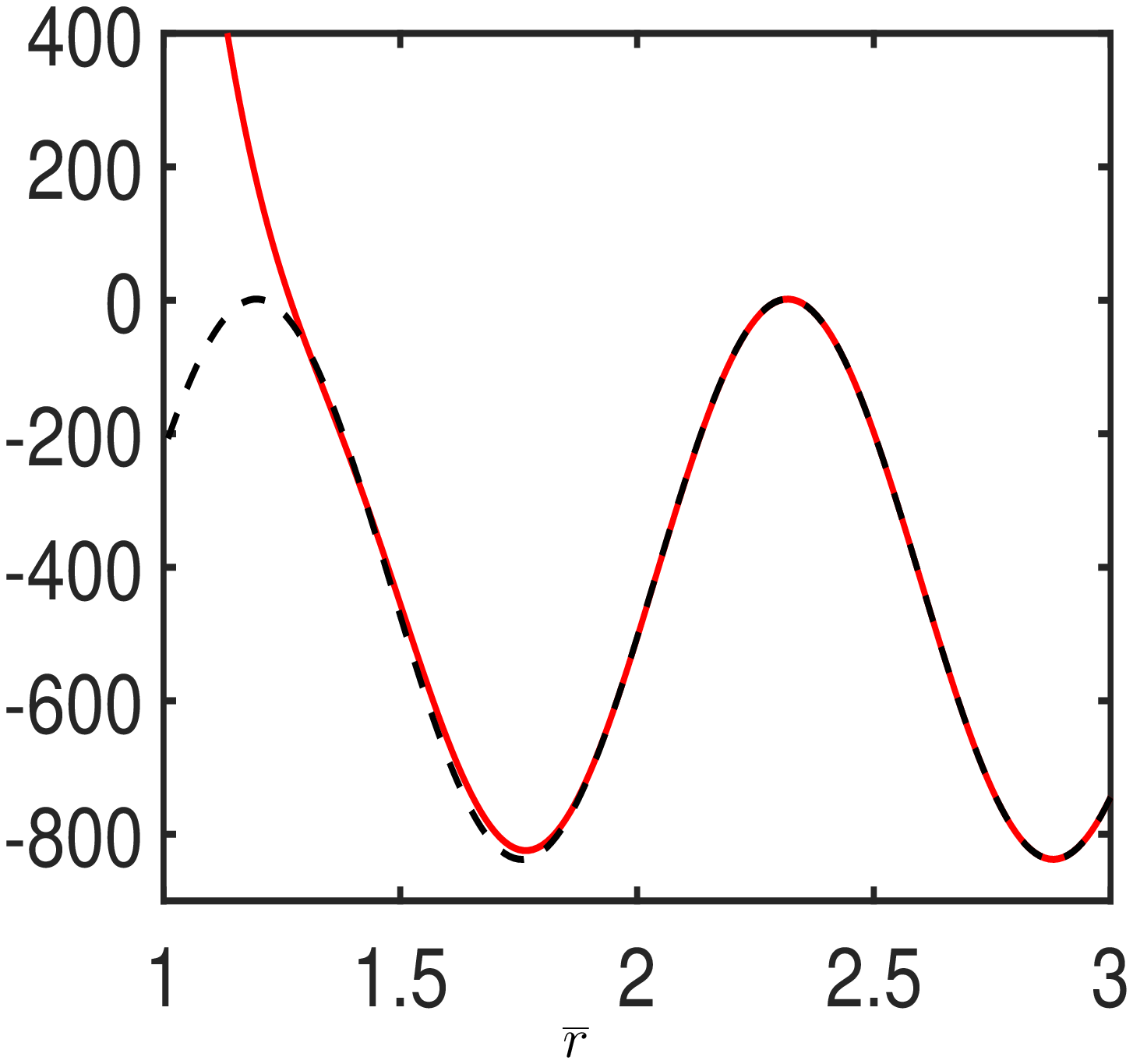}
(d)\includegraphics[width=0.45\textwidth,height=0.45\textwidth]{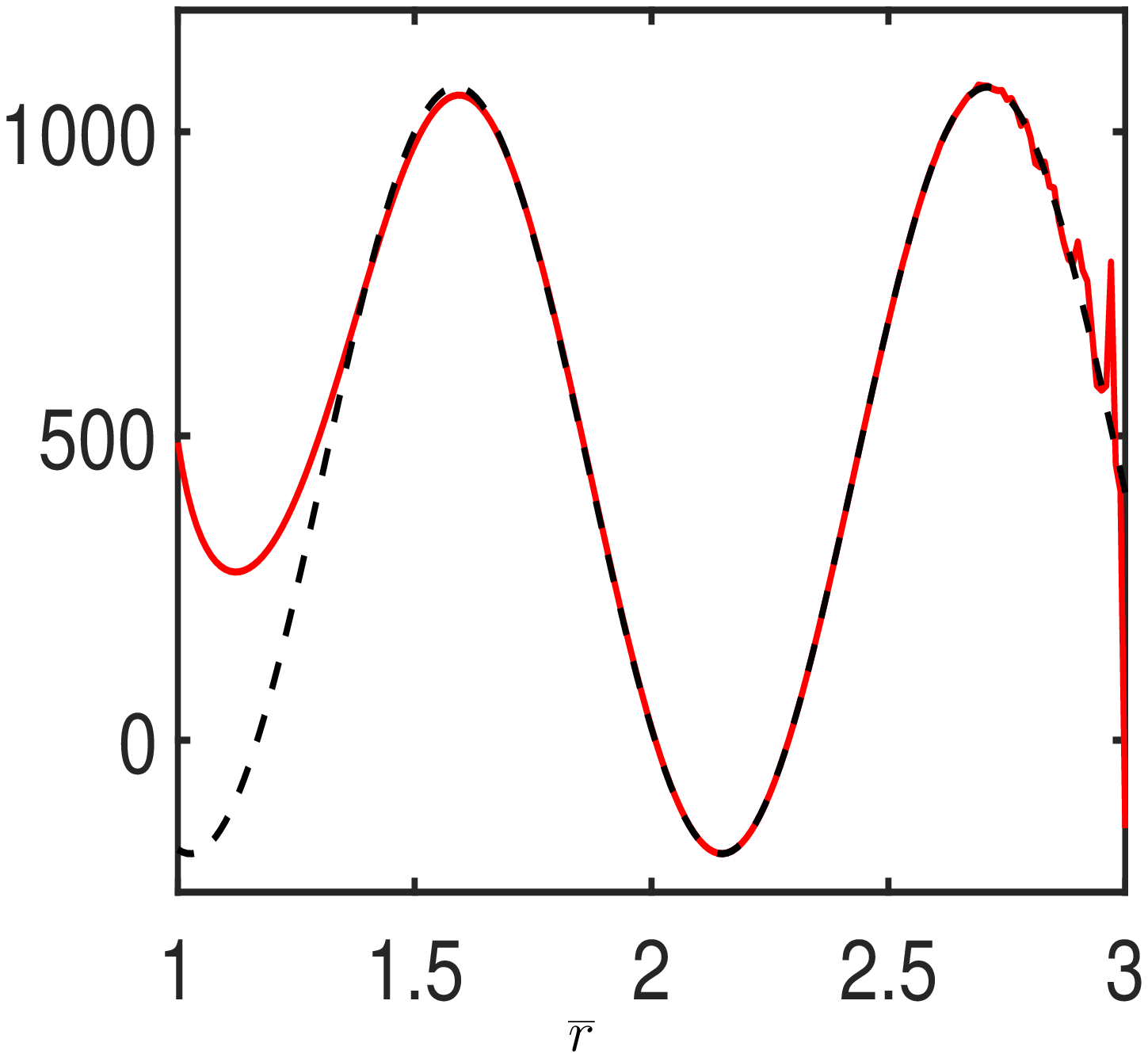}
\end{center}
\caption{(colour online) Numerical (solid line) and weakly-interacting asymptotic (dashed line) forms of the inner products (a) $e^{\lambda_r|\rbar|}\langle \varphi^r_{\xi_2}, \psi^r_{\xi_1} \rangle$, (b) $\langle \what, \partial_x\psi^r_{\xi_1} \rangle$, (c) $\langle \Psihat_1\what, \psi^r_{\xi_1} \rangle$ and (d) $\langle \Ghat, \psi^r_{\xi_1} \rangle$, with $\gbar=\sqrt{2}$.}
\label{fig:appen_approx}
\end{figure}

Plots of a selection of these inner products and their large $\rbar$ asymptotic approximations for the parameters \eqref{eqn:parameters}, taking central values of the confidence intervals, are given in figure \ref{fig:appen_approx}.
 All the inner products show excellent agreement over the region $\rbar\in[2.04,2.55]$ for which they are fitted, and beyond in most cases. One other thing to note is that for $\rbar>2.7$ the inner products involving $\Ghat$ become affected by running close to machine precision due to the fact that $e^{-2\lambda_r\rbar}=O(10^{-16})$ at this point, hence our choice of fitting region \eqref{eqn:fittingregion} is designed to stop us trying to fit coefficients in this region. The inner product in figure \ref{fig:appen_approx}(a) is the same as that in figure \ref{fig:appen_inner}, so this shows we can do an even better job of fitting these constants if more terms in the asymptotic expansion \eqref{2.tails2} are used.

Although we only present the calculation of these coefficients for the two-pulse case, it is worth noting that in the multi-pulse case, where the pulses are very-weakly-interacting, then we expect that all pulses, $V_i$ will only interact with their neighbouring pulses $V_{i-1}$ and $V_{i+1}$. Hence the inner products should have similar forms to those above, except with two terms, the first a function of $r_i-r_{i-1}$ and $g_i-g_{i-1}$ and the second a function of $r_i-r_{i+1}$ and $g_i-g_{i+1}$, similar to the interaction functions \eqref{eqn:interactionequations}.

\bibliographystyle{plain}
\bibliography{CGLpaperbib}

\end{document}